\RequirePackage{fix-cm}

\let\latexvec\vec

\documentclass[twocolumn]{svjour3}
\smartqed
\usepackage{graphicx}

\let\vec\latexvec

\usepackage{amssymb,amsmath,graphics,psfrag,graphicx,subfigure}
\usepackage[table]{xcolor}

\definecolor{lightblue}{rgb}{0.22,0.45,0.70}
\usepackage[colorlinks=true,breaklinks=true,linkcolor=lightblue,citecolor=lightblue]{hyperref}
\usepackage{stackengine}
\pdfoutput=1

\usepackage{widetext}

\numberwithin{equation}{section}
\numberwithin{figure}{section}
\numberwithin{table}{section}
\newcommand\cero{\boldsymbol{0}}
\newcommand\bB{\mathbf{B}}
\newcommand\bC{\mathbf{C}}
\newcommand\bD{\mathbf{D}}
\newcommand\bE{\mathbf{E}}
\newcommand\bF{\mathbf{F}}

\newcommand\bI{\mathbf{I}}
\newcommand\bP{\mathbf{P}}
\newcommand\bPi{\boldsymbol{\Pi}}
\newcommand\bR{\mathbb{R}}

\newcommand\bb{\boldsymbol{b}}

\newcommand\bn{\boldsymbol{n}}

\newcommand\bsigma{\boldsymbol{\sigma}}

\newcommand\btau{\boldsymbol{\tau}}

\newcommand\bu{\boldsymbol{u}}
\newcommand\bv{\boldsymbol{v}}
\newcommand\bw{\boldsymbol{w}}
\newcommand\bx{\boldsymbol{x}}

\newcommand\cG{\mathcal{G}}
\newcommand\cH{\mathcal{H}}

\newcommand\cR{\mathcal{R}}

\newcommand\cT{\mathcal{T}}

\newcommand\tr{\mathop{\mathbf{tr}}\nolimits}

\newcommand\bt{\boldsymbol{t}}

\newcommand{\mbs}[1]{\boldsymbol{#1}}

\newcommand{\jump}[1]{\left[\![#1\right]\!]}
\newcommand\fo{{\mbs{f}_0}}
\newcommand\so{{\mbs{s}_0}}
\newcommand\no{{\mbs{n}_0}}

\newcommand\mathcalbb[2][2]{%
  \stackengine{0pt}{$\mathcal{#2}$}{$\mkern#1mu\mathcal{#2}$}{O}{l}{F}{F}{L}}

\newcommand{\mypar}[1]{\medskip\noindent\textbf{#1}}

\allowdisplaybreaks
\journalname{Biomech Model Mechanobiol}

\begin{document}
\titlerunning{3D stress-assisted diffusion in cardiac electro-viscoelasticity}
\title{An orthotropic electro-viscoelastic model for the heart with stress-assisted diffusion}

\author{Adrienne Propp \and Alessio Gizzi \and \\ Francesc Levrero-Florencio \and Ricardo Ruiz-Baier}
\authorrunning{A. Propp et al.}

\institute{Adrienne Propp\at Mathematical Institute,
University of Oxford, A. Wiles Building, Woodstock Road, Oxford OX2 6GG, United Kingdom. 
\email{Adrienne.Propp@maths.ox.ac.uk}
\and 
Alessio Gizzi\at Nonlinear Physics and Mathematical Modeling, 
Department of Engineering, University Campus Bio-Medico, Rome, Italy. \email{A.Gizzi@unicampus.it}  
\and 
Francesc Levrero-Florencio \at Department of Computer Science, 
University of Oxford, 15 Parks Road, Oxford OX1 3QD, United Kingdom. \email{Francesc.Levrero-Florencio@cs.ox.ac.uk}   
\and Ricardo Ruiz-Baier (Corresponding author) \at Mathematical Institute,
University of Oxford, A. Wiles Building, Woodstock Road, Oxford OX2 6GG, United Kingdom; and 
Laboratory of Mathematical Modelling, Institute of Personalised Medicine, Sechenov University, Moscow, Russian Federation.  \email{Ricardo.RuizBaier@maths.ox.ac.uk}
}

\date{\today}
\maketitle

\begin{abstract}
We propose and analyse the properties of a new class of models for the electromechanics of cardiac tissue. The set of governing equations consists of nonlinear elasticity using a viscoelastic and orthotropic exponential constitutive law (this is so for both active stress and active strain formulations of active mechanics)  coupled with a four-variable phenomenological model for human cardiac cell electrophysiology, which produces an accurate description of the action potential. The conductivities in the model of electric propagation are modified according to stress, inducing an additional degree of nonlinearity and ani\-so\-tro\-py in the coupling mechanisms; and the activation model assumes a simplified stretch-calcium interaction generating active tension or active strain. The influence of the new terms in the electromechanical model is evaluated through a sensitivity analysis, and we provide numerical validation through a set of computational tests using a novel mixed-primal finite element scheme. 
\keywords{Orthotropic nonlinear elasticity \and mixed-primal finite element method \and Kirchhoff stress formulation \and
stress-assisted diffusion \and viscoelastic response \and cardiac electromechanics}
\subclass{\\65M60 \and 92C10 \and 74S05 \and 74F99 \and 74D10}
\end{abstract}

\section{Introduction}
In order to effectively combat cardiovascular disease, we need a robust scientific understanding of the mechanisms of the heart and the nature of such health conditions.
Recent progress in the field is encouraging; the concept of patient-specific treatment is no longer a distant dream, but a conceivable reality and a topic of ongoing research. However, a major difficulty is our incomplete knowledge about the relationship between processes at the cellular and sub-cellular level, and the performance of the organ as a whole \cite{augustin16}. Indeed, a great deal of treatment is still based on trial-and-error experimentation rather than a more fundamental scientific understanding of the changes responsible for the onset and progression of disease \cite{goktepe13}. Several treatments, such as resynchronisation therapy and anti-arrhythmic medications, for example, are known to be ineffective or even exacerbate pathological conditions in some patients, for reasons that are not yet well understood \cite{jaffe14}. 
One potential obstacle to deep understanding of cardiac function is the difficulty of acquiring sufficiently detailed data. Until recently, there were no experimental techniques capable of recording 3D cardiac activity with high enough spatio-temporal resolution to provide the required level of information. However, relatively recent studies (see e.g. \cite{christoph17}) have been able to use optical mapping to assess electromechanical waves with acceptable physiological accuracy. 

Computational models have thus been critical in allowing for extensive study of the heart even without sufficient data. The development of complex multi-scale and multi-physics models, accompanied by advances in simulation and imaging techniques, has enabled researchers to investigate the many different aspects of cardiac function and disease. The hope is that the knowledge gained from these models can contribute to new and improved treatment methods. Even though the problem of cardiac electromechanics has been the focus of a large number of modelling and computational studies (see for instance \cite{augustin16,cherubini12,colli16,costabal17,gizzi15,goktepe13,nobile12,quarteroni17,sundnes14} and the references therein), there still remain many challenges in the development of more accurate and detailed models and the accompanying methods.

In such a context, the large majority of the proposed approaches rely on continuum formulations of the complex microstructural interactions occurring among the heart tissue components, e.g. cardiomyocytes, involving different scales \cite{quarteroni17}. The study of single cell and cell-cell \cite{lenarda18} chemomechanical and electromechanical interactions has attempted to unveil some of the underlying 
complex features of the cardiac function, and different multi-field nonlinear models have been gradually generalising classical approaches as the monodomain equations and Fick's law of diffusion. In particular, fractional diffusion \cite{cusimano18}, nonlinear diffusion \cite{hurtado16}, and stress-assisted diffusion formulations \cite{cherubini17} were recently proposed to reproduce porous multiscale excitation phenomena within the framework of  homogenised models for cardiac tissue. These studies, in fact, paved the route towards new challenging theoretical and computational problems aiming at a reliable in silico prediction of heart rate variability, cardiac repolarisation and inducibility of life-threatening arrhythmias \cite{phadumdeo18}.
At the same time, macro-scale incompressibility, orthotropic and hysteretic mechanical features have been shown to fully characterise the human cardiac tissue under multiaxial loading tests (see e.g. \cite{gultekin16} and references therein). Viscosity properties, in particular, have been incorporated as one spring element coupled with Maxwell elements in parallel endowing the model with hysteretic characters describing the viscous response due to matrix, fibre, sheet and fibre-sheet couplings through four dedicated dashpots \cite{gultekin16}. Also in this case, a porous medium motivation has been advanced in \cite{yao12}, including the extracellular fluid filtrating through the elastic body, contributed by the active contractile behaviour of the muscle. However, complete agreement concerning the specific multiscale features involved in energy dissipation for the cardiac tissue, and soft biological tissues in general, is still lacking. The stress evolution equations for time-dependent viscous behaviour are based on finite strain viscoelasticity \cite{holzapfel01}, motivated by a rheological analogue from \cite{simo87}, and endowed with equilibrium and non-equilibrium contributions \cite{lubliner85} in which the usual assumption of volume-preserving deformations during time-dependent responses is made.

A distinguishing feature of our approach is the introduction of the mechanoelectrical feedback (MEF) in the electric conductivities, through a direct dependence on the Kirchhoff stress. This framework, known as stress-assisted diffusion (SAD), is widely employed in the modelling of gels and polymers \cite{klepach14}, but has only recently been adapted for active biological media undergoing reaction-diffusion excitation \cite{cherubini17}, and more tailored for cardiac models in \cite{loppini18}. While these contributions consider hyperelastic formulations coupled with multiphysics activation mechanisms, we also consider here the viscoelastic effects typical of soft microstructured fibre-reinforced biological tissues, and using realistic ventricular geometries.  
Fully mixed methods for the hyperelasticity of the cardiac tissue (that is, formulations involving stresses or strains in addition to simply displacement or displacement-pressure) are not yet widely employed. They have been introduced in \cite{ruiz15} and used more recently in  \cite{cherubini17,garcia19,ruiz18}. In our case, our model and our numerical method include a three-field elasticity formulation (variationally based on a modification of the Hu-Washizu principle \cite{lami06}) that states the governing equations in terms of stress-displacement-pressure, and that is motivated by the desire to avert volumetric locking and to solve directly for additional variables of interest. In particular, we solve for the Kirchhoff stress, which we use explicitly in our incorporation of SAD. This formulation 
includes a pressure-stabilisation term needed, in the lowest-order case, for triangular or for tetrahedral meshes. It constitutes a generalisation of the three-field formulation for nearly incompressible hyperelasticity, designed in \cite{chavan07} using quadrilateral meshes. Another difference in the {present} contribution is that we employ a more accurate 
cellular model, tailored for recovering human action potential dynamics, restitution features under constant pacing as well as sustained fibrillation behaviours and spiral waves breakup \cite{bueno08}. 
While the active strain approach is adopted in many instances in the literature and is often favoured due to the practicality of measuring strains directly using imaging techniques \cite{rossi14}, the active stress approach is somewhat simpler and more naturally incorporated in already existing models for passive deformation \cite{giantesio18}. In this work, we will adopt both  formulations, although we find that the active strain formulation better reproduces physiologically accurate deformation regimes in ventricular geometries. 
To the best of our knowledge, no previous attempts have been made incorporating both active stress and active strain within a generalised stress-assisted reaction-diffusion formalism and embedding orthotropy, incompressibility and viscoelasticity for human cardiac ventricular 
domains.

This paper has been structured in the following manner. 
Section \ref{sec:model} lays out the elements of the mathematical model 
that describes the electro-viscoelastic function of the heart, including 
 the active contraction of the cardiac muscle and the 
representation of the mechanoelectric feedback using 
stress-assisted conductivity, {as well as} a contribution from geometric nonlinearities (or geometrical 
feedback). 
The passive hyperelastic response of the tissue is 
described by an orthotropic exponential model, whereas 
the ionic activity which causes active contraction is  
incorporated through orthotropic active stress (active strain will {also} 
be addressed).
The specific structure of the governing equations (written  
in terms of stress, displacements, electric potential, activation 
generation, and ionic variables) suggests to cast the problem in a 
mixed-primal form, and to use a mixed-primal finite element method 
for its numerical approximation. This is precisely the method that 
we outline in Section~\ref{sec:FE}, which also includes a description of the 
consistent linearisation and implementation 
details. Our computational results in 2D and 3D, along with 
numerical validation and 
pertinent discussions on the modelling considerations are then 
presented in Section~\ref{sec:results}. We close with a summary and some 
remarks on model limitations and ongoing extensions, collected in Section~\ref{sec:concl}.

\section{Mathematical model}\label{sec:model}

\subsection{Finite-strain cardiac mechanics} Let $\Omega\subset\bR^d$, $d\in\{2,3\}$ denote a deformable body with
a piecewise smooth boundary $\partial\Omega$, considered in its reference
configuration, and {let $\bn$ denote} the outward unit normal vector on
$\partial\Omega$. The kinematical description of finite deformations
regarded on a time interval $t\in (0,t_{\text{final}}]$ is made
precise as follows.  A material point in $\Omega$ is denoted by $\bx$,
whereas $\bu(t):\Omega\to\bR^d$ will denote the displacement field
defining its new position in the deformed configuration.  
The tensor $\bF:=\bI + \nabla\bu$ is the
gradient (applied with respect to the fixed material coordinates) of
the deformation map; its determinant, denoted by $J=\det\bF$,
measures the solid volume change during the deformation; and
$\bC=\bF^{\tt t}\bF$ and $\bB=\bF\bF^{\tt t}$ are respectively the right and left Cauchy-Green deformation tensors on
which all strain measures will be based (here the superscript $()^{\tt t}$ denotes 
the transpose operator). The first isotropic
invariant ruling deviatoric effects is the scalar $I_{1}(\bC)=\rm tr \, \bC$, and for
generic unitary vectors $\fo, \so$, the scalars
$I_{4,f}(\bC)=\fo\cdot(\bC\fo)$, $I_{4,s}(\bC)=\so\cdot(\bC\so)$, 
$I_{8,fs} (\bC)=\fo\cdot(\bC\so)$ are
 pseudo-invariants of $\bC$ measuring direction-specific 
stretch \cite{ciarlet}.

The triplet $(\fo(\bx),\so(\bx),\no(\bx))$ represents a coordinate system pointing
in the local direction of the muscular cardiac fibres, transversal sheetlet compound,
and normal cross-fibre direction $\no(\bx) = \fo(\bx)\times\so(\bx)$.
Note that the system is restricted to $(\fo(\bx),\so(\bx))$ in the
two-dimensional case, and that these directions are defined in the reference configuration. 
Constitutive relations characterising the material
properties and underlying microstructure of the myocardial tissue will follow the orthotropic model proposed in
\cite{holzapfel09}, whose strain energy density 
(relating the amount of energy stored within the material in response [Joule/Volume]
to strain, and which assumes an additive decomposition into 
isotropic, volumetric and anisotropic contributions) and the
first Piola-Kirchhoff stress tensor (associated with a passive, elastic deformation) read, respectively
\begin{align}\nonumber
\Psi_{\mathrm{pas}}(\bF) & = \frac{a}{2b}e^{b(I_{1}-d)}
+
\sum_{i\in\{ f,s\} } \dfrac{a_{i}}{2b_{i}}\bigl[e^{b_{i}(I_{4,i}-1)_+^{2}}-1\bigr]\\ 
& \qquad +\dfrac{a_{fs}}{2b_{fs}}\bigl[e^{b_{fs} (I_{8,fs})^2 }  -1\bigr], \label{piolatensor}
\\ 
\bP_{\mathrm{pas}} & =\frac{\partial\Psi_{\mathrm{pas}}}{\partial\bF}-pJ\bF^{-{\tt t}},\nonumber
\end{align}
where $a,b$ are material constants associated with the isotropic matrix response, 
$a_f$ and $b_f$ rule the directional behaviour of the material along myocardial fibres, 
$a_s$ and $b_s$ account for the cross-contribution of the fibre sheet directions, and 
$a_{fs},b_{fs}$ encapsulate the shear effects in the fibre-sheet plane. 
Moreover, the field $p$ denotes the solid hydrostatic pressure, and we have used the notation 
$(u)_+:= u$ if $u>0$ or  zero otherwise, for a generic real-valued function $u$.   
This modelling choice aligns with the fact that fibres have a quite different behaviour 
under compression or tension regimes. In addition, 
taking the positive part of the exponents in the anisotropic energy   
results in excluding anisotropic energetic contributions for compressed fibre 
configurations, which in the case of passive fibres should have an effect only during 
extension \cite{pezzuto14}. We remark here that the particular mechanisms of
 soft tissue anisotropic mechanical behaviour is still under investigation \cite{humphrey14}. Moreover, 
full incompressibility of the tissue will be enforced in the present framework, and this 
has some advantages associated with the mathematical and numerical structure of the system. Although biological tissues 
possess a complex porous structure, compression features are still being systematically investigated 
ex-vivo, and a more comprehensive answer on the subject is still needed \cite{mcevoy18}.

\subsection{Active stress and active strain} 
In physiological scenarios, 
the mechanical deformation is also actively influenced by microscopic tension generation. 

\mypar{Active stress model.}
A simple description is given in terms of active stresses (see for instance \cite{sundnes14}): 
we assume that the first Piola-Kirchhoff 
stress tensor decomposes as 
\begin{equation}\label{eq:total-stress}
\bP = \bP_{\mathrm{pas}} + \bP_{\mathrm{act}},
\end{equation}
where  
the active stress component acts differently on each 
local direction with an intensity depending on the scalar field of active tension $T_a$, 
that synthesises the biochemical state of myocytes (and whose 
dynamic behaviour will be specified later on). Then 
\begin{align}
\nonumber	\bP_{\text{act}} & = 
	J\bsigma_{\text{act}} \bF^{-{\tt t}}, \quad \text{with} \\    
	\bsigma_{\text{act}} &=  \frac{T_a}{J\lambda_f} \bF\fo\otimes\bF\fo +
 	\frac{\kappa_{sn} T_a}{J\lambda_s\lambda_n} \mathrm{sym}(\bF\so\otimes\bF\no)\nonumber  \\    
 	&\quad +\frac{\kappa_{nn}T_a}{J\lambda_n} \bF\no\otimes\bF\no,\label{eq:active-stress}
\end{align}
where $\kappa_{sn},\kappa_{nn}$ are positive constants representing the variation of 
activation on each specific 
direction, as proposed in \cite{dorri06}, and $\lambda_f=\sqrt{I_{4,f}},\lambda_s=\sqrt{I_{4,s}},\lambda_n=\sqrt{\no\cdot(\bC\no)}$ are the fibre, sheetlet, and cross - fibre 
stretches.
Setting appropriate models for $ \bsigma_{\text{act}}$ is not a trivial task since the active contribution 
to the force should account for the geometric properties of deformation, and these undergo substantial 
changes during contraction in the finite strain regime \cite{pezzuto14}.  
Details of other anisotropic activation forms can be found, for instance, in \cite{rossi14} 
for active strain and in \cite[Appendix B]{usyk00} for active stress descriptions, but they 
are basically responsible for additional deformation effects such as wall thickening, radial constriction and 
torsion, as well as longitudinal shortening. Note that the active Cauchy stress does not include a contribution 
on the diagonal entry associated with the local sheetlet direction $\so$ since a stress component on this direction would 
counteract wall thickening mechanisms \cite{dorri06}. Moreover, the intensity of the active tension effect 
on the cross-fibre direction $\no$ is assumed to be substantially 
smaller than that appearing on the off-diagonal component $\mathrm{sym}(\bF\so\otimes\bF\no)$, see Table~\ref{table:params}. 
Also note that some references do not include a rescaling with local stretches in each term of $ \bsigma_{\text{act}}$. 

\mypar{Active strain model.}
Next we recall the active strain model for ventricular
electromechanics (see e.g. \cite{cherubini08}). There, the
contraction of the tissue results from activation mechanisms governed
by internal variables and incorporated into the finite elasticity
context using a multiplicative decomposition 
of the deformation gradient into
a passive (purely elastic) and an active part, $\bF =\bF_{\!E}\bF_{\!A}$, with
\begin{align}\nonumber
\bF_{\!A} &=\bI+\gamma_f\fo(\bx)\otimes\fo(\bx)  \\ &   \label{eq:Fa}
\qquad   +\gamma_s\so(\bx)
\otimes\so(\bx)+\gamma_n\no(\bx)\otimes\no(\bx).
\end{align}
The 
coefficients $\gamma_i$, with $i=f,s,n$, are smooth scalar functions encoding the
macroscopic stretch in specific directions, whose precise definition
will be postponed to \eqref{eq:gammas}.  The inelastic contribution to the deformation 
modifies the length and maybe also the shape of the cardiac fibres, and then compatibility 
of the motion is restored through an elastic deformation accommodating 
the active strain distortion. A physiological motivation for the active strain approach is related 
to the shortening of sarcomeres as a response to the sliding filaments of the actin-myosin 
molecular motor: such shortening is encapsulated in $\bF_{\!A}$, which determines a new 
(and fictitious, or virtual) intermediate configuration that is regarded as a reference for the 
elastic deformation \cite{pezzuto14}. Therefore, the strain energy function and the
first Piola-Kirchhoff stress tensor (after applying the active strain
decomposition) are functions of $\bF_{\!E}$ only, and they read respectively 
\begin{align}
\nonumber
\widehat\Psi(\bF_E) & = \frac{a}{2b}e^{b(I^E_{1}-d)} +
\sum_{i\in\{ f,s\} } \dfrac{a_{i}}{2b_{i}}\bigl[e^{b_{i}(I^E_{4,i}-1)_+^{2}}-1\bigr]\label{piola:astrain}\\
&\quad  +\dfrac{a_{fs}}{2b_{fs}}\bigl[e^{b_{fs}(I_{8,fs}^E)^{2}}-1\bigr] , 
\\ 
\bP& =\frac{\partial\widehat\Psi}{\partial\bF}-pJ\bF^{-{\tt t}}.\nonumber
\end{align}
As in the description of \eqref{piolatensor} above, we again note that one switches off the anisotropic 
 contributions under compression. An additional advantage is that the associated terms in the 
 strain energy function (in both the pure passive and active-strain formulations) 
 can be shown to be strongly elliptic \cite{pezzuto14} (these will be the terms 
appearing on the second diagonal block of the weak formulation from Section~\ref{sec:FE}, 
the block corresponding to displacements), however 
the overall problem will remain of a saddle-point structure. 
The modified elastic invariants $I_i^E$ are functions of
the coefficients $\gamma_i$, as well as of the invariant and pseudo-invariants in the 
following manner \cite{rossi12}
%
{
\begin{align*}
I_{1}^{E} & = 
I_1 - \gamma_{f}\frac{\gamma_{f}+2}{(1+\gamma_{f})^2} I_{4,f}
- \gamma_{s}\frac{\gamma_{s}+2}{(1+\gamma_{s})^2} I_{4,s} \\&
- \gamma_{n}\frac{\gamma_{n}+2}{(1+\gamma_{n})^2} I_{4,n}  \,,
\\
I_{4,f}^E & = \left(1+\gamma_{f}\right)^{-2} I_{4,f} ,\qquad 
I_{4,s}^E  = \left(1+\gamma_{f}\right) I_{4,s},\\ 
I_{8,fs}^E & = \left(1+\gamma_{f}\right)^{-1/2} I_{8,fs}\,.
\end{align*}
Such dependencies are a consequence of assuming isochoric active 
deformations \cite{pezzuto14}, i.e. $\det \bF_{\!A}=1$, justified by the fact that 
the volume of the cardiomyocytes does not vary substantially during contraction.
Besides, following Rossi et al. \cite{rossi12}, previous expressions are obtained 
by assuming $\gamma_s=\gamma_n$ and
making use of the fact that $I_1=I_{4,f}+I_{4,s},I_{4,n}$, as well as that 
$\bF_{\!E}=\bF\bF_{\!A}^{-1}$, with
\begin{align*}
\det \bF_{\!A} &= (1+\gamma_f)(1+\gamma_s)(1+\gamma_n) \,, \\
\bF_{\!A}^{-1} &= \bI -
\dfrac{\gamma_f}{1+\gamma_f} \fo\otimes\fo \\&- 
\dfrac{\gamma_s}{1+\gamma_s} \so\otimes\so -
\dfrac{\gamma_n}{1+\gamma_n} \no\otimes\no \,.
\end{align*}}
Accordingly, the active strain, and consequently the force associated to the 
active part of the total stress, will receive contributions acting 
on the three main directions. The calcium-based activation signal travels up 
to four times faster along the fibre axis than in the sheet and normal directions, 
and this fact further motivates the use of orthotropic active strain \cite{rossi14}.

\subsection{Viscoelasticity and equations of motion}  
Extension and shear tests demonstrate the importance of incorporating viscoelastic 
effects in models for cardiac passive mechanics \cite{gultekin16}. 
In the heart, the extracellular fluid filtrating through the elastic solid is one of the main 
generators of the viscoelastic effects of the tissue \cite{yao12}. {Viscous effects are also tied to crossbridge processes identified in Ca$^{2+}$ activated fibres \cite{maughan98}}, and have a well-established literature 
as well as a consistent methodology for their implementation (the stress update algorithm that uses 
a convolution integral representation) developed for general soft tissues \cite{holzapfel01}.
From the viewpoint of kinematics, it suffices 
to relate stress to strain rates. Decomposition of the spatial velocity gradient $\bw = \dot{\bu}$ into 
the rate of deformation and spin tensors yields the relation 
$$\dot{\bB}  = \nabla\bw \bB + \bB (\nabla \bw)^{\tt t},$$
and a simplified rheological Kelvin-Voigt model for the viscous component of the Cauchy stress can be defined 
as follows (see e.g. \cite{karlsen17}) 
\begin{equation}\label{eq:beta}
\bsigma_{\text{visc}} = \delta e^{\beta \dot{I_1}} \dot{\bB},
\end{equation}
which depends on the history of the isotropic contribution to the Cauchy stress. Here 
$\delta,\beta >0$ are model parameters. In this way, 
after a pull-back operation, we see that 
\begin{equation}
\bP_{\text{tot}}  = \bP + J\bsigma_{\text{visc}} \bF^{-\tt t},\label{eq:viscoelastic}
\end{equation}
is the total first Piola-Kirchhoff stress tensor that includes $\bP$ defined from either 
 \eqref{piolatensor}-\eqref{eq:total-stress} or \eqref{piola:astrain}, and the viscoelastic 
contributions. 

More advanced rheologies can be easily incorporated in the context of active stress formulations 
as done in e.g. \cite{katsnelson04}, as the generalised Hill-Maxwell model recently 
proposed in \cite{cansiz17}, {as in the perturbed equations of harmonic wave motion using springpot-based models 
with fractional order derivatives from \cite{capilnasiu19}}, 
or as the thermodynamical electro-viscoelastic models 
that use statistical fibre distributions \cite{gizzi18b}.
We will, however, confine the presentation to \eqref{eq:viscoelastic} without introducing 
stochasticity of the anisotropic components. 

Irrespective of the activation formalism one adopts (active strain or active stress), the balance of linear momentum and the 
incompressibility constraint (allowing only isochoric motions) are 
written together in the following way, when posed in the inertial reference frame and under 
transient mechanical equilibrium, 
\begin{subequations}\label{eq:mech}
\begin{align}
\rho\partial_{tt}\bu- \boldsymbol{\nabla}\cdot \bP_{\text{tot}}  & = \rho_0 \bb & \text{in } \Omega\times(0,t_{\mathrm{final}}], \label{eq:momentum}\\ 
\rho J - \rho_0 & = 0 & \text{in } \Omega\times(0,t_{\mathrm{final}}], \label{eq:mass}
\end{align}
\end{subequations}
where  
$\rho_0, \rho$ are the reference and current 
medium density, $\bb$ is a smooth vector field of imposed body loads, 
$\partial_{tt}$ denotes the second time derivative, and 
the divergence operator in \eqref{eq:momentum} applies on 
the tensor fields row-by-row. The 
balance of angular momentum translates into the condition that the 
Kirchhoff stress tensor $\bPi= \bP_{\text{tot}}\bF^{\tt t}$ must be symmetric, which is in turn 
encapsulated into the momentum and constitutive relations \eqref{eq:momentum}, \eqref{piolatensor}, \eqref{piola:astrain}.

Defining 
$$\cG = \begin{cases}
\displaystyle \cG(\bu,T_a) & := \frac{\partial\Psi}{\partial\bF}\bF^{\tt t} + J\bsigma_{\text{visc}} + \bP_{\text{act}}\bF^{\tt t}\\
&\qquad \qquad  \qquad  \text{for active stress},\\[2ex]
\displaystyle \cG(\bu,\gamma)& := \frac{\partial\widehat\Psi}{\partial\bF}\bF^{\tt t} + J\bsigma_{\text{visc}}  \\
&\qquad \qquad  \qquad  \text{for active strain},
\end{cases}
$$ 
as the contribution to the Kirchhoff stress that does not involve pressure, we then have 
\begin{equation}\label{eq:Pi}
\bPi = \cG - p J\bI.
\end{equation}

\begin{figure*}
\begin{center}
\includegraphics[height=0.42\textwidth]{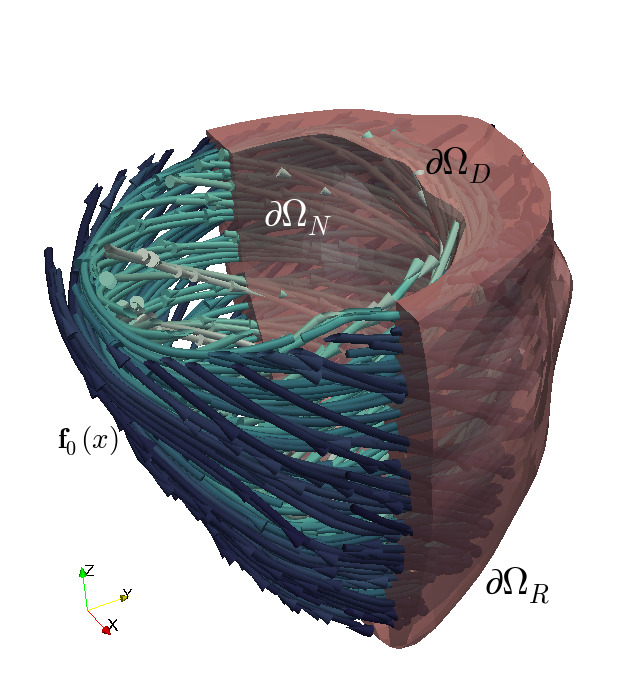}
\raisebox{3mm}{\includegraphics[height=0.4\textwidth]{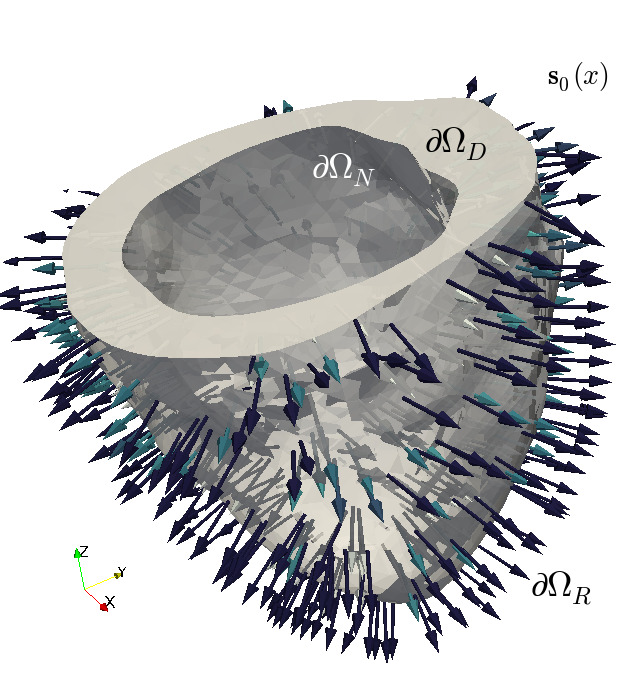}}
\end{center}
\vspace{-3mm}
\caption{Schematic representation of a mono-ventricular domain where 
\eqref{eq:dirichletBC} is imposed on the basal cut, \eqref{eq:pressureBC} on 
the endocardial surface, and \eqref{eq:robinBC} on the epicardium. The left panel 
depicts the fibre field and the right panel the sheetlet directions (in this case, parallel to the normal direction of the 
epicardium).} \label{fig:sketch}
\end{figure*}

Stating the balance equations in terms of Kirchhoff stress, 
displacements, and pressure suggests that, at the level of 
writing finite element schemes, we will use mixed methods 
respecting the same structure. 
Additionally, setting boundary 
conditions for the motion of the left ventricle is not trivial, as the 
organ is known to slightly move and twist during the heartbeat. 
In our case, equations \eqref{eq:momentum}-\eqref{eq:mass}-\eqref{eq:Pi} will be supplemented 
 with mixed normal displacement and traction boundary conditions 
\begin{subequations}
\begin{align}
\label{eq:dirichletBC} \bu\cdot \bn &= 0 & \text{on}\ \partial\Omega_D\times(0,t_{\mathrm{final}}], \\
\label{eq:pressureBC} \bPi\bF^{-\tt t} \bn-  p_N J \bF^{-\tt t}\bn&= \cero   & \text{on}\ \partial\Omega_N\times(0,t_{\mathrm{final}}],\\
\label{eq:robinBC} \bPi\bF^{-\tt t} \bn + \eta J \bF^{-\tt t} \bu & = \cero  &\text{on}\ \partial\Omega_R \times(0,t_{\mathrm{final}}],
\end{align}
\end{subequations}
where $\partial\Omega_D$, $\partial\Omega_N$, $\partial\Omega_R$ conform a disjoint partition of the boundary. 
{The condition \eqref{eq:dirichletBC} means that 
we constrain the normal motion along the normal direction with respect to the surface  $\partial\Omega_D$.
 The  term $p_N$ in \eqref{eq:pressureBC} denotes} a (possibly time dependent) prescribed boundary pressure associated with 
endocardial blood pressure, which is uniform over the deformed 
counterpart of $\partial\Omega_N$ and it is applied in the normal direction to the epicardium in the 
deformed configuration. However, owing to Nanson's formula \cite{ciarlet}, this contribution 
regarded on the reference configuration depends on the cofactor of the deformation gradient and therefore the 
boundary condition is nonlinear in the undeformed configuration; moreover, the traction written in terms 
of the Kirchhoff stress tensor is $\bt = \bPi\bF^{-\tt t}\bn$. 
 Also note that the 
Robin 
conditions \eqref{eq:robinBC} account for stiff springs 
connecting the cardiac medium with the surrounding soft tissue 
and organs (whose stiffness is encoded in the scalar $\eta$). More sophisticated boundary conditions that 
consider an interaction with the pericardium can be also imposed \cite{fritz14}. A sketch of a mono-ventricular 
domain specifying boundary surfaces and fibre directions is depicted in Figure~\ref{fig:sketch}. 

\subsection{Monodomain equations} In the context of 
electromechanical processes, the propagation of electric 
potential $v$  is governed by the following 
reaction-diffusion system, known as the 
monodomain equations (see e.g. \cite{ColliBook}), which are cast here in the reference configuration.  The current 
conservation is written only in terms of the transmembrane potential and the coupling with additional ionic quantities are encoded in the vector 
$\vec{r}$ {(here we use $\vec{\cdot}$ instead of bold to denote vector fields of dimension other than $d$)}  
\begin{subequations} \label{eq:monodomain}
 \begin{align}
 \chi {\partial_t v} -  
	\nabla \cdot \{ \bD(v,\bF,\bPi)\, \nabla v \} 
&=  g(v,\vec{r}) + I_{\rm ext}\nonumber \\
& \quad\text{in } \Omega\times(0,t_{\mathrm{final}}],\label{eq:monodomain1}\\
	\frac{d\vec{r}}{dt} & = \vec{m}(v,\vec{r}) \nonumber \\
&\quad \text{in } \Omega\times(0,t_{\mathrm{final}}].\label{eq:monodomain2}
\end{align}
\end{subequations}
Here $\chi$ is the
ratio of membrane area per tissue volume, and 
$I_{\rm ext}$ is a spatio-temporal external stimulus applied 
to the medium. We will adopt the minimal model for human ventricular action potential, 
proposed in 
\cite{bueno08} and fitted to capture 
restitution curves, conduction velocity, spiral/arrhythmic dynamics, and 
complex behaviour typical to nonlinear dynamical systems (used later for cardiac 
alternans in \cite{gizzi13}). That model was, however, tailored originally for  the  case of isotropic conductivity $\bD=D\bI$, 
and 
so the extended fully-coupled model discussed below will be able to 
accommodate a wider class of propagation 
patterns, and will also constitute a generalisation over other 
recent models for stress-assisted diffusion \cite{cherubini17,loppini18}. 

{Specification of the ionic currents and gating variables can be found in Appendix~\ref{sec:app-minimal}.}

Boundary and initial conditions for \eqref{eq:monodomain} correspond to 
\begin{subequations}
\begin{align}
\bD(v,\bF,\bPi)\, \nabla v \cdot \bn & = 0 \qquad \text{on $\partial\Omega \times (0,t_{\text{final}}]$},\label{bc-mono}\\
v(0) = 0, \quad \vec{r} & = [1,1,0] \qquad \text{in $\Omega \times \{0\}$},\label{init-mono}
\end{align}\end{subequations}
and \eqref{init-mono} can be combined with suitable initial pacing, especially needed in more complex and more physiologically accurate cell models. 
The minimal model, as proposed in \cite{bueno08}, has a heterogeneous character that we do not consider in our study. Their 
description contains separate parameter sets that are able to reproduce experimental results for the epicardium, mid-myocardium and endocardium, as well as parameter sets that mimic the results of two more complicated ionic models for human ventricular cells. For simplicity (and also as a consequence of lack of personalised experimental data) we use the parameter set developed for the epicardium {(see values in Table~\ref{table:params})}, assuming that it is consistent throughout the cardiac wall. Extension to the heterogeneous case can be readily incorporated.  

\subsection{Stress-assisted conduction} 
The mechanoelectrical feedback (the process where the current 
mechanical state of the deforming solid modifies both the excitability and electrical conduction of the tissue) is here 
introduced in the conductivity tensor, through a direct 
dependence on the Kirchhoff stress (which constitutes 
one of the main novelties in our approach, stemming as a generalisation of 
the anisotropy induced by stress proposed in \cite{cherubini17} 
and later used for simplified 2D cardiac electromechanics in \cite{loppini18}). 
In addition, due to the Piola transformation (yielding a transformation  
of the diffusion tensor using the deformation gradients), the
conductivity tensor also depends 
nonlinearly on the deformation gradient (actually,  the term $J\bC^{-1}$ 
constitutes a strain-based modification of tissue
conductivity, also referred to as geometric feedback in 
\cite{colli16}) 
\begin{align}\nonumber
\bD(v,\bF,\bPi) = & [D_0 + D_1 v]J\bC^{-1} + D_0/2J \fo\otimes\fo \\
&\qquad + D_2 J\bF^{-1} \bPi \bF^{\tt -t},\label{def:D}
\end{align}
where the nonlinear conductivity (self diffusion depending on $v$) accounts 
for porous media electrophysiology following the development in \cite{hurtado16}, 
but appropriately modified to incorporate information 
about preferred directions of diffusivity according to the microstructure of the tissue (encoded in the 
second term defining $\bD$). 
The parameter $D_0$ signifies the usual diffusion for isotropic materials{,} whereas 
$D_1$ and $D_2$ represent the intensity of the porous media electrophysiology 
and of the stress-assisted diffusion, respectively. An additional term in the 
nonlinear self-diffusion (e.g. $D_3v^2$, as in \cite{gizzi17,ruiz18}) eventually leads to very slight 
modifications in conduction patterns and we have therefore decided not to include it. 
Tuning $D_1$ is sufficient to, if needed, calibrate the speed and action potential duration 
at the depolarisation plateau phase. 

It is useful to point out that both the nonlinear self-diffusion term and the SAD argument derive from rigorous thermodynamical principles, formulated under specific assumptions for porous materials. In particular, nonlinear self-diffusion is naturally related to the 
 transport of chemicals within porous media, while classical models of stress-assisted diffusion for general materials \cite{aifantis80} also consider the transport of diffusing species within solids exhibiting finite strains. 
For the specific case of cardiac tissue, both approaches are justified by the multiple scales involved in the transport of ions and generation and propagation of action potential within the cell and across different cells \cite{lenarda18}. In particular, we can mention the role of intercalated discs and gap junctions between communicating cells or the presence of the ephathic couplings in the extracellular space \cite{ly18}, as well as micro-invaginations on the cell membrane known as microtubules and microdomains \cite{miragoli16}. All of these emerging effects contribute to the macroscopic nonlinearities {and additional anisotropies} considered in the diffusion tensor herein and which could be further analysed through a consistent multiscale homogenisation study.

It is important to remark that the solvability of the monodomain equations \eqref{eq:monodomain1}-\eqref{eq:monodomain2} 
depends on the properties of $\bD$. In particular,  the stress-assisted diffusion tensor 
needs to remain symmetric and uniformly elliptic, which is a non-trivial condition, given 
the dependence on stress and on voltage. A thorough sensitivity analysis (but for a simpler 
dependence on stress) can be found in \cite{cherubini17}. Here we perform 
a much lighter calibration, as mentioned later in Section~\ref{sec:results}. 
Comparisons between the 
effects of SAD and the more conventional mechanoelectrical feedback through stretch-activated 
currents have been reported in \cite{loppini18}.

\subsection{Activation and excitation-contraction coupling} 
When using the active stress approach, we will adopt a simple description where the active tension is generated 
by ionic quantities (calcium) as well as by local fibre stretch. That is, we {propose} a regularised active tension 
model of the form 
\begin{equation}\label{eq:Ta}
\partial_t T_a =  \hat\alpha \Delta T_a + \ell (T_a,\vec{r}, I_{4,f}) \qquad  \text{in } \Omega\times(0,t_{\mathrm{final}}],
\end{equation}
with  $\hat{\alpha} = \alpha_1 D_0$, and 
$\ell (T_a,\vec{r}, I_{4,f}) = T_a - \alpha_2 r_3 + \alpha_3 I_{4,f}$, where 
$\alpha_1,\alpha_2,\alpha_3 = 0.1\alpha_2$ are positive model constants. 
  As calcium concentration is not readily available in the phenomenological 
cellular model we are employing, we use $r_3$ as a proxy for intracellular calcium \cite{bueno08}. In 
addition,  a linear dependence on the calcium proxy and on the 
local stretch are sufficient in our setting to qualitatively capture the dynamics of active tension.
  
On the other hand, in the framework of active strain, a 
constitutive equation for the activation functions $\gamma_i$ in 
terms of the microscopic cell shortening $\xi$ is expressed as follows (see e.g. \cite{barbarotta18}) 
\begin{align}\nonumber
\gamma_f(\xi) &= \xi, \quad  \gamma_s(\xi) = (1+\xi)^{-1}(1+K_0\xi)^{-1}-1, \\
 \gamma_n(\xi) &= K_0\xi,\label{eq:gammas}
\end{align}
and the specific relation between the myocyte shortening $\xi$ 
and the dynamics of slow ionic quantities (in the context 
of our phenomenological model, only $\vec{r}$) is made precise 
using the law 
\begin{equation}\label{eq:xi}
\frac{d\xi}{dt}  =\hat{\ell}(\xi,\vec{r}) \qquad \text{in } 
\Omega\times(0,t_{\mathrm{final}}],
\end{equation}
which does not require an explicit dependence on local fibre stretch{,} as the 
sliding of myofilaments is driving the dynamics of the functions $\gamma_i$. We employ  
the nonlinear reaction term 
$\hat{\ell}(\xi,\vec{r})= K_1 (1+r_3)^{-1} + K_2 \xi$, and we make the distinction that $\ell$ and $\hat\ell$ 
characterise the evolution of the activation in the approaches of active stress and active strain, respectively.

\section{Numerical method and implementation}\label{sec:FE}
\subsection{Mixed-primal weak form}
Restricting to the case of an active strain model with 
Robin conditions \eqref{eq:robinBC} on the whole boundary for the mechanical layer 
(that is $\partial\Omega_R = \partial\Omega$) 
and the boundary and initial conditions \eqref{bc-mono}-\eqref{init-mono} for the electrical layer, 
 we proceed to take the inner product of  the differential equations \eqref{eq:momentum}, \eqref{eq:mass}, \eqref{eq:Pi}, 
\eqref{eq:monodomain}, \eqref{eq:xi} with adequate test functions, and to integrate by parts 
whenever appropriate. We then arrive at the following weak 
form of the problem: 
For 
$t>0$, find $(\bPi,\bu,p)\in \mathbb{L}^2_{\text{sym}}(\Omega)\times \mathbf{H}^1(\Omega) 
\times \mathrm{L}^2(\Omega)$ and $(v,\vec{r},\xi)\in \mathrm{H}^1(\Omega)\times \mathrm{L}^2(\Omega)^3\times\mathrm{L}^2(\Omega)$ such that 
\begin{widetext}\begin{align}
\int_\Omega [\bPi - \cG + pJ\bI] : \btau & = 0 & \forall \btau \in \mathbb{L}^2_{\text{sym}}(\Omega),\nonumber\\
\int_{\Omega} \rho\partial_{tt}\bu\cdot\bv +\int_\Omega \bPi \bF^{\tt -t}  : \nabla\bv+ \int_{\partial\Omega} \eta  \bF^{-\tt t} \bu\cdot\bv 
& = \int_\Omega \rho_0 \bb\cdot\bv &  \forall \bv \in \mathbf{H}^1(\Omega),\nonumber\\
\int_\Omega [J  - 1 ] q & = 0 &  \forall q\in \mathrm{L}^2(\Omega),\label{eq:weak-form} \\
\int_\Omega {\partial_t v}\, w + \int_\Omega  \bD(v,\bF,\bPi)\, \nabla v \cdot \nabla w 
  - \int_\Omega \bigl[ g(v,\vec{r}) + I_{\rm ext}\bigr] w&=0  & \forall w \in  \mathrm{H}^1(\Omega),\nonumber\\
\int_\Omega \bigl({\partial_t \vec{r}} \cdot \vec{s} 
+  {\partial_t \xi}\, \varphi\bigr) - \int_\Omega 
\bigl( \vec{m}(v,\vec{r}) \cdot \vec{s} + \hat\ell(\xi,\vec{r}) \varphi \bigr) & = 0 & \forall (\vec{s},\varphi)\in \mathrm{L}^2(\Omega)^3\times L^2(\Omega),\nonumber
\end{align}
\end{widetext}
where $\mathbb{L}^2_{\text{sym}}(\Omega) := \{ \btau \in  \mathbb{L}^2(\Omega): \btau = \btau^{\tt t}\}$, and {where
 the case for an active stress formulation necessitating an active tension model is addressed similarly} (however,  
the regularity of $T_a(t)$ is then $\mathrm{H}^1(\Omega)$). 
Theoretical aspects regarding the coupling of elasticity and stress-assisted diffusion problems has been recently addressed in the context of mixed-primal and 
mixed-mixed formulations in \cite{gatica18b}, but only for the case of simplified linear three-field elasticity and steady diffusion. 

{The spatial discretisation will 
follow a mixed-primal Galerkin approach based on the weak formulation 
\eqref{eq:weak-form}. Details on the finite-dimensional spaces and 
linearisation procedure are laid out in Appendix \ref{sec:app-mixed}. }

The motivation for using three-field elasticity formulations {is the need to produce} 
robust solutions with balanced convergence orders for all variables. In addition, 
these methods are robust in the incompressible regime{;} they are not subject to volumetric 
locking \cite{lamichhane12}; and most importantly, they provide direct approximation 
of variables of interest, albeit at a higher computational cost. 
Another advantage of using the Kirchhoff stress is that this tensor is symmetric, and, for 
simpler material laws, is a polynomial function of the displacements (whereas first and second 
Piola-Kirchhoff stresses are rational functions of displacement) \cite{chavan07}. 
Solving in terms of stresses proves particularly useful, as this 
variable participates actively in the electro\-mechanical cou\-pling through the stress-a\-ssis\-ted diffusion.
Moreover, for the lowest-order method characterised by $l=0$, the matrix 
system associated with \eqref{eq:FE} has fewer unknowns than the discretisation 
that uses piecewise qua\-dratic and continuous displacement approximations and 
piecewise linear and discontinuous pressure approximations (and which is a popular  
locking-free scheme for hyperelasticity in the dis\-place\-ment-pre\-ssure for\-mulation, utilised 
for stress-a\-ssis\-ted di\-ffusion problems in the recent work \cite{loppini18}). 
The importance of casting the 
equations of motion in terms of the coupling variables  has been already emphasized in \cite{ruiz15} in 
the context of cardiac electromechanics, which demonstrates 
that the computation of output indicators of interest (such as conduction velocities) 
may suffer from loss of accuracy if one simply postprocesses 
stress or strain from discrete displacements as approximations in the geometric feedback.

\subsection{Solver structure and implementation details}
According to the fixed-point separation between 
electrophysiology and solid deformation solvers, 
nonlinear mechanics will be solved using the Newton-Raphson method stated above, 
and an operator splitting algorithm will separate an implicit diffusion 
solution (where another Newton iteration handles the nonlinear
self-diffusion) from an explicit reaction step for the kinetic
equations, turning the overall solver into a semi-implicit
method. Such a strategy is feasible since the Jacobians associated with the reaction 
and excitation-contraction models do not possess highly varying eigenvalues (otherwise one would need to 
include these terms in the Newton iteration). 
Updating and storing of the internal variables $\xi$ and $\vec{r}$ will be done locally at the quadrature points. 
We solve the linear systems  arising at each Newton iterate by the Krylov iterative method GMRES,  
preconditioned with an incomplete LU(0) factorisation (except for the linear systems in the convergence tests in Section~\ref{sec:accuracy}, 
which will be solved with the direct method SuperLU),
and the iterates are terminated once a tolerance of $10^{-6}$ 
(imposed on the $\ell^\infty$ norm of the non-pre\-con\-ditioned residual) has been a\-chieved. 
The mass matrices associated with the discretisation of the monodomain equations 
are assembled in a lumped manner, which reduces the amount of artificial 
diffusion and violation of the discrete maximum principle \cite{quarteroni17}. 
All routines have been 
implemented using the finite element library FEniCS \cite{fenics}. 

\section{Computational results}\label{sec:results}

\begin{table*}
\rowcolors{1}{gray!30}{gray!10}
\begin{center}
\subfigure[Hyperelasticity variables]{
\begin{tabular}{|l|c|c|c|c|c|c|c|}
\hline
DoF  & $h$ &  $\|\bPi-\bPi_h\|_{0,\Omega}$  &  \texttt{rate}  &  $\|\bu-\bu_h\|_{1,\Omega}$  &  \texttt{rate} &  $\|p-p_h\|_{0,\Omega}$  
&   \texttt{rate}   \\
\hline
\multicolumn{8}{|c|}{$l=0$}\\
\hline
77   & 0.7071 & 43.252 &    --         & 0.0576 &    --         & 30.161 &  --  \\
253 & 0.3536 & 27.137 & 0.6725 & 0.0342 & 0.6345 & 19.030 & 0.6647  \\
917 & 0.1768 & 12.535 & 1.1140 & 0.0216 & 0.7615 & 9.2110 & 1.0471 \\
3493 & 0.0884 & 6.2636 & 1.0012 & 0.0118 & 0.8751 & 4.8012 & 0.9401  \\
13637 & 0.0442 & 1.9169 & 1.1727 & 0.0071 & 0.9516 & 1.9631 & 1.3817 \\
53893 & 0.0221 & 0.9841 & 0.9907 & 0.0042 & 0.9737 & 0.9206 & 0.9858 \\
\hline
\multicolumn{8}{|c|}{$l=1$} \tabularnewline
\hline
221 & 0.7071 & 19.481 &        -- & 0.0146 &   --           & 6.0355 &    -- \\
789 & 0.3536 & 7.9032 &  1.3034 & 0.0037 &  1.7593 & 1.5809 &   1.4581 \\
2981 & 0.1768 & 2.6409 &  1.8079 & 0.0011 &  1.7809 & 0.4120 &   1.7269 \\
11589 & 0.0884 & 0.7277 &  1.9033 & 4.11E-4 & 1.8065 & 0.1353 &   1.8813 \\
45701 & 0.0442 & 0.2063 &   1.9182 & 1.09E-4 &  1.9330 & 0.0382 &   1.9602 \\
181509 & 0.0221 & 0.0569 &  1.9466 & 3.12E-5 &  1.9522 & 0.0094 &   1.9571 \\
\hline
\end{tabular}}
\subfigure[Electrophysiology variables]{\begin{tabular}{|l|c|c|c|c|c|c|c|}
\hline
DoF  & $h$ &  $\|v-v_h\|_{1,\Omega}$  &  \texttt{rate} & $\|r-r_h\|_{1,\Omega}$ &  \texttt{rate} & $\|T_a-T_{a,h}\|_{1,\Omega}$ &  \texttt{rate} \\  
\hline
\multicolumn{8}{|c|}{$l=0$}\\
\hline
77 & 0.7071 & 0.1528 &  -- & 0.1926 &    -- & 0.1623 &  -- \\
253 & 0.3536 & 0.0902 & 0.7601 & 0.1069 & 0.8499 & 0.0847 & 0.8824 \\
917 & 0.1768 & 0.0491 & 0.8769 & 0.0573 & 0.8968 & 0.0433 & 0.9673 \\
3493 & 0.0884 & 0.0282 & 0.8016 & 0.0317 & 0.9536 & 0.0218 & 0.9896 \\
13637 & 0.0442 & 0.0153 & 0.9304 & 0.0172 & 0.9612 & 0.0121 & 0.9446\\
53893 & 0.0221 & 0.0084 & 0.9587 & 0.0091 & 0.9843 & 0.0067 & 0.9562\\ 
\hline
\multicolumn{8}{|c|}{$l=1$}\\
\hline
221 & 0.7071 & 0.0329 &  --                & 0.0583 &  --        & 0.0469 &  -- \\
789 & 0.3536 & 0.0102 &  1.5043       & 0.0152 & 1.7317 & 0.0133 &   1.6300 \\
2981 & 0.1768 & 0.0029 &  1.7608     & 0.0039 &  1.8809  & 0.0035 &   1.7095 \\
11589 & 0.0884 & 8.03E-4 &  1.7849  & 0.0010 &   1.9021 & 9.25E-4 &   1.8822 \\
45701 & 0.0442 & 2.31E-4 &  1.8964  & 2.70e-4 &  1.8966 & 2.41E-4 &   1.8907 \\
181509 & 0.0221 & 6.11E-5 & 1.9598  & 7.05e-5 &  1.9604 & 6.86E-5 &   1.9649 \\
\hline
\end{tabular}}
\end{center}
\caption{Test 1: Error history (errors on a sequence of successively refined grids and convergence rates) associated with the mixed finite element method \eqref{eq:FE} applied to a steady-state electromechanical coupling under 
active stress, and using different 
polynomial degrees $l\in\{0,1\}$.} \label{table01}
\end{table*}

\subsection{Mesh convergence}\label{sec:accuracy}  We begin with the numerical validation of our mixed-primal method 
on a problem slightly simpler than \eqref{eq:mech} - \eqref{eq:monodomain} - \eqref{eq:Ta}, but 
that still retains the main ingredients of the model. These include orthotropic active mechanics, nonlinear reaction-diffusion 
with stress-assisted diffusion, and a nonlinear excitation-contraction coupling.

A convergence test is generated by computing errors between smooth exact solutions 
and approximate solutions using the first-order and the second-order methods discussed 
 in Section~\ref{sec:FE}. Let us consider the following closed-form solutions to 
 a steady-state 
 counterpart of the variational form \eqref{eq:weak-form} for the electromechanics equations, also 
 assuming the absence of viscoelastic 
 effects, and  
defined on the domain $\Omega=(0,1)^2$ with the fibres/sheetlets defined as $\fo =(0,1)^{\tt t},\so = (-1,0)^{\tt t}$ 
\begin{align*}
\bu(x,y) & = 0.1\begin{pmatrix}
\sin(\pi x)\cos(\pi y)\\
\cos(\pi x)\sin(\pi y)\end{pmatrix},\\ 
p(x,y) &= 0.1\sin(\pi x)\sin(\pi y),\\ 
 v(x,y) &= 1+0.1\cos(\pi x)\cos(\pi y),\\
r(x,y) &= 0.1\cos(\pi x)\sin(\pi y)\sin(\pi x),\\ T_a(x,y) &= 1+0.1\cos(\pi x)\sin(\pi y){.}
\end{align*}
Then the Kirchhoff stress $\bPi$, as well as suitable forcing terms (volumetric load{,} an additional external stimulus, and 
the active tension source) are computed from these smooth solutions, 
the balance equations, relations \eqref{eq:total-stress}, \eqref{eq:active-stress}, \eqref{def:D}, 
and using the following simplified constitutive equations 
\[
m(v,r) = v-r^2, \  g(v,r) = (v-1)vr, \ \ell(T_a,r) = - T_a + r.
\]
Note also that the incompressibility constraint for this test is $J=J_{ex}$, where $J_{ex}$ is computed from 
the exact displacement. {Here we also prescribe Dirichlet boundary conditions for displacements, 
transmembrane potential, and active tension (incorporated in the discrete trial spaces).} Errors due to fixed-point iterations are avoided by 
taking a full monolithic coupling and computing solutions using Newton-Raphson iterations 
with an exact Jacobian. 
On a sequence of six uniformly refined meshes, we proceed to compute 
errors between the exact and approximate solutions computed with methods using $l=0$ and 
$l=1$. Kirchhoff stress and pressure errors 
are measured in the $L^2-$norm, whereas for the remaining variables the errors are measured 
in the $H^1-$norm. The obtained error history is reported in
Table~\ref{table01}, where we observe an 
asymptotic $O(h^{l+1})$ decay of the error for each field variable. This behaviour   
corresponds to the optimal convergence according to the interpolation properties of the 
employed finite element subspaces \cite{chavan07}.

\begin{table*}
\setlength{\tabcolsep}{3pt}
\rowcolors{1}{gray!30}{gray!10}
\bigskip{}\centering{}
\begin{tabular}{|rlc| lrlcl|rlcl|rlc|}
\hline
\multicolumn{15}{|c|}{Viscoelasticity constants}\\
\hline
$a=$    & 0.236 & [N/cm$^2$] && \quad $a_{f}=$ & 1.160 & [N/cm$^2$] & & \quad $a_{s}=$ & 3.724 & [N/cm$^2$] && \quad $a_{fs}=$ & 4.010 & [N/cm$^2$] \\
$b=$    & 10.81& [--]        && \quad $b_{f}=$ & 14.15 & [--] & & \quad $b_{s}=$ & 5.165& [--] && \quad $b_{fs}=$ & 11.60 & [--] \\
$p_0=$& 0.1 & [N/cm$^2$] &&\quad $\beta=$ & 10 & [ms] & & \quad $\delta=$& 22.6 & [N/cm$^2$\,ms] &&\quad $\zeta_{\text{stab}}=$ & 0.25 & [--]\\  
 $\eta_a=$ & 0.001 & [N/cm$^2$] &&\quad $\eta_b=$& 0.01 & [N/cm$^2$] && \quad $\kappa_{sn}=$     & 0.6 &[--] && \quad $\kappa_{nn}=$ & 0.03 & [--] \\
 $\rho_0=$ & 0.001 & [N/cm$^2$] && && && && &&&&\\
\hline
\multicolumn{15}{|c|}{Electrophysiology constants}\\
\hline
$v_0=$   & 0 & [--]&&\quad $v_v=$ & 1.55 & [--] & & \quad $v_2^-=$   & 0.03 &[--] &&\quad $v_{so}=$ & 0.65 & [--]\\
$v_3=$   & 0.908 & [--]&&\quad $\theta_1=$ & 0.3 & [--] & & \quad $\theta_1^-=$   & 0.006 &[--] &&\quad $\theta_{o}=$ & 0.006 & [--]\\
$\theta_2=$   & 0.13 & [--]&&\quad $k_2^-=$ & 65 & [--] & & \quad $k_3=$   & 2.099 &[--] &&\quad $k_{so}=$ & 2.045 & [--]\\
${r^*_{2,\infty}}=$   & 0.94 & [--]&&\quad $\tau_{2,\infty}=$ & 0.07 & [--] & & \quad $\tau_{1,1}^-=$   & 60 &[--] &&\quad $\tau^-_{1,2}=$ & 1150 & [--]\\
$\tau^-_{2,1}=$   & 60 & [--]&&\quad $\tau^-_{2,2}=$ & 15 & [--] & & \quad ${\tau_{fi}}=$   & 0.11 &[--] &&\quad ${\tau_{o,1}}=$ & 30.02 & [--]\\
${\tau_{o,2}}=$   & 0.996 & [--]&&\quad $\tau_{so,1}=$ & 2.046 & [--] & & \quad $\tau_{so,2}=$   & 0.65 &[--] &&\quad $\tau_{3,1}=$ & 2.734 & [--]\\
$\tau_{3,2}=$   & 16 & [--]&&\quad ${\tau_{si}}=$ & 1.888 & [--] & & \quad $\tau^+_{1}=$   & 1.451 &[--] &&\quad $\tau^+_{2}=$ & 200 & [--]\\
$\chi=$ & 1 & [--] && && && && &&&&\\
\hline
\multicolumn{15}{|c|}{Activation and excitation-contraction coupling constants}\\
\hline
$D_0=$ & 1.171 & [mm$^2$/s] &&\quad $D_1=$ & 0.9 & [mm$^2$/s] &&\quad $D_2=$ & 0.01 & [mm$^2$/s]   &&$K_0=$   & 5 &[--] \\
$K_1=$ & -0.015 & [--] &&\quad  $K_2=$ & -0.15 & [--] && $\alpha_1=$ & 10 & [--] &&\quad  $\alpha_2=$ & 0.5 & [--]  \\
\hline\noalign{\smallskip}
\end{tabular}
\caption{Model parameters for the electro-viscoelastic model \eqref{eq:mech}, \eqref{eq:monodomain}, \eqref{eq:xi}, \eqref{eq:Ta}. Values 
are taken from \cite{cherubini17,gao,rossi14,bueno08}, and the transmembrane potential $v$ is  dimensionless.}
\label{table:params}
\end{table*}

\subsection{Parameter calibration}\label{sec:caliber} 
\begin{figure*}
\begin{center}
\subfigure[$t=800\,\text{ms}$]{\includegraphics[width=0.325\textwidth]{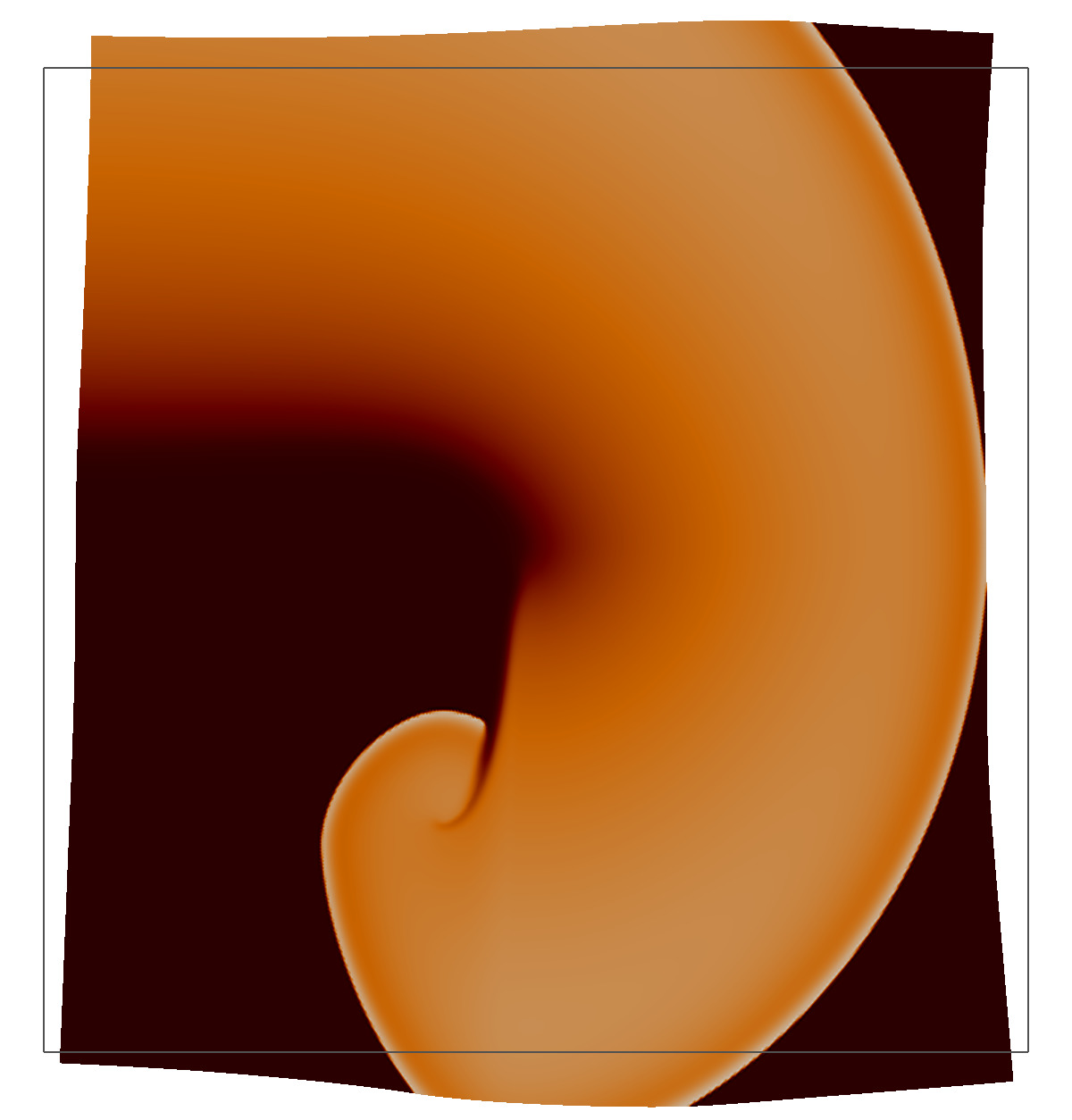}}
\subfigure[$t=900\,\text{ms}$]{\includegraphics[width=0.325\textwidth]{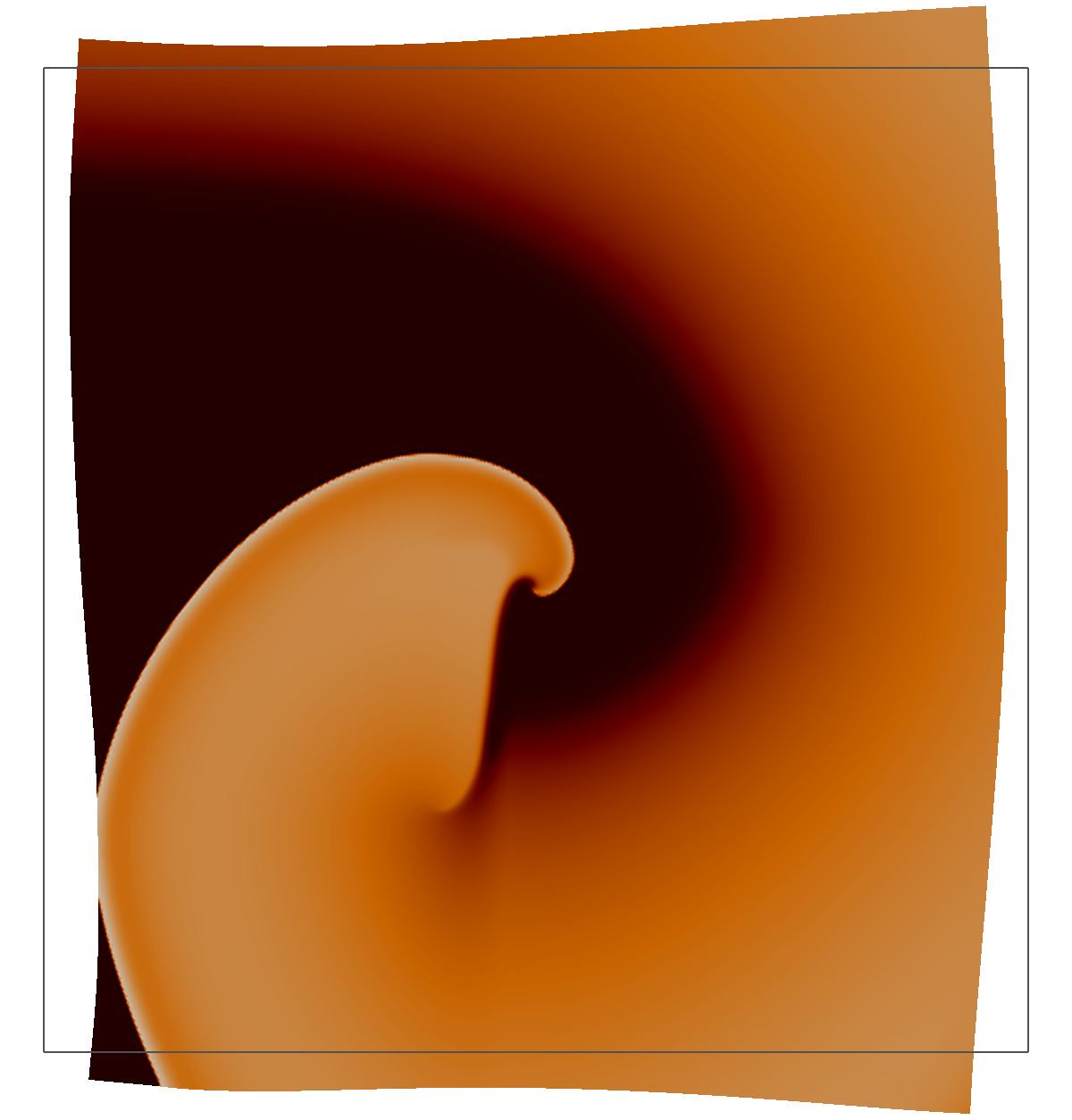}}
\subfigure[$t=1000\,\text{ms}$]{\includegraphics[width=0.325\textwidth]{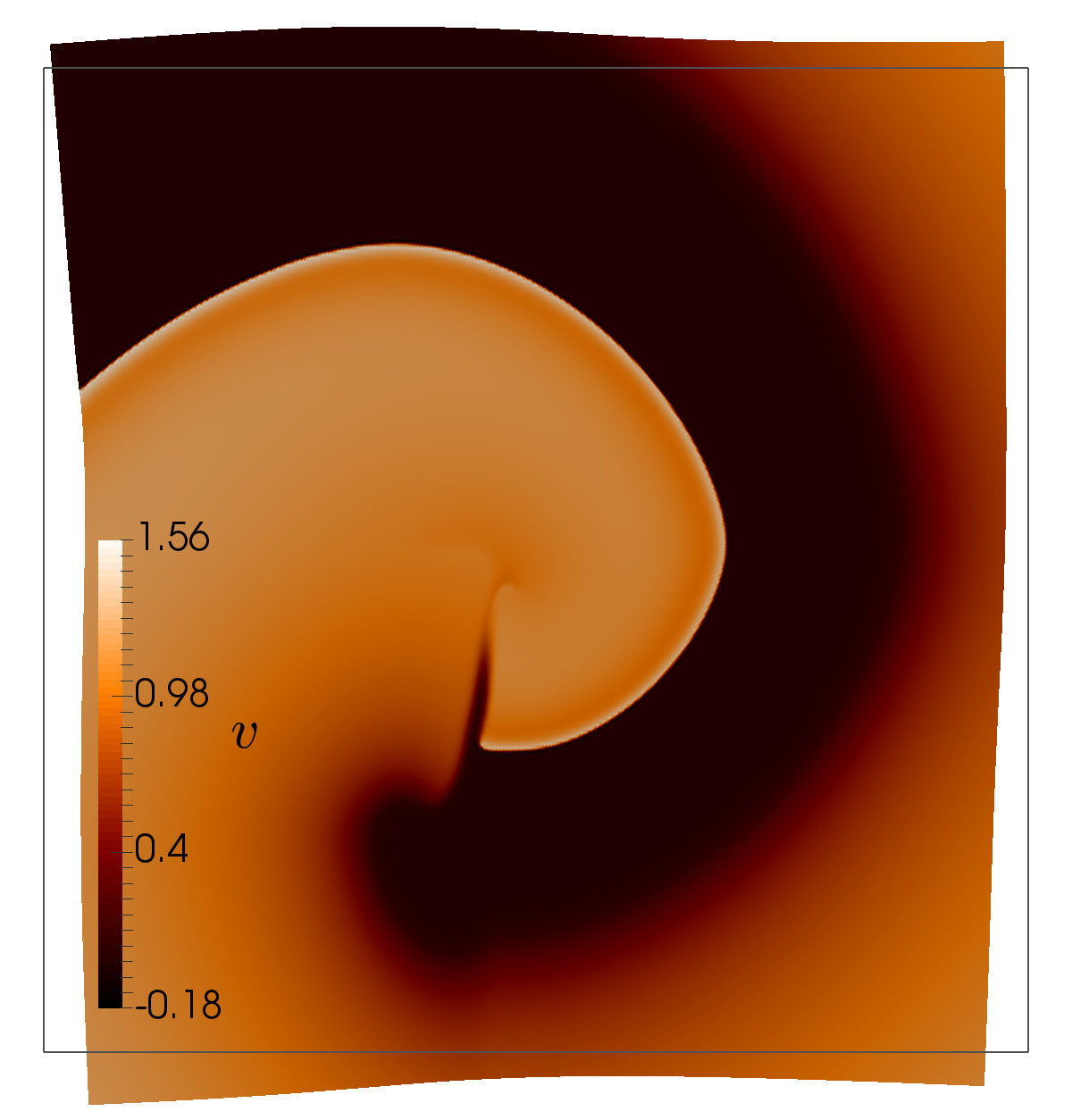}}
\end{center}
\vspace{-4mm}
\caption{Evolution of voltage after S2 stimulus, showing formation of a reentrant spiral wave on the deforming viscoelastic tissue, {computed using the active stress approach}.}
\label{fig:2Dspirals}
\end{figure*}

For the following 2D simulations, we will initially consider tissue slabs of $50 \times 50 $ mm$^2$, 
and set fibre and sheetlet directions simply as $\fo = (1,0)^{\tt t}$, $\so = (0,-1)^{\tt t}$. 
The initiation, maintenance, prevention and treatment of so-called reentrant waves is a major focus of current research due to their implication in atrial and ventricular fibrillations \cite{ColliBook}. We are thus interested in investigating the formation of spiral reentrant waves in our model setup, following the S1-S2 stimulation protocol. We first excite the tissue with a symmetric stimulus labelled S1. An asymmetric stimulus labelled S2 is then applied during the vulnerable window near the end of the refractory period, when some of the tissue has recovered excitability but depolarisation is still blocked elsewhere. This causes unbalanced excitation, which can lead to the formation of a spiral wave. We will define the spiral front as the edge of the spiral wave, where the excitation front meets the repolarisation waveback of the action potential. In our simulations, both waves have nondimensional amplitude 3 and duration 3\,ms. The S1 stimulus is a planar wave created by exciting the entire left edge of the tissue in 2D, and the entire bottom section (below some value of the $z-$coordinate) in 3D. The S2 stimulus is a square wave created by exciting the bottom left quadrant  at $t = 330$ ms and the bottom left octant at $t = 335$ ms in 2D and 3D, respectively. {Here we use the active stress approach, and the boundary conditions for the visco-elastodynamic equations 
correspond to \eqref{eq:robinBC}.}
The formation and evolution of the spiral wave on a deforming domain can be seen in Figure~\ref{fig:2Dspirals}. The spiral is initiated by the diffusion of voltage and transport of ionic entities from the S2 stimulus into the leftmost section of the tissue, which has recovered enough excitability after S1. The wave then spreads outwards in all directions, occupying the entire tissue except for the region that was just excited by the S2 wave. 

\begin{figure*}
\begin{center}
\subfigure[]{\includegraphics[height=0.395\textwidth]{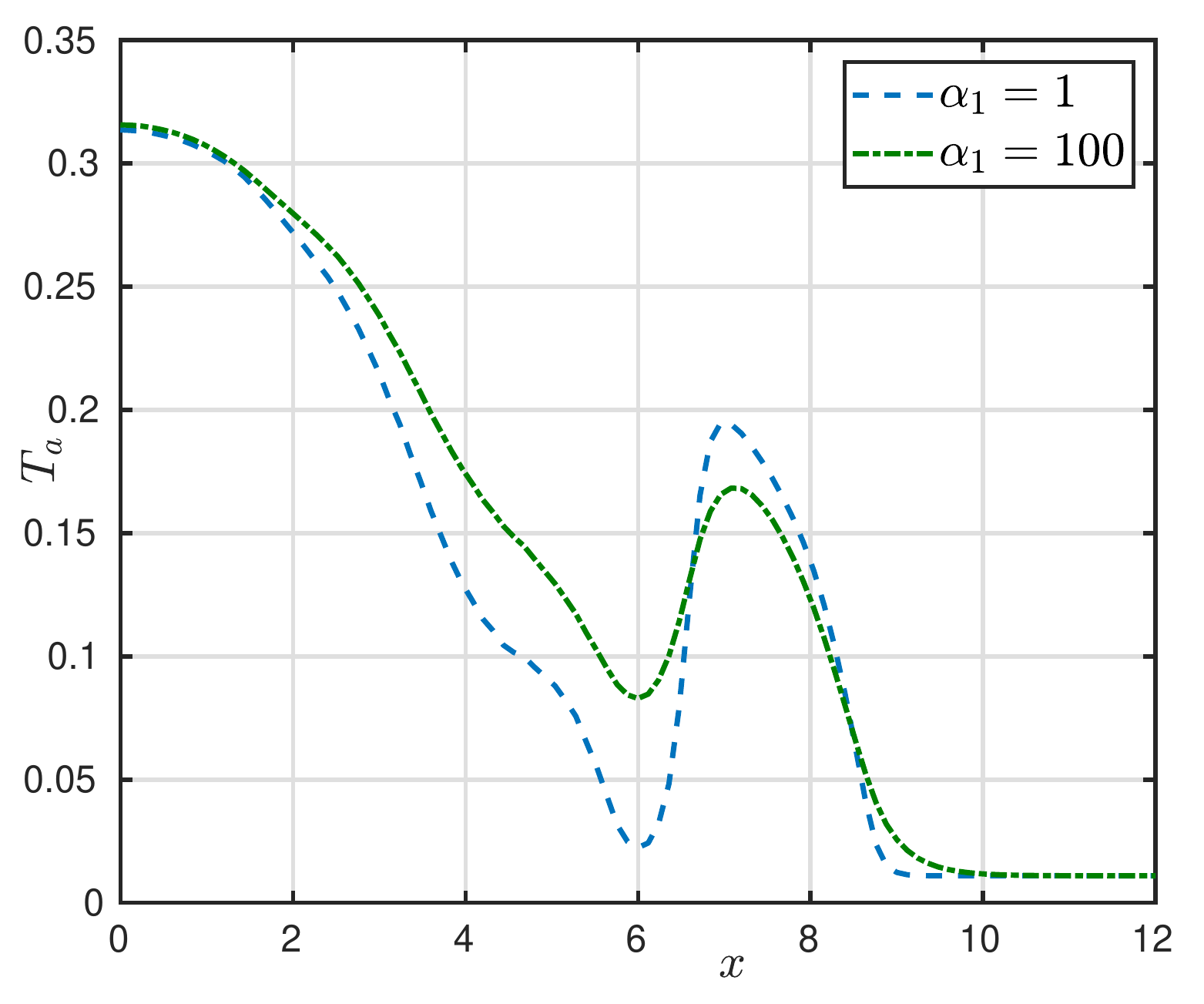}}
\subfigure[]{\includegraphics[height=0.395\textwidth]{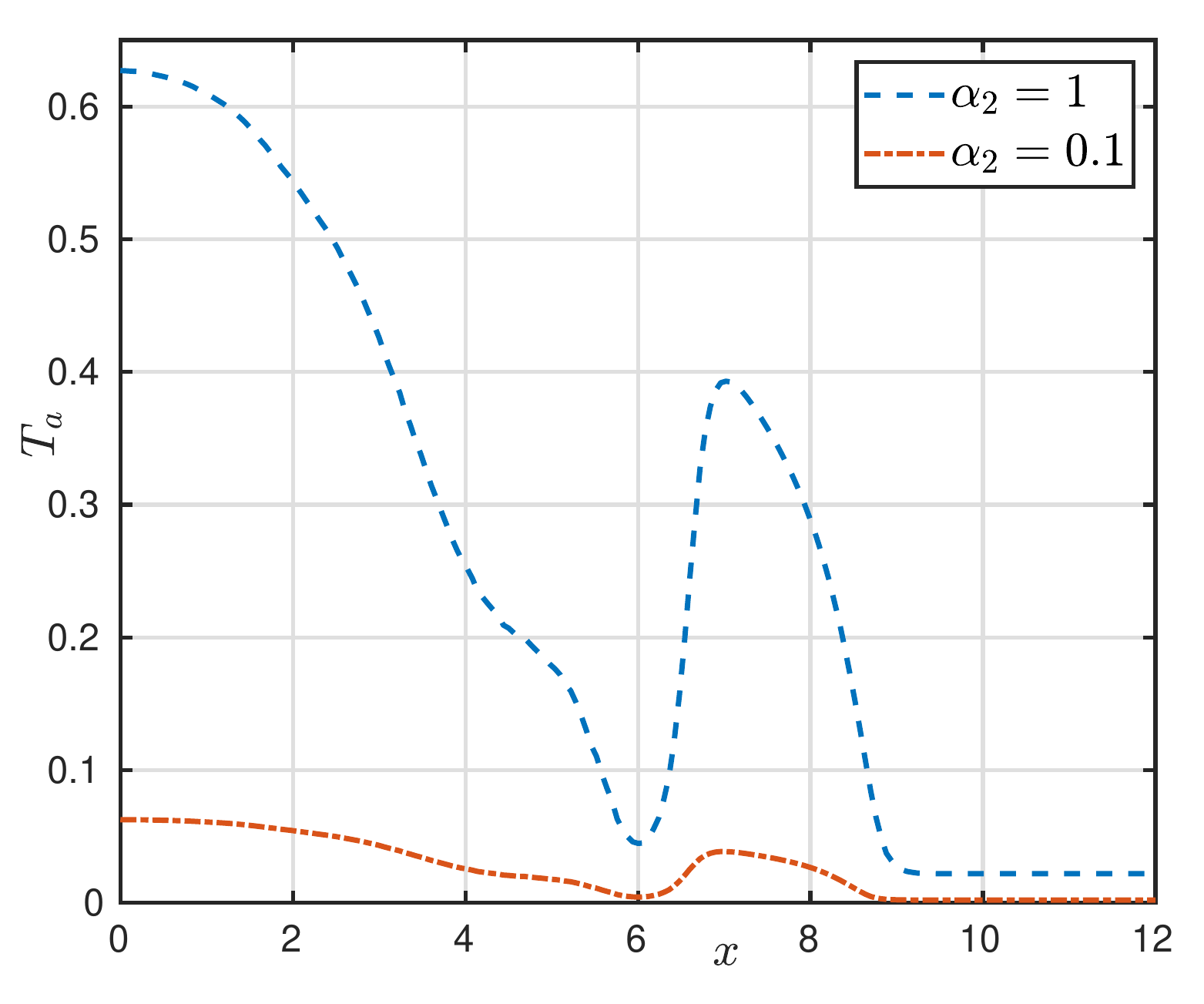}}
\end{center}
\vspace{-4mm}
\caption{Profiles of $T_a$ taken across {a smaller slab of} tissue at $y=6\,\text{mm}$ and $t=432\,\text{ms}$. These plots evaluate the effect of  $\alpha_1$ and $\alpha_2$.}
\label{fig:alphas}
\end{figure*}

Next, {since we are using the active stress formulation in this case}, we proceed to evaluate $\alpha_1, \alpha_2$, the parameters governing active tension in \eqref{eq:Ta}, and $\eta$, the stiffness parameter from \eqref{eq:robinBC}. We conduct a simple sensitivity analysis by increasing or decreasing either $\alpha_1, \alpha_2$ or $\eta$ by one order of magnitude, holding the others constant at their reference values ($\alpha_1 = 10$, $\alpha_2 = 0.5$, and $\eta = \eta_a = 0.001$N$/$cm$^2$, as listed in Table~\ref{table:params}). 
This simple analysis therefore does not test for compounding or interaction effects. We also consider a smaller slab of size $12 \times 12 $ mm$^2$. The parameter $\alpha_1$ contributes to producing 
smoother active tension profiles, while $\alpha_2$ controls the 
range of their magnitude. These effects are visible in Figure~\ref{fig:alphas}. We found that larger values of $\alpha_1$ produced smoother gradients in pressure and stress, while larger values of $\alpha_2$ produced, in average, higher magnitude displacement, Kirchhoff stress, and pressure, as well as some more subtle changes in ionic quantities.
Parameter $\eta$ determines the stiffness of the springs supporting the tissue, and so decreasing $\eta$ resulted in an increase in the maximum values of magnitude of displacement, stress, and pressure, as expected. However, these 
differences were minimal, even across the three orders of magnitude tested ($\eta = 1$E-4 to $\eta = 0.01$). The effects on ionic entities were even smaller,  for both the hyperelastic and viscoelastic cases, and therefore plots are not shown. 

Computational experiments reveal a window of values of $D_2$ for which our method converges. In the 2D hyperelastic case, we found that the upper bound for $D_2$ is approximately $D_2 = 2.1$E-2\,mm$^2/$s, with the linear solver failing to converge for larger values. In these simulations, the Kirchhoff stress achieved an $L^2-$norm of between 0.006 and 0.6. In turn, the viscoelastic case was able to accept slightly larger values of $D_2$, up to $D_2 = 2.2$E-2\,mm$^2/$s, with the $L^2-$norm of stress falling between 0.001 and 0.5. A possible explanation is the loss of coercivity or monotonicity in the 
stress-assisted diffusion coupling, as explored in \cite{cherubini17}. 

{The numerical method used for these tests is characterised by the time step, mesh size, polynomial degree, and stabilisation constant $\Delta t= 0.01$\,ms, $h=0.3534$\,mm, $\ell = 0$, $\zeta_{\mathrm{stab}} = 2.5$, respectively}. 

\subsection{Locking-free property}
We next proceed to assess the performance of the proposed mixed formulation for the mechanical problem. In this example we solve only 
for \eqref{eq:mech} without the acceleration term (otherwise present in all other simulations),  
using the active stress approach with a fixed value for the active tension  and 
without the contribution from the viscous stress \eqref{eq:beta}. We proceed to compare the deformation achieved by the mixed 
formulation with that of an 
asymptotic solution and the approximate solution generated by a more standard pressure-displacement finite element formulation. We consider different {stabilisation} parameter values and mesh refinements.

\begin{figure*}
\begin{center}
\subfigure[]{\includegraphics[height=0.4\textwidth]{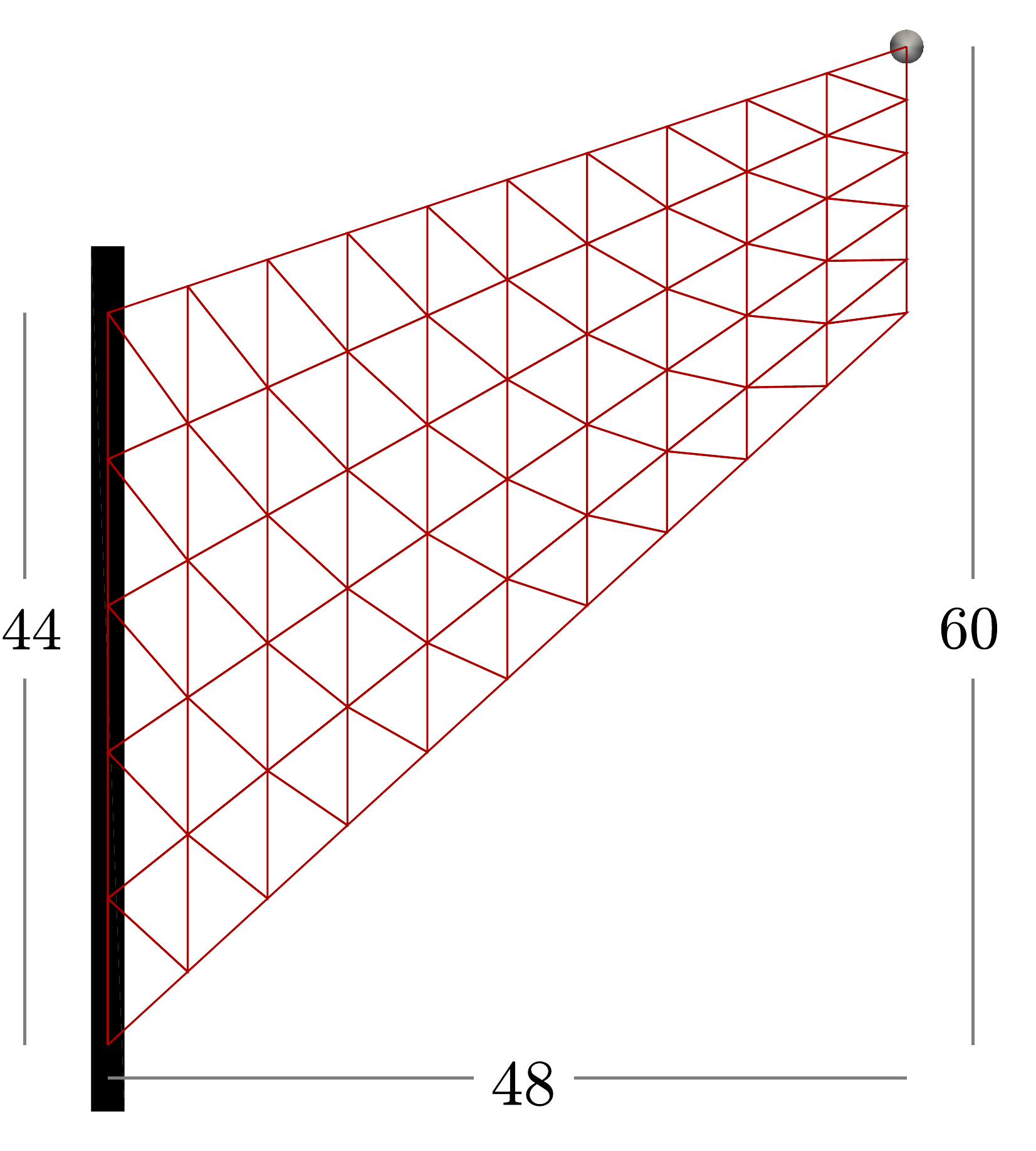}}
\subfigure[]{\includegraphics[height=0.4\textwidth]{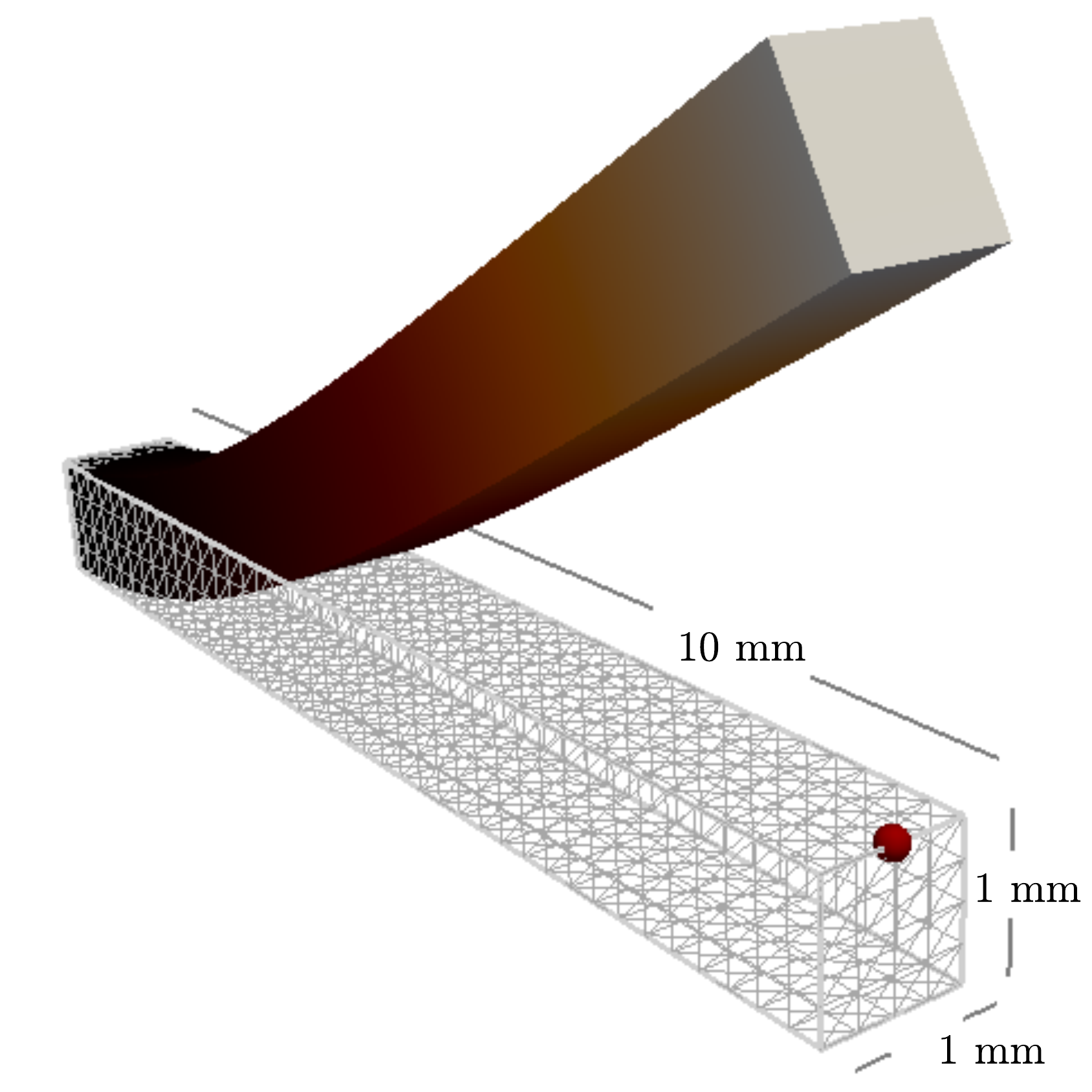}}
\end{center}
\vspace{-4mm}
\caption{Domain sketches and sample meshes for the deflection of Cook's membrane for an 
Holzapfel-Ogden material with constant active stress (a) and 
deflection of a 3D beam for a Guccione-Costa-McCulloch material with the active stress component set to zero (b). }
\label{fig:locking-sketch}
\end{figure*}

\begin{figure*}[t!]
\begin{center}
\subfigure[]{\includegraphics[height=0.375\textwidth]{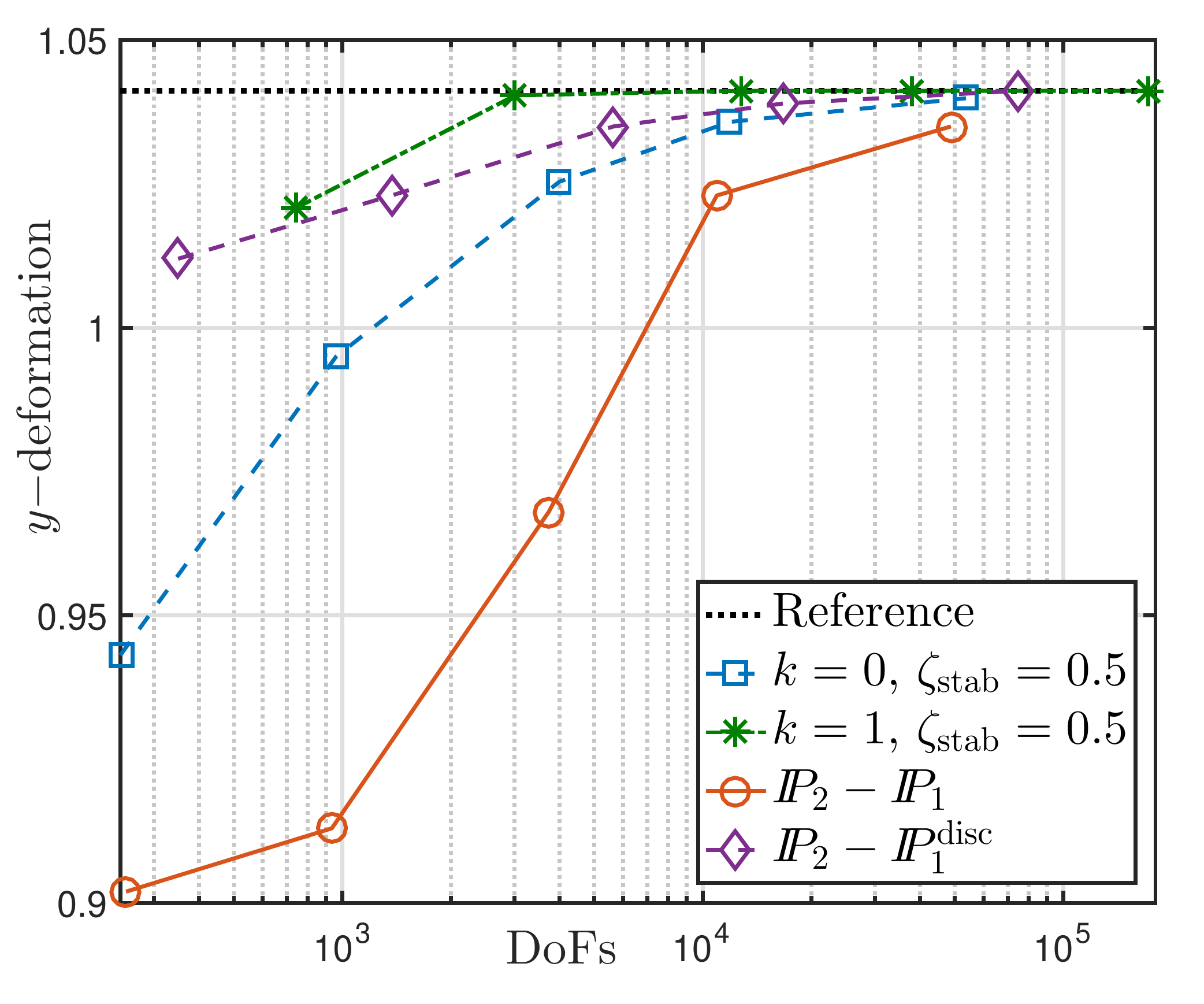}}
\subfigure[]{\includegraphics[height=0.375\textwidth]{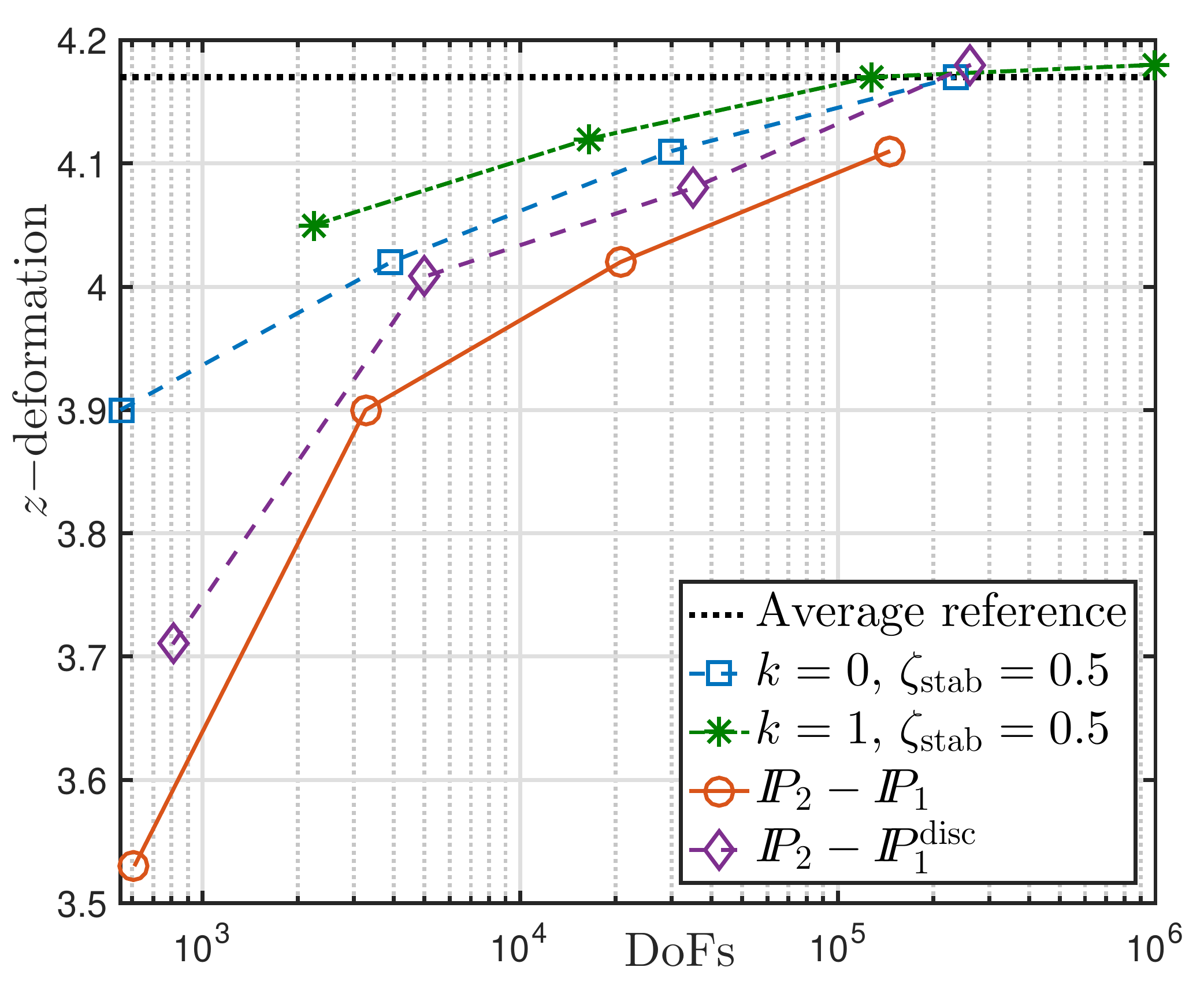}}\\
\subfigure[]{\includegraphics[height=0.375\textwidth]{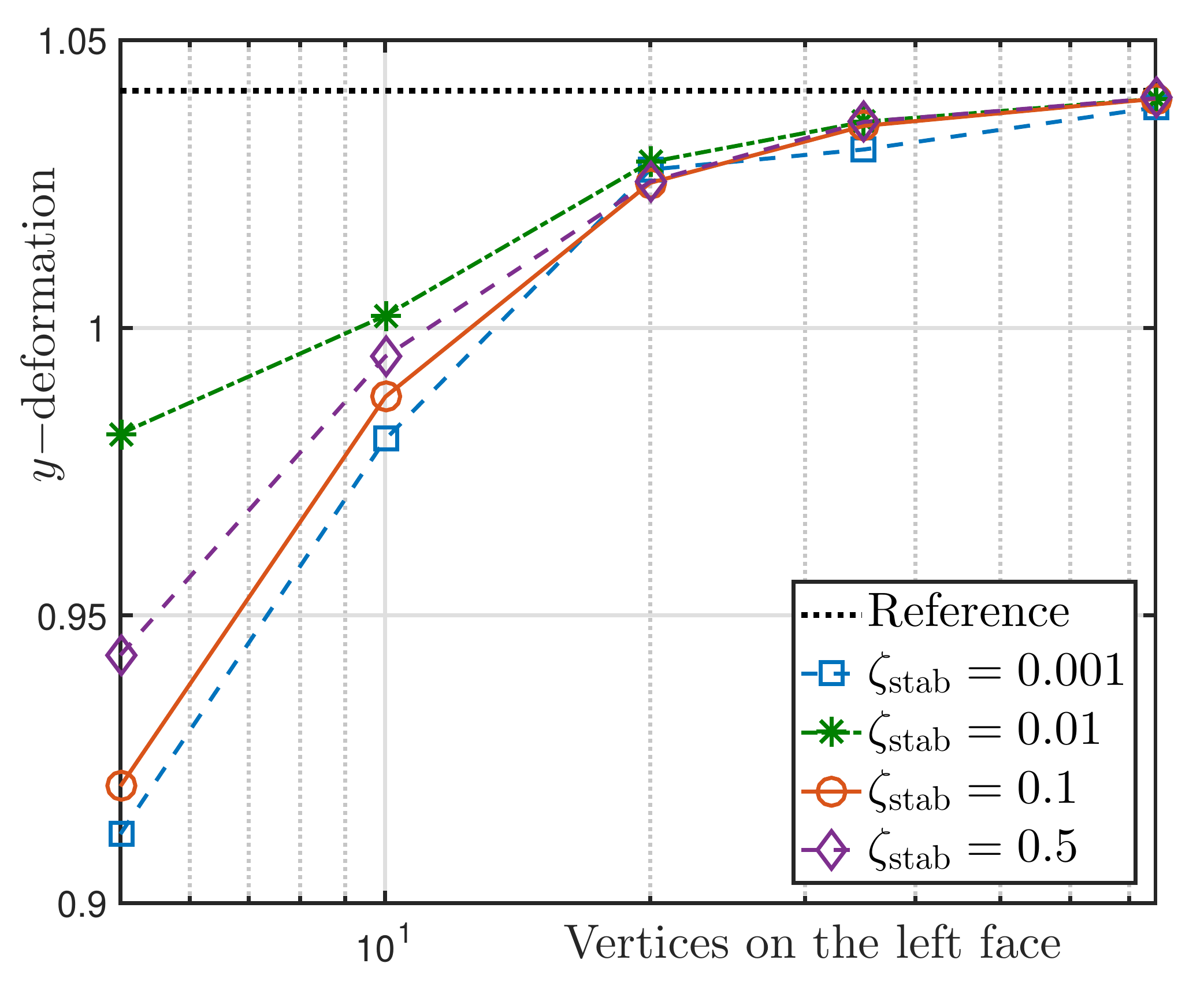}}
\subfigure[]{\includegraphics[height=0.375\textwidth]{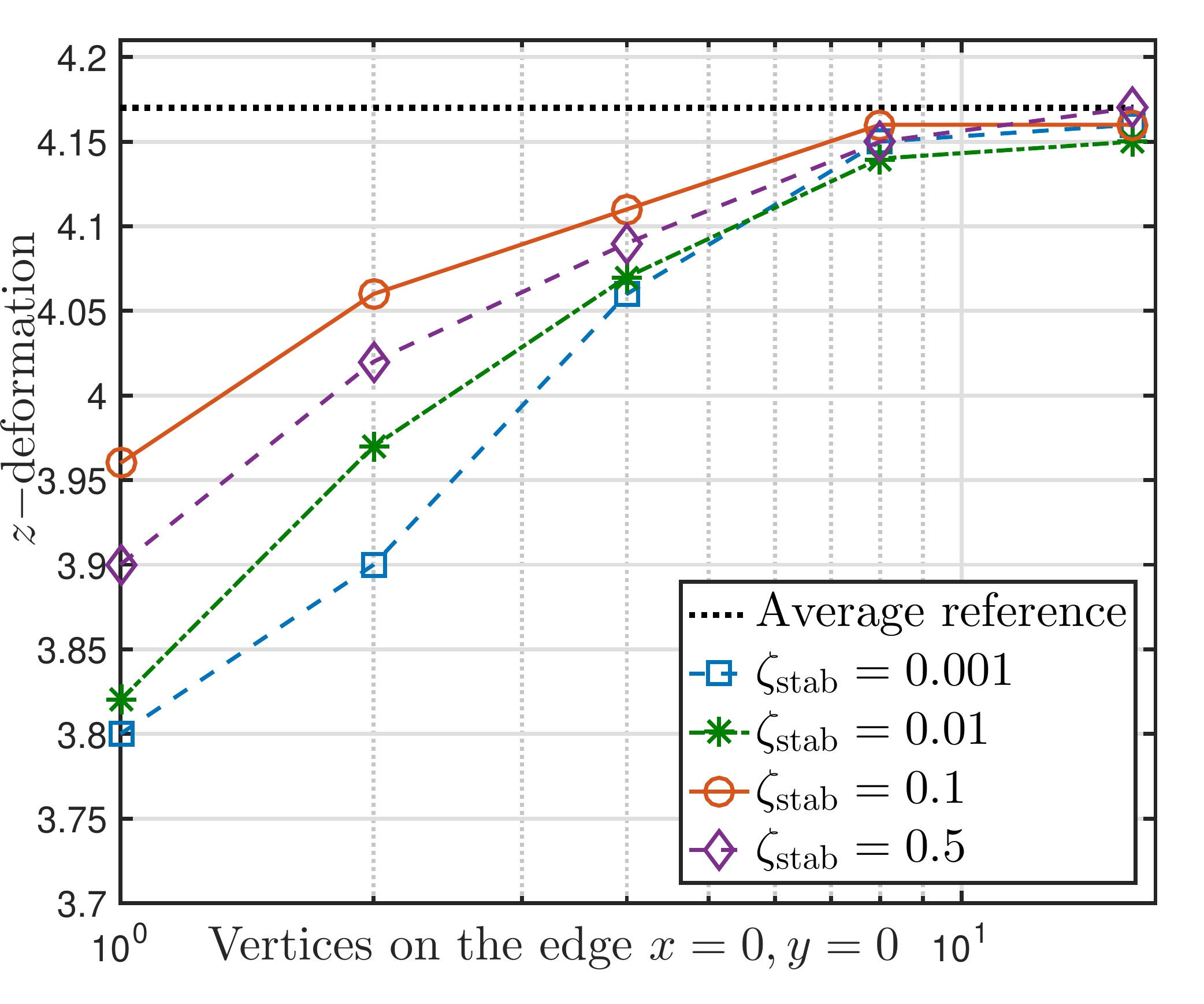}}
\end{center}
\vspace{-4mm}
\caption{Convergence of the deflection of Cook's membrane for an Holzapfel-Ogden material with constant active stress (a,c) and 
deflection of a 3D beam for a passive Guccione-Costa-McCulloch material (b,d). Maximal vertical deflection 
with respect to the mesh resolution for 
different numerical schemes (a,b), and different values of the stabilisation constant (c,d). }
\label{fig:locking}
\end{figure*}

We perform two sets of computations. First,  we undertake   
Cook's membrane benchmark test for a fully incompressible Holzapfel-Ogden material ({as was similarly done} for nearly incompressible Saint Venant-Kir\-chhoff and Neo-Hookean 
solids in \cite[Test II]{chavan07}), where we set an active tension of $T_a = 0.07$. 
This test involves applying an upward in-plane shear load $\boldsymbol{t}=(0,100)^{\tt t}$ to the right edge of a tapered panel with a clamped left edge, and measuring the vertical deformation of the upper right vertex. The domain is defined as the convex hull of the set $\{(0,0),(48,44),(48,60),(0,44)\}$ (see the 
sketch in Figure~\ref{fig:locking-sketch}(a)), and the fibre and sheetlet fields are $\fo=(1,0)^{\tt t}$ and $\so = (0,-1)^{\tt t}$, respectively. Secondly, we consider a 3D system 
suggested in \cite[Test I]{land15} as a simple benchmark for passive cardiac mechanics, and therefore we set $T_a=0$. The problem consists in computing the deformation 
of a {point at the right end of a beam defined} by the domain $\Omega = (0,10)\times(0,1)\times(0,1)$\, mm (see the 
sketch in Figure~\ref{fig:locking-sketch}(b)), where the fibre direction is $\fo=(1,0,0)^{\tt t}$. Instead of \eqref{piolatensor}, the material 
is characterised by the transversally isotropic strain energy density proposed by Guccione et al. \cite{gucc} (which is the 
material law used in the benchmark test from \cite{land15}): $\Psi_{\text{pas}} = a/2 (e^Q-1)$, with $Q = b_f E_{ff}^2 + b_t (E_{ss}^2+ E_{nn}^2 + E_{sn}^2+E_{ns}^2) 
+ b_{fs} (E_{fs}^2+ E_{sf}^2 + E_{fn}^2+E_{nf}^2)$, 
where $a=2$\,kPa, $b_f=8$, $b_t=2$, $b_{fs}=4$, and the $E_{ij}$ denote entries of the Green-Lagrange strain tensor $\bE$, 
rotated with respect to a local coordinate system aligned with $\fo,\so,\no$. 
The beam is clamped at the face $x=0$, a pressure of $p_N=0.004$\,kPa is imposed on the bottom 
face $z=0$, and the remainder of the boundary is considered with traction-free conditions.  According to \eqref{eq:pressureBC}, the pressure boundary condition changes with the deformed 
surface orientation, and its magnitude scales with the deformed area.

The outcome of these tests in Figures~\ref{fig:locking}(a,b) shows a rapid convergence of our first- and second-order
methods, while the computations using a pre\-ssure-dis\-place\-ment formulation and the Taylor-Hood finite elements 
({the well-known $\mathbb{P}_2-\mathbb{P}_1$ pair of} continuous and piecewise quadratic approximations of displacements and continuous and piecewise linear 
approximations for pressure) display a 
somewhat slower 
convergence to the asymptotic deflection of the membrane. Using discontinuous pressures (the $\mathbb{P}_2-\mathbb{P}_1^{\text{disc}}$ 
pair) rectifies the convergence, but at 
a higher computational cost. Quite similar results were obtained for the beam (where the reference value is the average of the reported simulations 
from the study in \cite{land15}). 
Moreover, Figures~\ref{fig:locking}(c,d) show the vertical deflections as a function of the number of vertices discretising the left side of 
the membrane and of the small edge of the beam, respectively. They 
indicate that the obtained results  
are consistent for varying values of the stabilisation parameter, $\zeta_{\text{stab}}$, and 
the observed behaviour {also confirms} that our method is locking-free.


\subsection{Stress-assisted diffusion and conduction velocity} 
In addition to determining a suitable parameter range for $D_2$ that 
ensures solvability of the discrete monodomain equations, we also investigated the effect of $D_2$ on the 
tissue's response to spiral wave dynamics. As in the second part of Section~\ref{sec:caliber}, this time the domain 
is a square slab of width $12\,{\rm mm}$ aligned with the canonical axes. We employ the active stress approach and use Robin boundary conditions 
for the viscoelasticity problem. The fibres assume the fixed direction $\fo=(1,0)^{\tt t}$ and the sheetlets $\so = (0,-1)^{\tt t}$, and 
the Holzapfel-Ogden material law is considered. The meshsize is approximately $h = 0.085\,$mm and the timestep is $\Delta t = 0.1\,$ms. 
We use the lowest-order finite element method and the stabilisation parameter is  $\zeta_{\text{stab}} = 2.5$.

\begin{figure*}
\begin{center}
\subfigure[]{\includegraphics[width=0.325\textwidth]{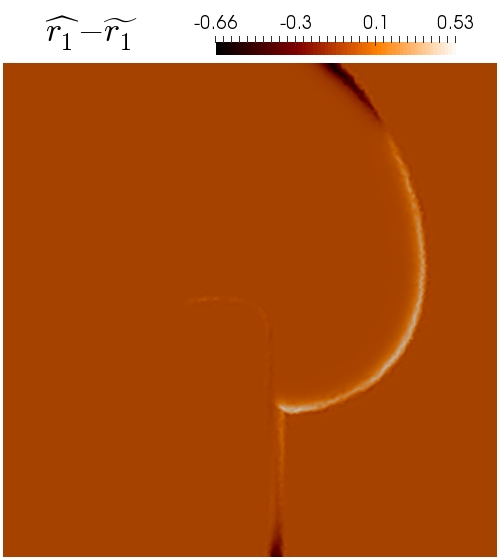}}
\subfigure[]{\includegraphics[width=0.325\textwidth]{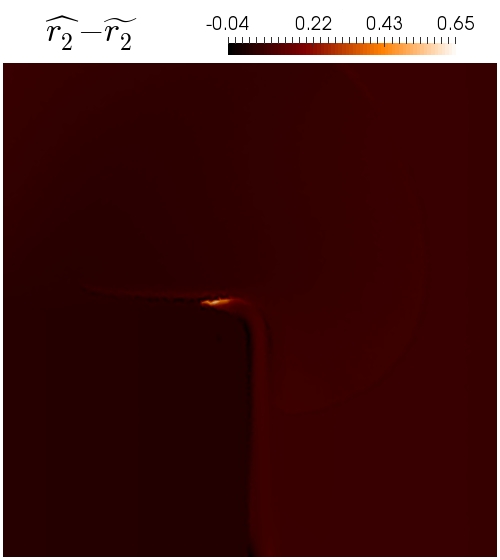}}
\subfigure[]{\includegraphics[width=0.325\textwidth]{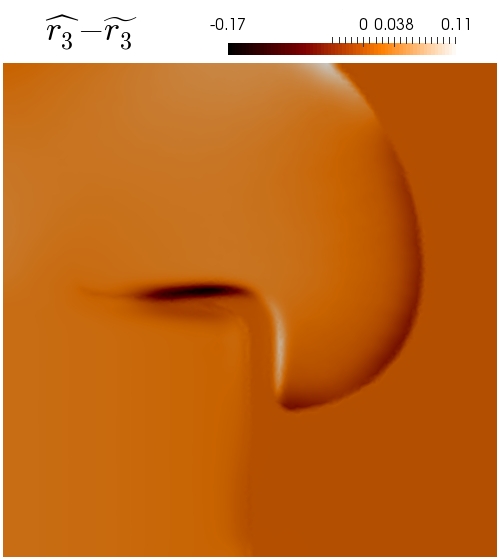}}
\end{center}
\vspace{-4mm}
\caption{Differences in ionic quantities from varying SAD parameter $D_2$ at $t=444\,\text{ms}$. Quantities $\widehat{r_i}$ indicate the profiles with $D_2=$7.5E-3, and $\widetilde{r_i}$ the profiles associated with $D_2=$1.0E-5.}\label{fig:diff-ionic}
\end{figure*}

\begin{figure*}
\begin{center}
\subfigure[$t=380\,\text{ms}$ ]{\includegraphics[height=0.35\textwidth]{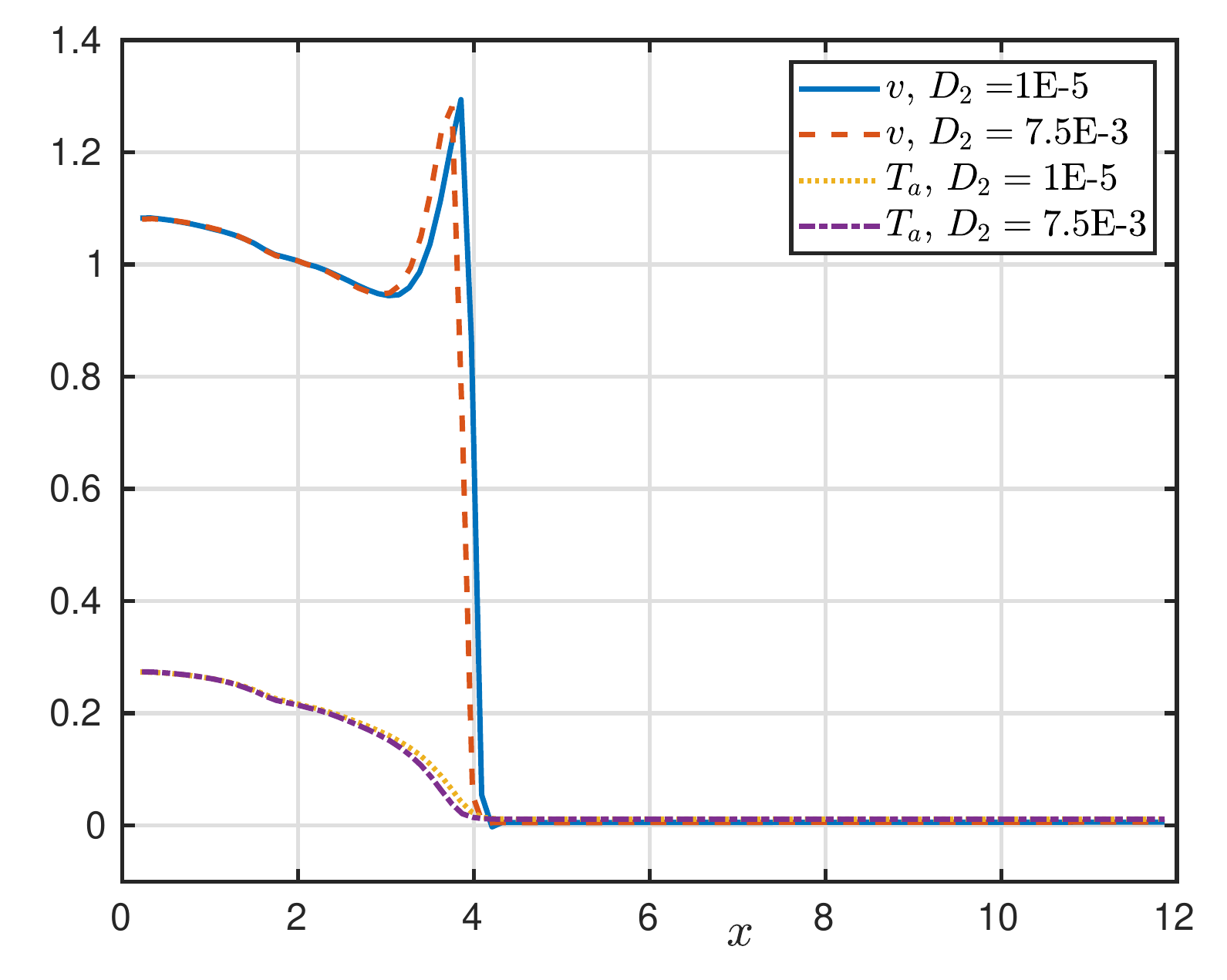}}
\subfigure[$t=456\,\text{ms}$ ]{\includegraphics[height=0.35\textwidth]{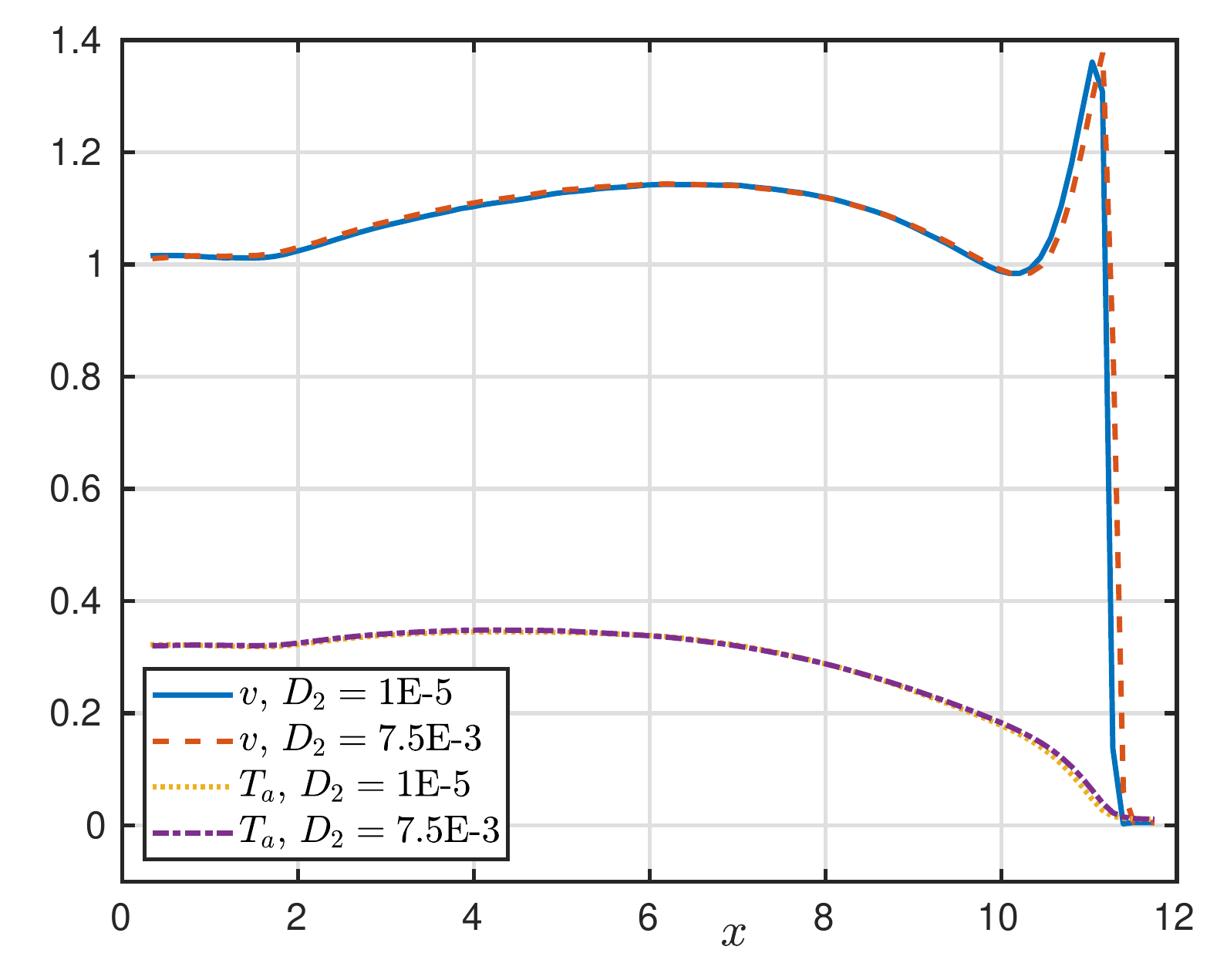}}\\[2ex]
\includegraphics[width=0.2\textwidth]{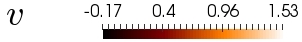}\\
\subfigure[$D_2=$1.0E-5]{\includegraphics[height=0.35\textwidth]{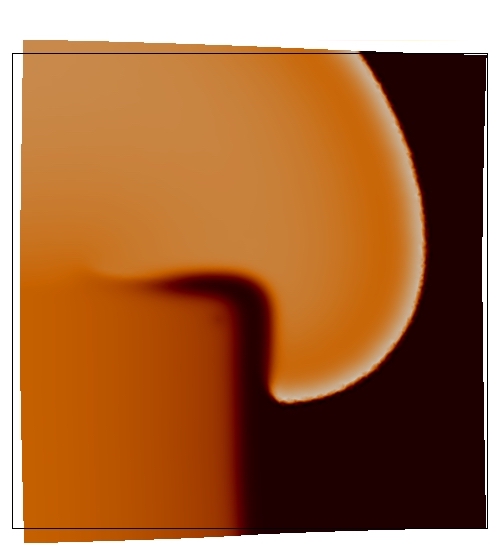}}
\subfigure[$D_2=$7.5E-3]{\includegraphics[height=0.35\textwidth]{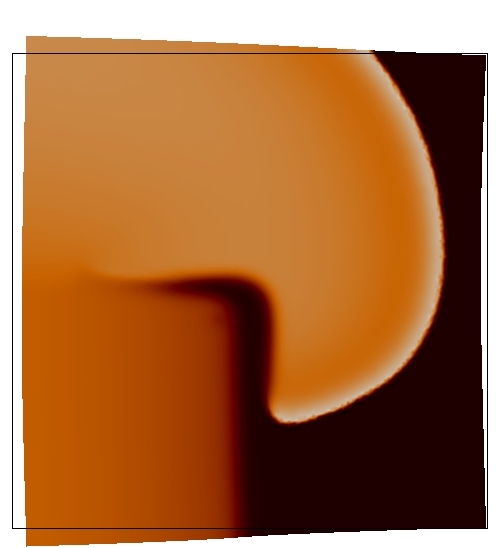}}
\subfigure[$D_2=$5.0E-2]{\includegraphics[height=0.35\textwidth]{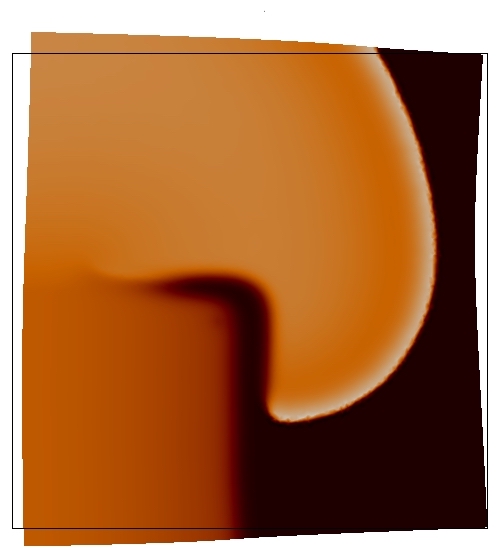}}
\end{center}
\vspace{-4mm}
\caption{(a,b): Propagation of action potential $v$ and active tension $T_a$, measured by taking the profile over a horizontal line segment crossing the upper half of the tissue at $y=7\,\text{mm}$. Comparison is provided for two different values of $D_2$. (c,d,e): Effect of $D_2$ 
on the potential wave at $t=444\,\text{ms}$ in the viscoelastic case.}\label{fig:diff-1d}
\end{figure*}

Figure~\ref{fig:diff-ionic} shows the differences 
in the ionic quantities between simulations with a very small contribution of SAD ($D_2 =$ 1 E-5\,mm$^2/$s) 
 and a more prominent, but still mild SAD contribution ($D_2 = $7.5E-3\,mm$^2/$s).  
 The snapshots correspond to the time $t = 444$ ms, when the spiral tip has not yet formed. 
A closer inspection suggests that these contrasts were due to a difference in conduction velocity 
 induced by SAD. In Figure \ref{fig:diff-1d}(a,b), we see that conduction velocity was higher for larger values of $D_2$ (meaning a larger SAD contribution). When the wave first emerged, the peak action potential was more advanced for the case of reduced $D_2$, but the large $D_2$ peak eventually caught up to and surpassed it, which is a phenomenon also observed in the active tension curves. The ionic quantities followed the same trend. Indeed, an analysis similar to that which produced Figure \ref{fig:diff-1d}(a,b) revealed that the overall profiles of the ionic quantities were highly similar between the two cases compared in Figure~\ref{fig:diff-ionic}, but differed in the speed at which they are transported through the tissue.

We also remark that the effect of changing conduction velocities was not spatially consistent. SAD increases conduction velocity in the fibre (horizontal) direction, but actually decreased conduction velocity in the vertical and diagonal directions. This resulted in a noteworthy effect on the growth of the spiral wave. Figures~\ref{fig:diff-1d}(c,d,e) show a comparison of the spiral wave in the viscoelastic case for three different values of $D_2$. The upper right area of the spiral is {slightly more vertical} in the simulation with a larger value of $D_2$ {than in the other cases}, suggesting that propagation of the voltage  {is somewhat hindered in the fibre} direction. {We also observe a slightly more pronounced deformation of the right side of the domain due to the two-way coupling between tissue motion and electrophysiology}. A similar effect was seen in the {hyperelastic} case. 

As in other studies, here we observe that conduction velocity is sensitive to spatio-temporal discretisation. In Table~\ref{table:CVconv}, we include the results of a simple convergence test for conduction velocity, similar to the benchmark test conducted in \cite{ruiz18}. We calculated the horizontal propagation of the action potential using different time steps and mesh refinements. Differently to  the case of nonlinear diffusion without SAD from \cite{ruiz18}, the experiment reveals that lower resolutions produce larger conduction velocities than the physiological values. 
This test also confirms that with our time step and mesh resolution ($0.1\,\text{ms}$, and above $200,000$ DoF, respectively), conduction velocity is in the expected physiological range, whereas larger time steps will systematically fail to capture the dynamics 
of the ionic model.

\begin{table}[t!]
\setlength{\tabcolsep}{3.5pt}
\rowcolors{1}{gray!30}{gray!10}
\begin{center}
\begin{tabular}{|l|c|ccccc|}
\hline
\multicolumn{7}{|c|}{Convergence of Conduction Velocity, mm/ms}\\
\hline
 {DoF}  &  h(mm)    &   $\Delta t=$ & $0.3\,\text{ms}$ & $0.1\,\text{ms}$ &  $0.05\,\text{ms}$ &  $0.01\,\text{ms}$  \\
\hline
27038 & 0.3817  && 0.1130	 & 0.1032 &	0.1015&	0.0994 \\
108576 & 0.1909 && 0.0754 & 0.0705 &	0.0654 &	0.0637 \\
170919 & 0.1527 && 0.0733 & 0.0657 & 0.0632& 0.0620 \\
246456 & 0.1273 && 0.0701 & 0.0632	 & 0.0601& 0.0589\\
554960 & 0.0849 && 0.0649 & 0.0553 & 0.0551& 0.0550 \\
1204362 & 0.0768 && 0.0610 & 0.0552 & 0.0550& 0.0547 \\
\hline
\end{tabular}
\caption{Convergence of conduction velocity with respect to temporal and spatial discretisation.}\label{table:CVconv}
\end{center}
\end{table}

\begin{figure*}[t!]
\begin{center}
\includegraphics[width=0.95\textwidth]{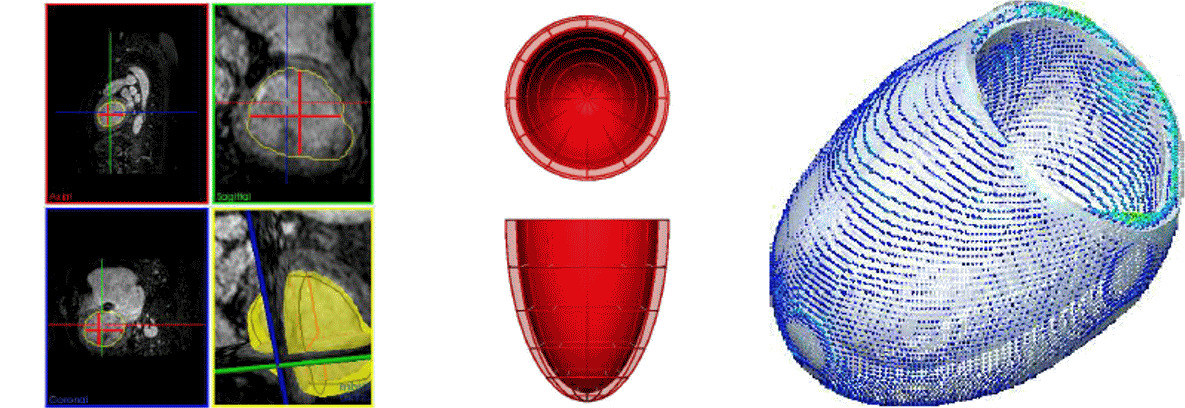}\\
\vspace{-4mm}
\subfigure[]{}\qquad\qquad\qquad\qquad\qquad\qquad\qquad\qquad\qquad
\subfigure[]{}\qquad\qquad\qquad\qquad\qquad\qquad\qquad\qquad\qquad
\subfigure[]{}
\end{center}
\vspace{-4mm}
\caption{Segmentation and mesh personalisation process from \cite{warriner18,lamata18}. Semi-automatic segmentation by 3D extrapolation (yellow surface and contours) of 2D segmentation contours (red contours and projections) (a); surface mesh template (b); and resulting mesh (white surface) overlaid with the segmentation surface colour coded by the distance between them (jet colour map, from 0 mm in blue to 1 mm in red) (c). Used with permission.}\label{fig:segmentation}
\end{figure*}

\subsection{Scroll waves on mono-ventricular geometries}\label{sec:scroll}
For the ventricular geometries, we test both the active strain and active stress formulations.
We start from patient-specific left ventricular geometries (available from 
\cite{warriner18,lamata18}) and rescale them using approximately the same dimensions as the idealised ventricles studied 
in \cite{ruiz18}. 
The segmentation process is outlined in Figure~\ref{fig:segmentation}. From there we 
define boundary labels and produce volumetric tetrahedral meshes of varying resolutions.  
The domain boundaries are set as sketched in Figure~\ref{fig:sketch}: The basal cut corresponds to 
$\partial\Omega_D$,  the epicardium to $\partial\Omega_R$, 
where the Robin boundary conditions \eqref{eq:robinBC} are defined with a spatially 
varying stiffness 
$$\eta(y) = \frac{1}{y_b-y_a}[\eta_a(y_b-y) + \eta_b(y-y_a)],$$
and the endocardium to $\partial\Omega_N$, where we set $p_N(t) = p_0\sin^2(\pi t)$, representing 
the variation of endocardial pressure. The constants 
$y_a,y_b$ are the vertical components of the apical and basal locations, 
and $\eta_a<\eta_b$ denotes the stiffness sought at the apex and base, respectively 
(assuming that the contact of the muscle with the aortic root is more resistant to traction than the more flexible pericardial sac and surrounding organs). In addition, since fibre and sheetlet fields for mono-ventricular geometries are not usually extracted from 
MRI data, we generate them using a mixed-form adaptation to the Laplace-Dirichlet rule-based method proposed in \cite{wong13,rossi14}.

\begin{figure*}[t!]
\begin{center}
\includegraphics[width=0.245\textwidth]{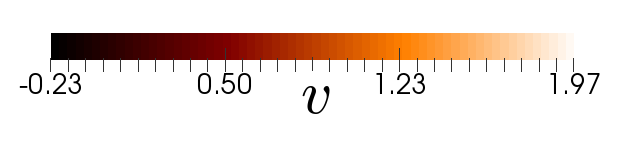}\\
\vspace{-2mm}
\subfigure[$t=400\,\text{ms}$]{\includegraphics[width=0.32\textwidth]{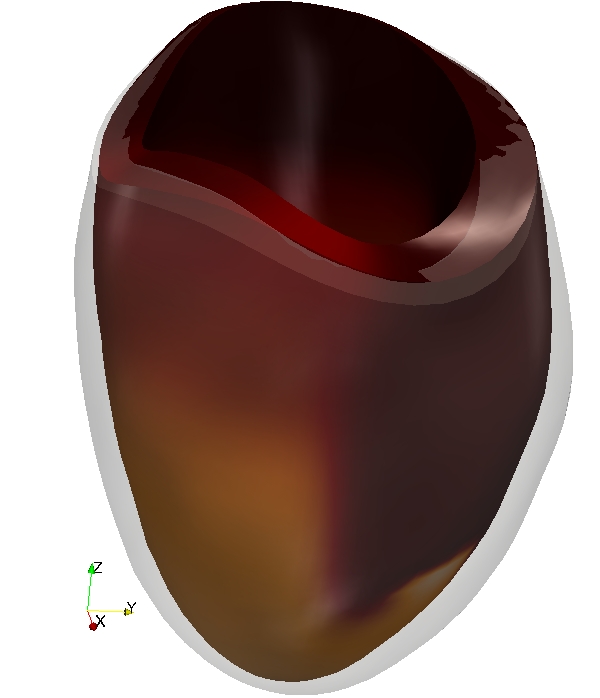}}
\subfigure[$t=500\,\text{ms}$]{\includegraphics[width=0.32\textwidth]{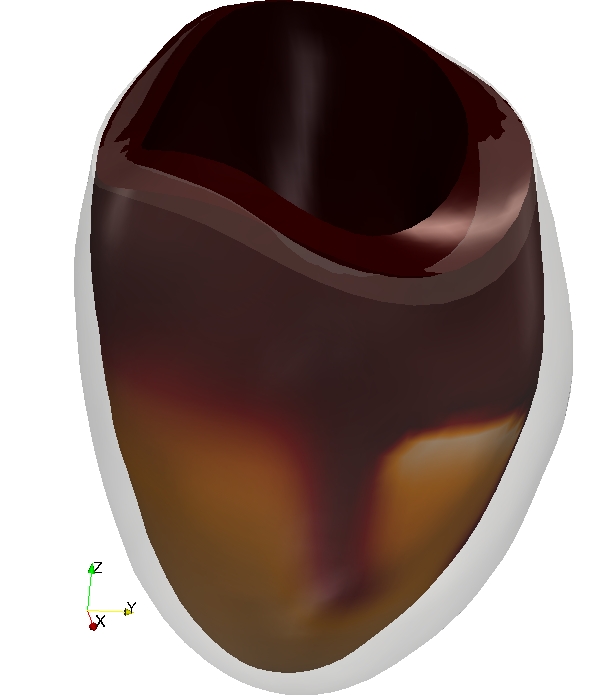}}
\subfigure[$t=600\,\text{ms}$]{\includegraphics[width=0.32\textwidth]{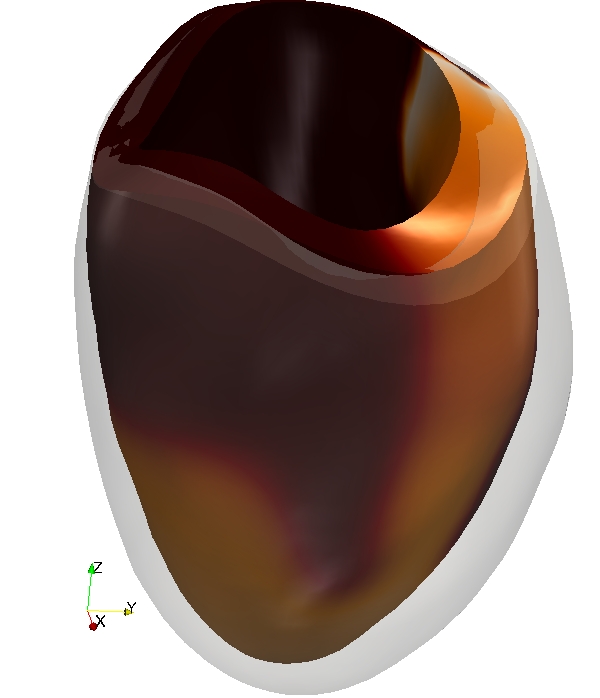}}
\end{center}
\vspace{-4mm}
\caption{Evolution of voltage after S2 stimulus (at $t=335\,\text{ms}$), showing formation of a scroll wave on a contracting ventricle, {using the active strain model}. The shadow of the undeformed ventricle geometry is shown for comparison.}\label{fig:3Dspiral}
\end{figure*}

\begin{figure*}[t!]
\begin{center}
\subfigure[]{\includegraphics[width=0.325\textwidth]{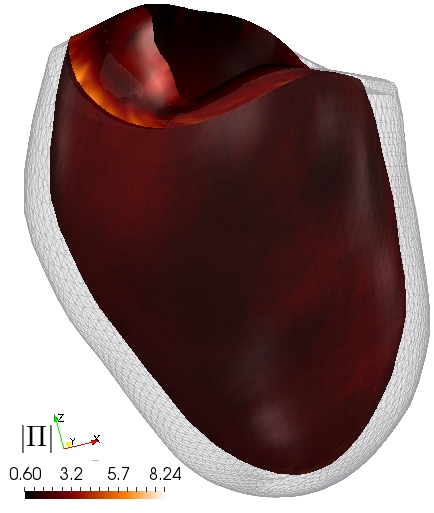}}
\subfigure[]{\includegraphics[width=0.325\textwidth]{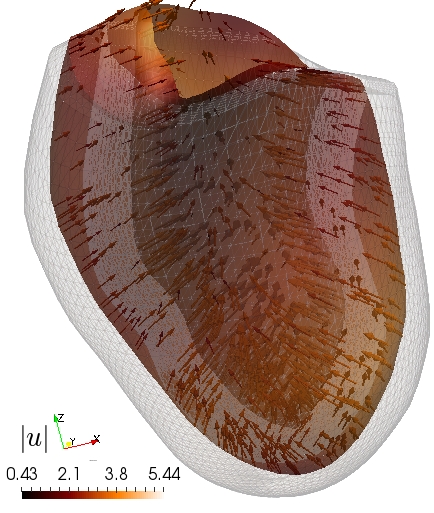}}
\subfigure[]{\includegraphics[width=0.325\textwidth]{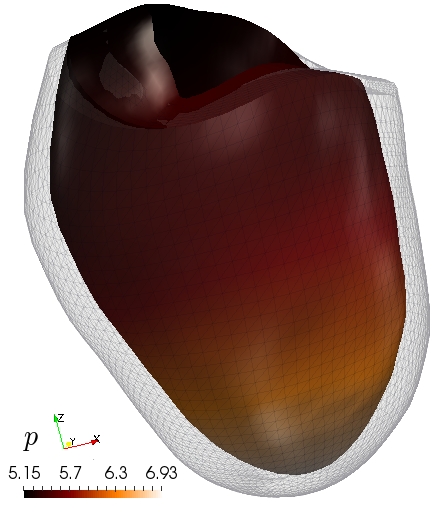}}\\
\subfigure[]{\includegraphics[width=0.325\textwidth]{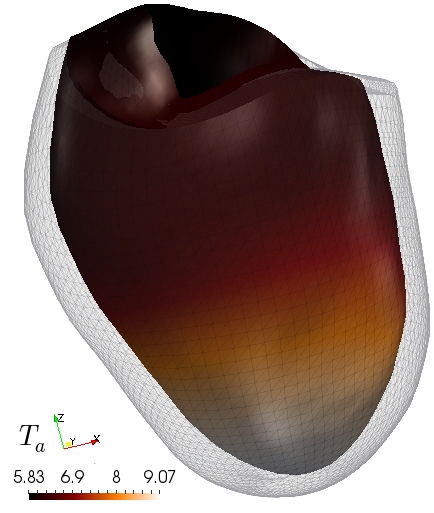}}
\subfigure[]{\includegraphics[width=0.325\textwidth]{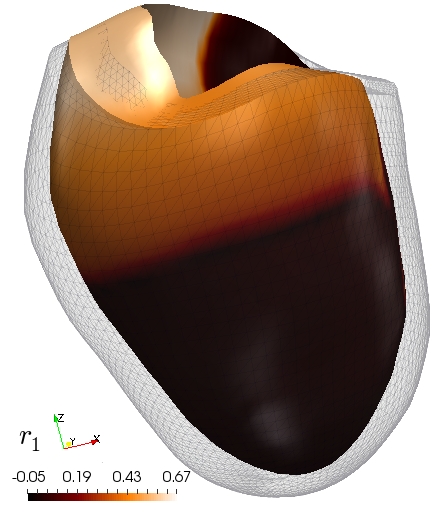}}
\subfigure[]{\includegraphics[width=0.325\textwidth]{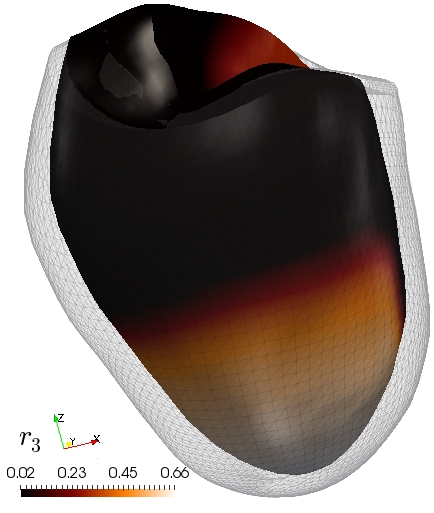}}
\end{center}
\vspace{-4mm}
\caption{Snapshot at $t=600\,\text{ms}$ of field variables plotted on the deformed domain and less opaque undeformed mesh. {Here we have also used the active strain approach}.}\label{fig:3Dfields}
\end{figure*}

After the S2 stimulus excites a group of cells in the lower left octant at $t = 335$ ms, a spiral wave forms and sweeps around both sides of the ventricle, the two sides merging at approximately $t = 415$ ms. Simultaneously, we see contraction of the apical region in the upwards direction, complemented by torsion and thickening of the ventricle wall. Figure~\ref{fig:3Dspiral} shows the propagation of the action potential on the deforming ventricle, with the original ventricle geometry shown with reduced opacity for comparison. The S2 stimulus occurs on the apex and the nascent scroll wave is not visible until the two arms of the wave interact. {For these tests the meshsize was approximately $h = 0.24\,$mm and the timestep $\Delta t = 0.1\,$ms. We have employed the lowest-order finite element method $l=0$ and the stabilisation parameter is taken as $\zeta_{\text{stab}} = 25$.}

\begin{figure*}
\begin{center}
\includegraphics[width=0.23\textwidth]{ex04-vbar}\\
\vspace{-4mm}
\subfigure[]{\includegraphics[width=0.325\textwidth]{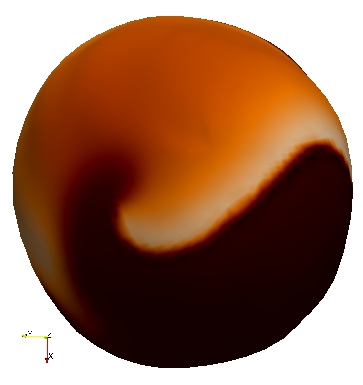}}
\subfigure[]{\includegraphics[width=0.325\textwidth]{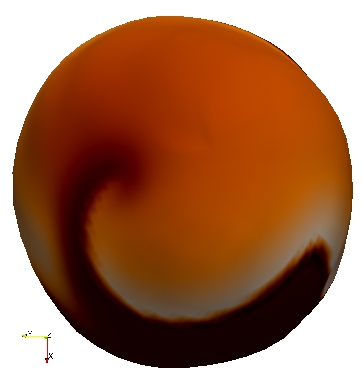}}
\subfigure[]{\includegraphics[width=0.325\textwidth]{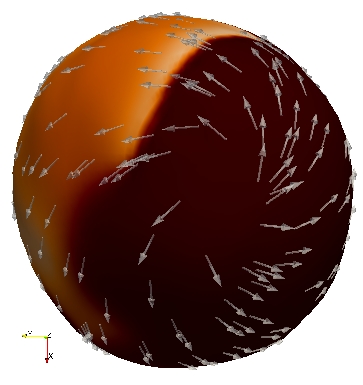}}
\end{center}
\vspace{-4mm}
\caption{Effect of SAD on spiral wave propagation, using the active strain formulation. Panels (a,b) show voltage and (c) shows the difference between the SAD and non-SAD cases $v_\text{SAD}-v_\text{non-SAD}$ (which has a different scale). The action potential wave using SAD moved along the fibre direction ahead of the non-SAD case.}\label{fig:SAD3DFibers}
\end{figure*}

In addition, and similarly to the 2D case, incorporating SAD impacted the propagation of the spiral wave anisotropically. In the fibre direction, SAD led to earlier advancement of the spiral. In the transverse direction, the non-SAD case advanced earlier. Figure \ref{fig:SAD3DFibers} shows the difference in voltage for the two cases (along with the actual voltage profile, for reference). 
The effect seen in the fibre direction (indicated by the white arrows) was not seen in the other directions. {For these tests 
we have used the active strain formulation, we have included viscoelastic effects, as well as inertial contributions}. 

\begin{figure*}
\begin{center}
\subfigure[]{\includegraphics[width=0.325\textwidth]{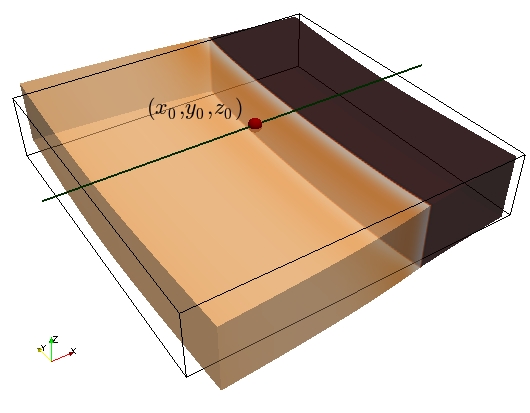}}
\subfigure[]{\includegraphics[width=0.325\textwidth]{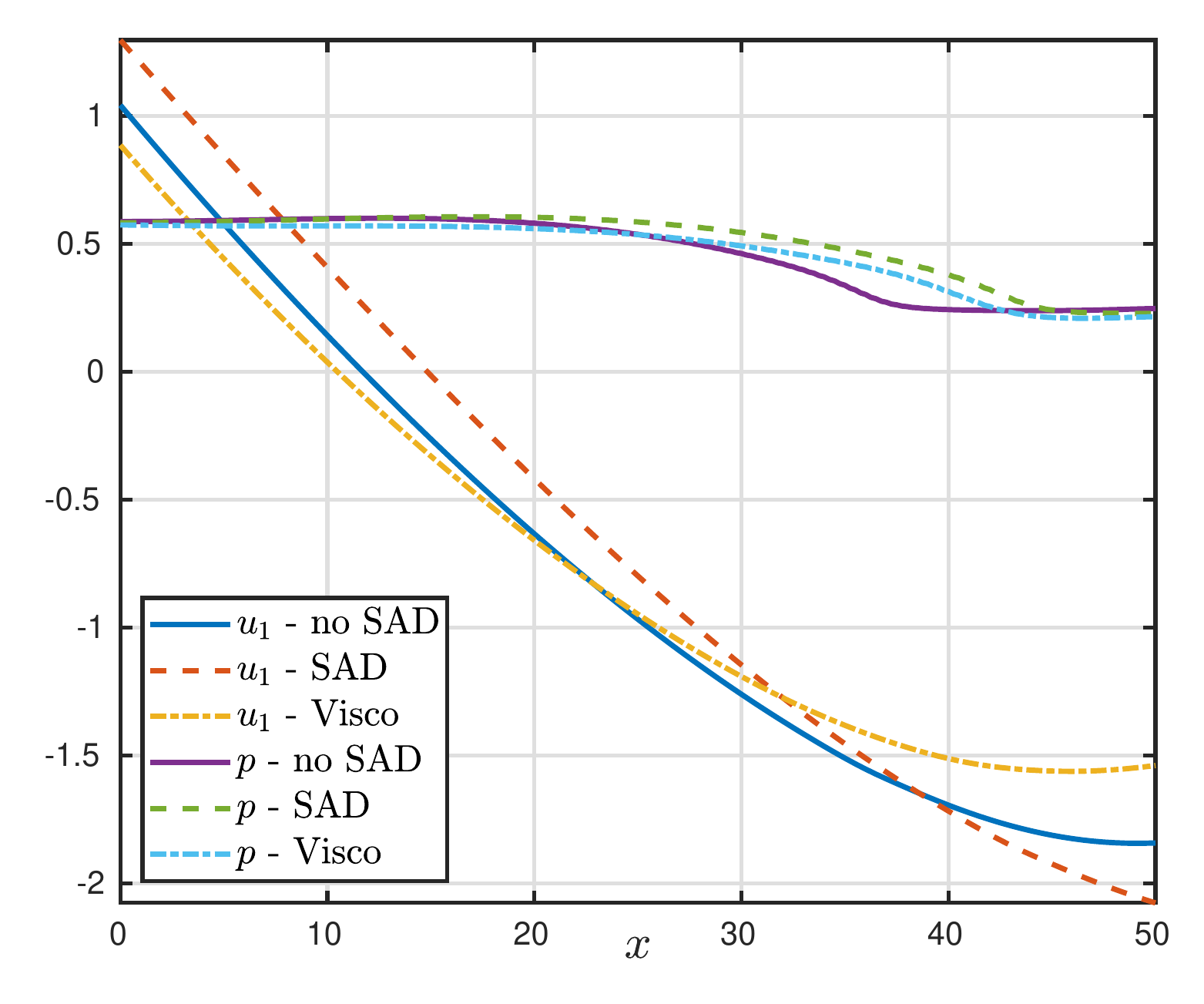}}
\subfigure[]{\includegraphics[width=0.325\textwidth]{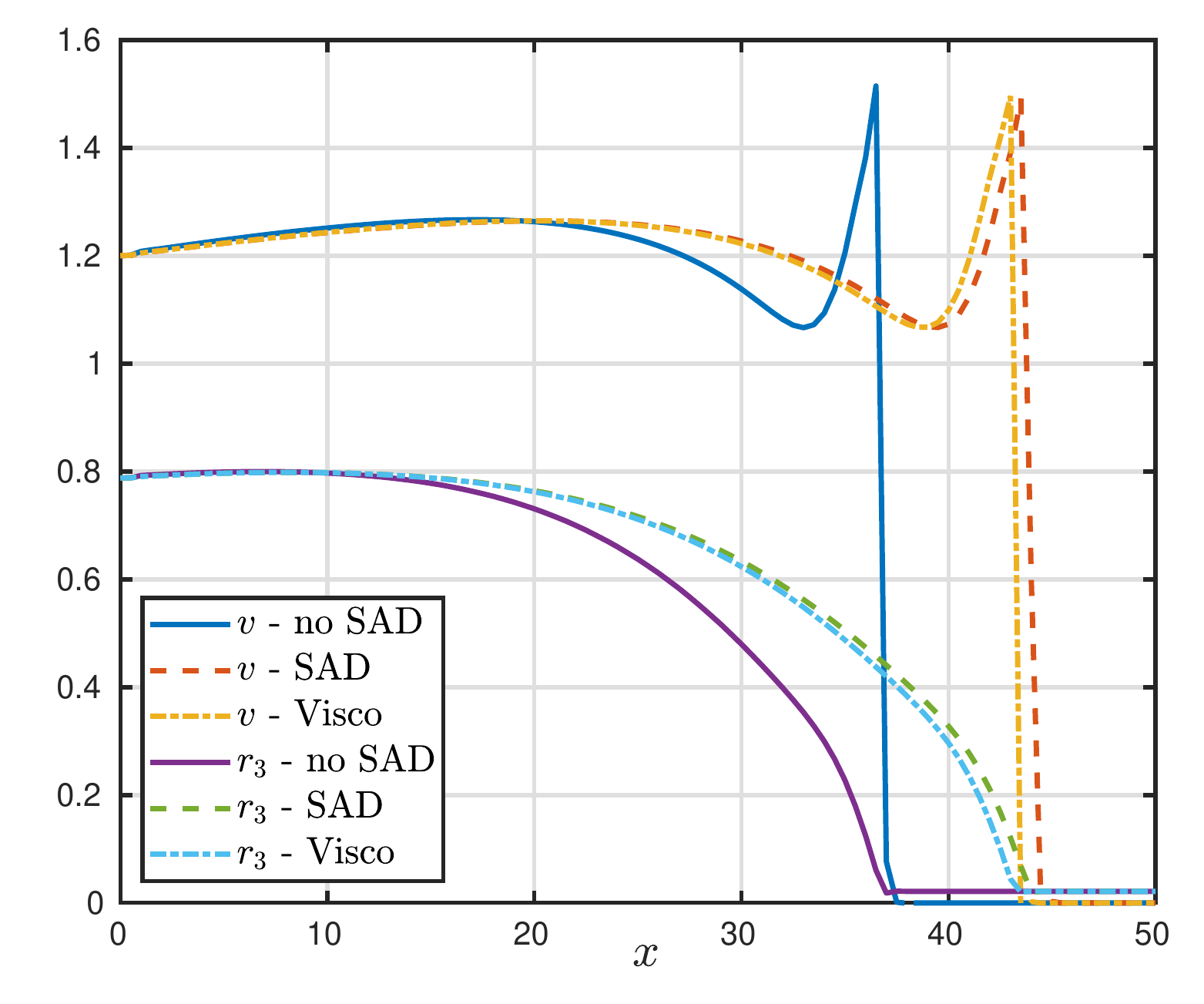}}\\
\subfigure[]{\includegraphics[width=0.325\textwidth]{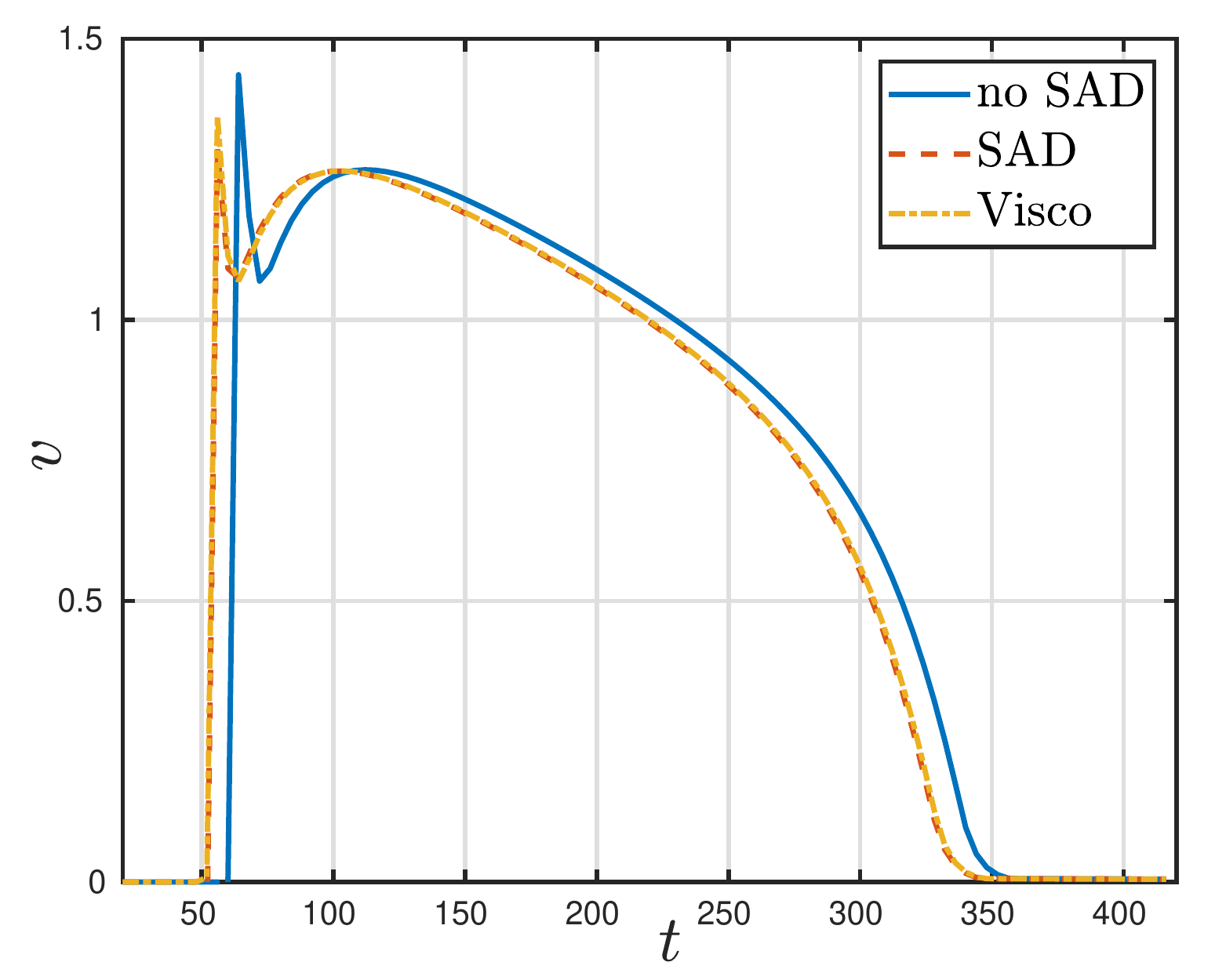}}
\subfigure[]{\includegraphics[width=0.325\textwidth]{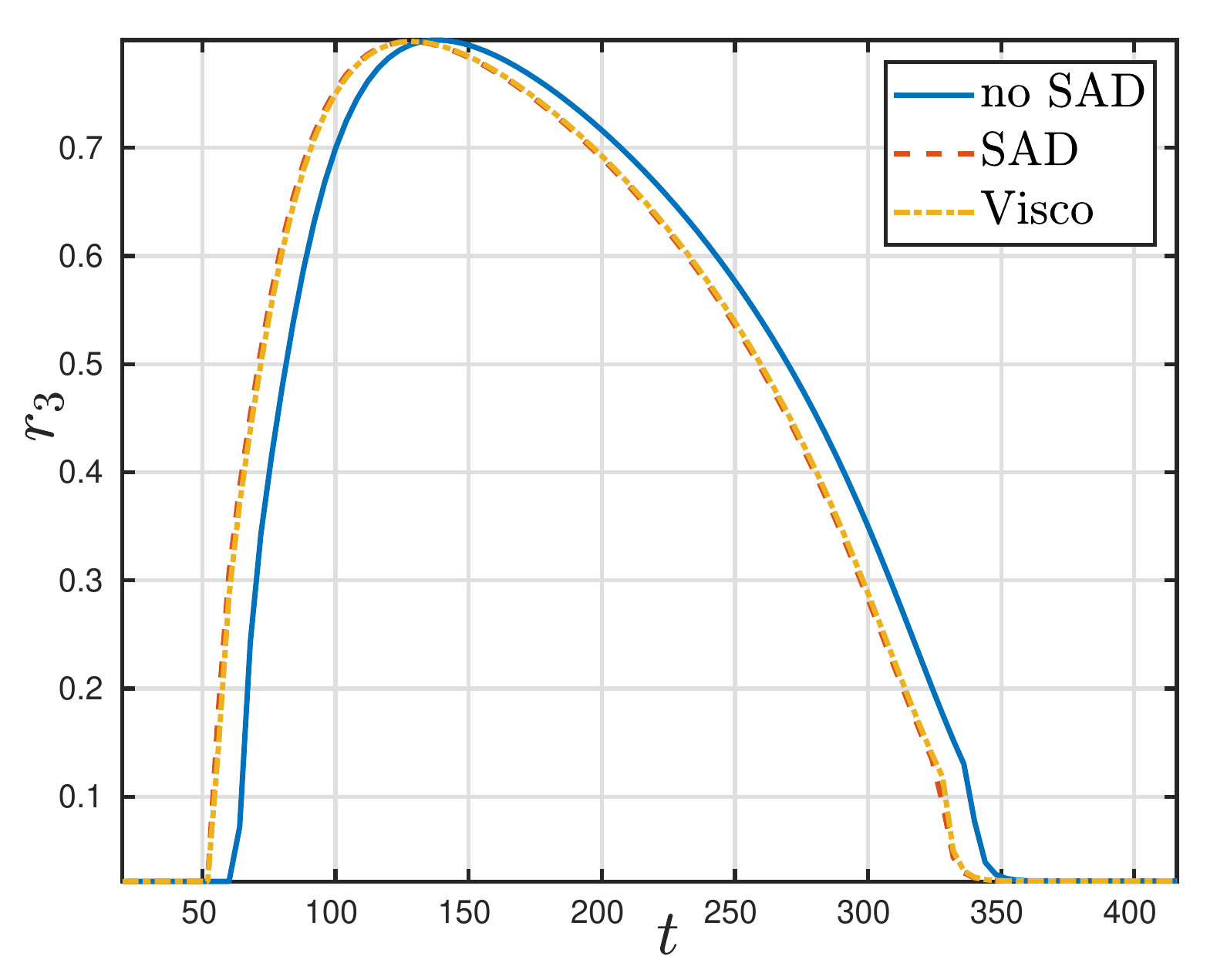}}
\subfigure[]{\includegraphics[width=0.325\textwidth]{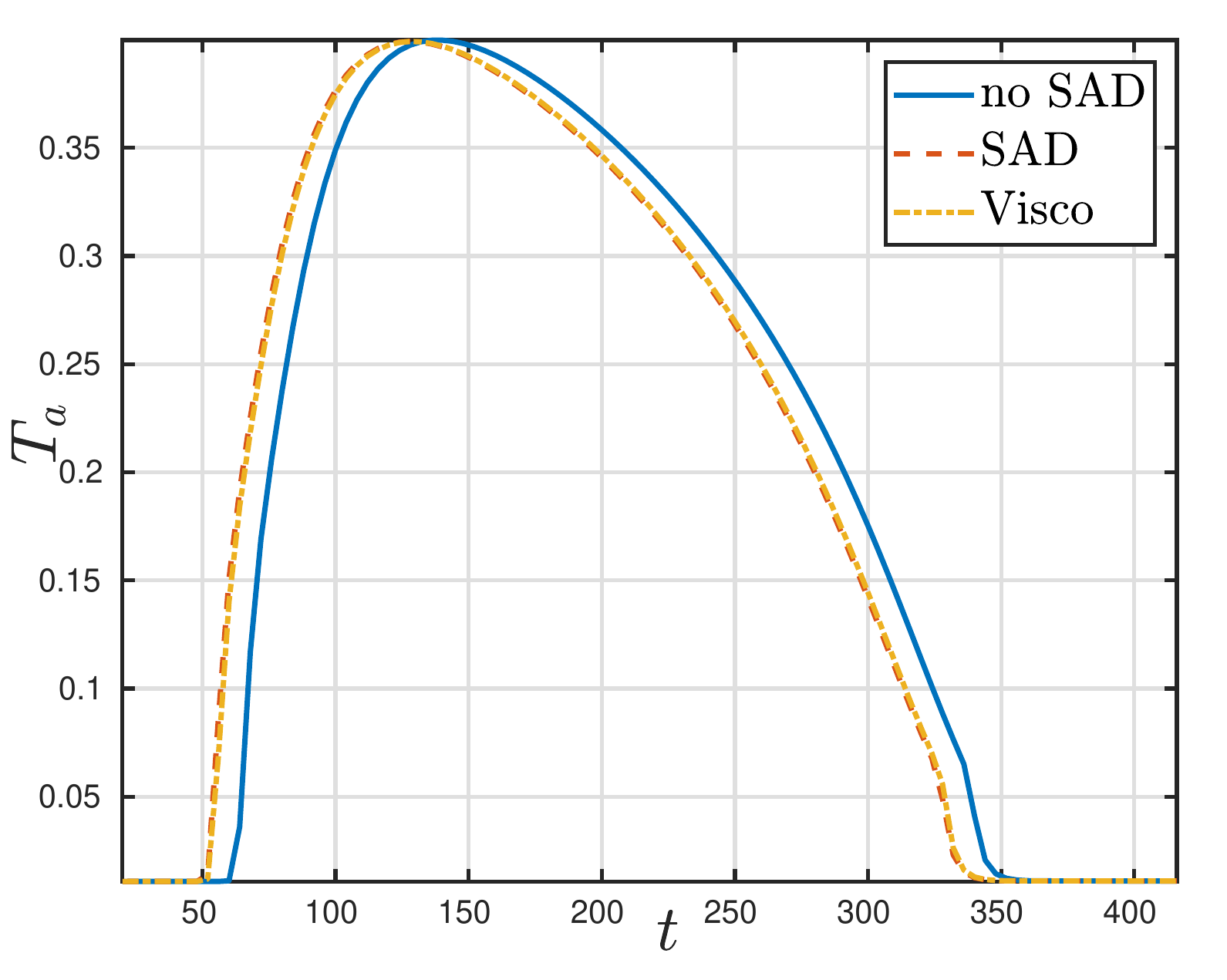}}\\
\subfigure[]{\includegraphics[width=0.325\textwidth]{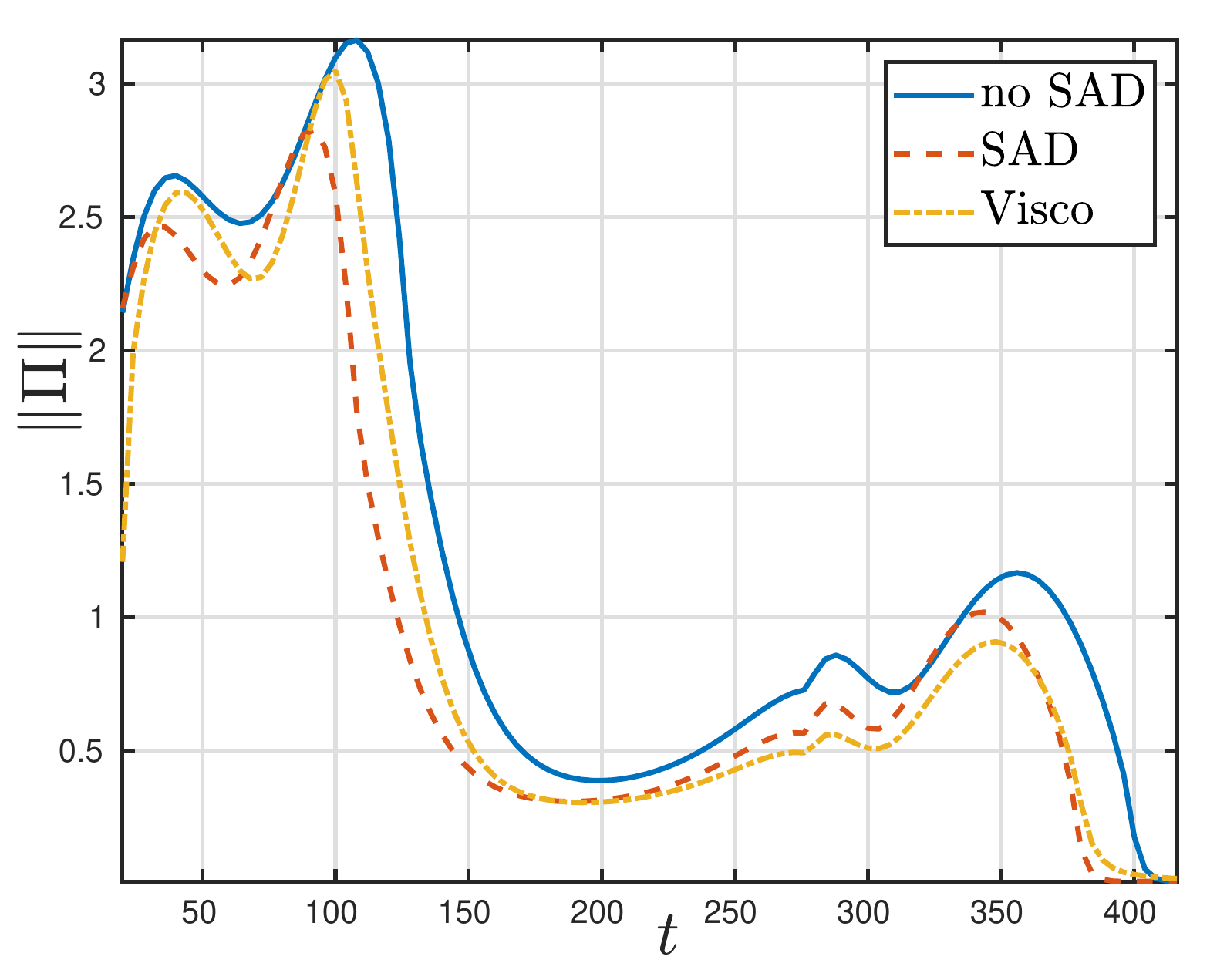}}
\subfigure[]{\includegraphics[width=0.325\textwidth]{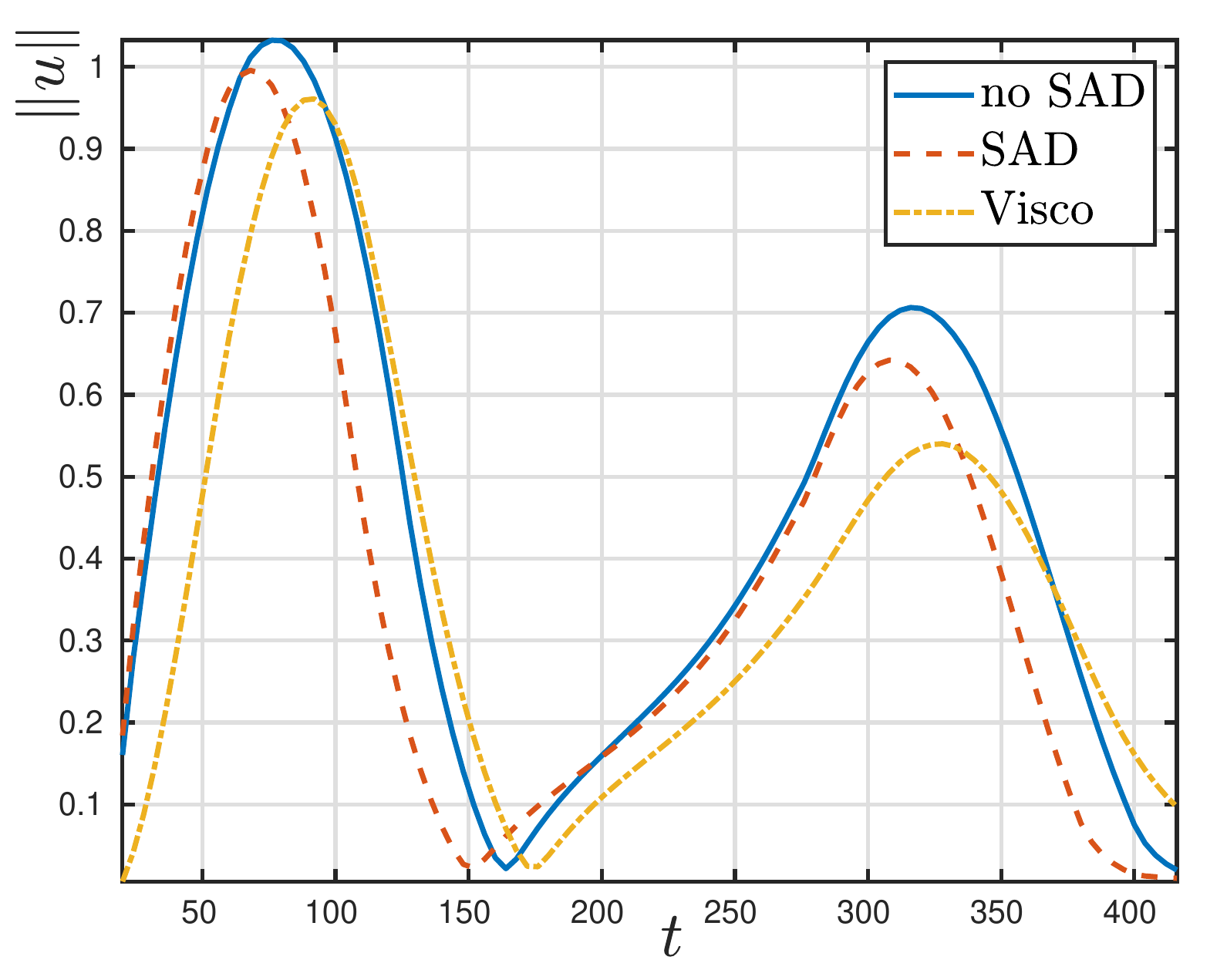}}
\subfigure[]{\includegraphics[width=0.325\textwidth]{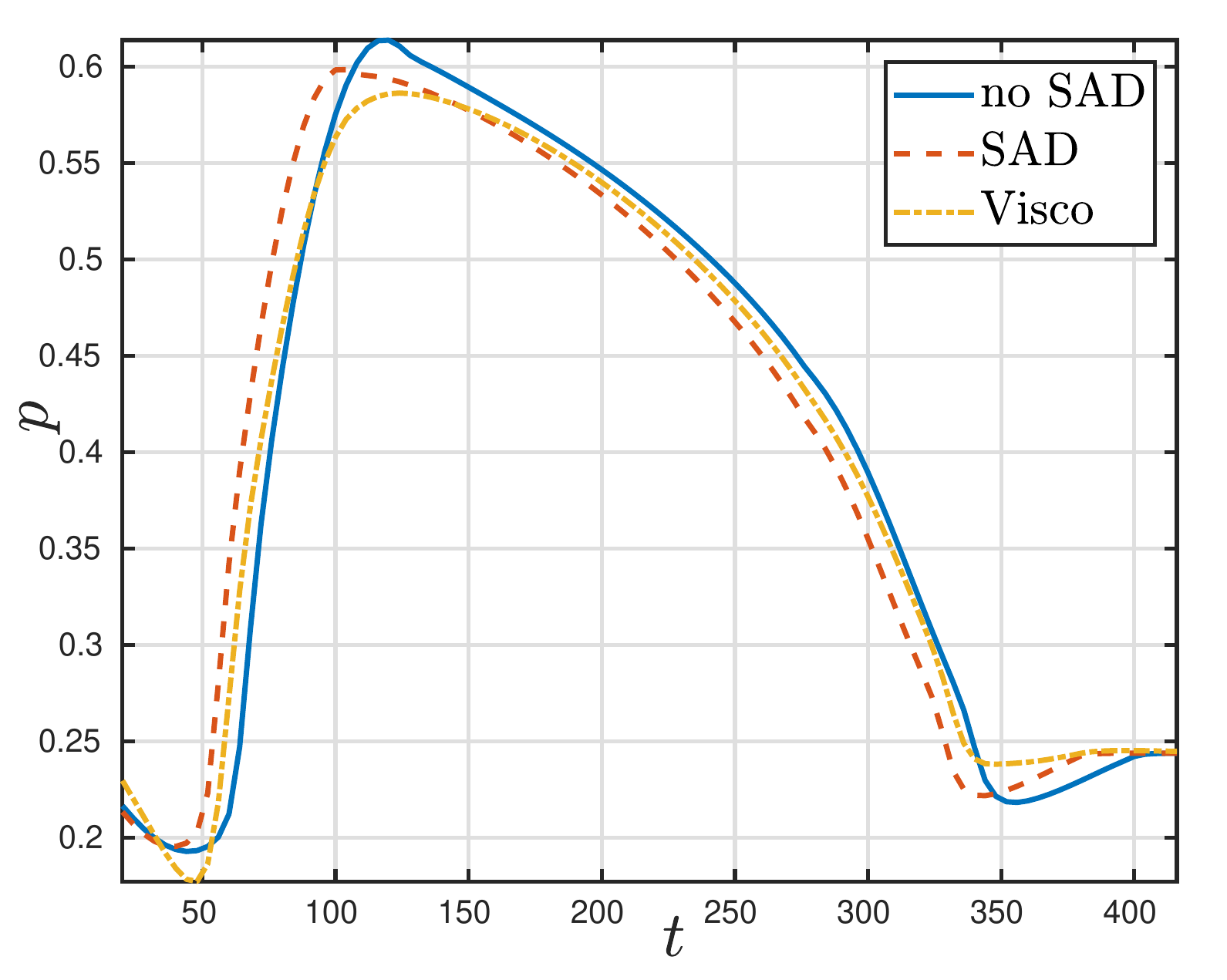}}
\end{center}
\vspace{-4mm}
\caption{Comparison of field variables between hyperelastic and viscoelastic cases on a line parallel to the $x$ axis (sketched in 
(a)) taken at $t=92\,\text{ms}$ (b,c); and point-wise evolution of field variables on the point $(x_0,y_0,z_0)$ (d-i) for the cases of 
hyperelasticity without SAD, with the baseline case of SAD but without viscous stresses, and the viscoelastic case 
(line, dashed, and dashed-dotted curves, respectively). {For these tests we have used the active stress formulation}.  }
\label{fig:visc-slab}
\end{figure*}

\subsection{Effects due to viscoelasticity} 
In order to quantify the discrepancies between hyperelastic and viscoelastic effects, we conduct a series of simulations using the coupled 
model on a 3D slab of dimensions $50\times 50\times 10$ mm$^3$ using a tetrahedral mesh of $h = 0.25$\,mm, {also setting $\Delta t = 0.1\,$ms, $\zeta_{\text{stab}} = 25$, and} 
$\fo=(1,0,0)^{\tt t}$, $\so =(0,1,0)^{\tt t}$. {These tests are conducted using the active stress formulation, and we consider inertial effects}. 
We apply a S1 stimulus 
on the face $x=0$ and after $t= 92\,\text{ms}$ the propagation front has reached the state shown in Figure~\ref{fig:visc-slab}(a), 
plotted on the deformed configuration (which was computed with a full electro-viscoelastic model). {The boundary conditions for 
the viscoelasticity are of Robin type everywhere on the boundary}. 
At that time, in panels (b,c) we depict snapshots of the approximate solutions obtained 
using the hyperelastic and viscoelastic models with their base-line parameter values as reported in Table~\ref{table:params}, and 
shown over a line segment crossing the tissue slab parallel to the $x-$axis. We 
show profiles of the mechanical entities ($x-$components of displacement and pressure), as well as potential and $r_3$. For reference, we also include the results 
obtained using a model without SAD contributions  (that is with $D_2=0$). 
We note that the curves produced without SAD 
are substantially lagged (as expected from the choice of diffusion parameters) with respect to the two other cases, that display no major discrepancies. The remaining 
panels in the figure show point-wise transients of the main mechanical and electrical fields measured on the point 
$(x_0,y_0,z_0)=(25,25,10)$. The evolution of the electric and activation fields remains very similar in all three cases; for instance 
the shape of the action potential is almost not modified after adding SAD or viscous contribution and for the other fields also 
very subtle differences are observed (the calcium concentration was slightly shifted to the left in the hyperelastic and viscoelastic cases). 
The changes are more pronounced 
in the Frobenius norm of the Kirchhoff stress, the displacement magnitude and the pressure 
(panels g,h,i). These computations suggest that viscous effects will result in a decreased displacement, stress, and pressure (similar 
conclusions were drawn in \cite{pandolfi17}, but not in the context of models for ventricular viscoelasticity).  
These discrepancies, however, are qualitatively small, and this observation was robust 
to every parameter combination that we tested, consistent spatially and in time. The application of  
a viscous model also had consequences related to performance. For instance, in the tests mentioned above, 
the average number of Newton iterations needed to reach convergence 
was systematically lower in the viscous case than in the hyperelastic case. This behaviour is expected as 
for simple viscoelastic models the tangent problem is essentially a rescaled version of the elastic stiffness, 
which contributes to improving the stability of the tangent problem.

\begin{figure*}
\begin{center}
\subfigure[]{\includegraphics[width=0.32\textwidth]{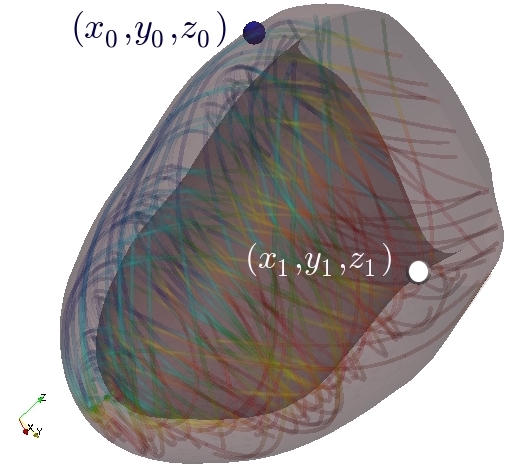}}
\subfigure[values on $(x_0,y_0,z_0)$]{\raisebox{0.2cm}{\includegraphics[width=0.32\textwidth]{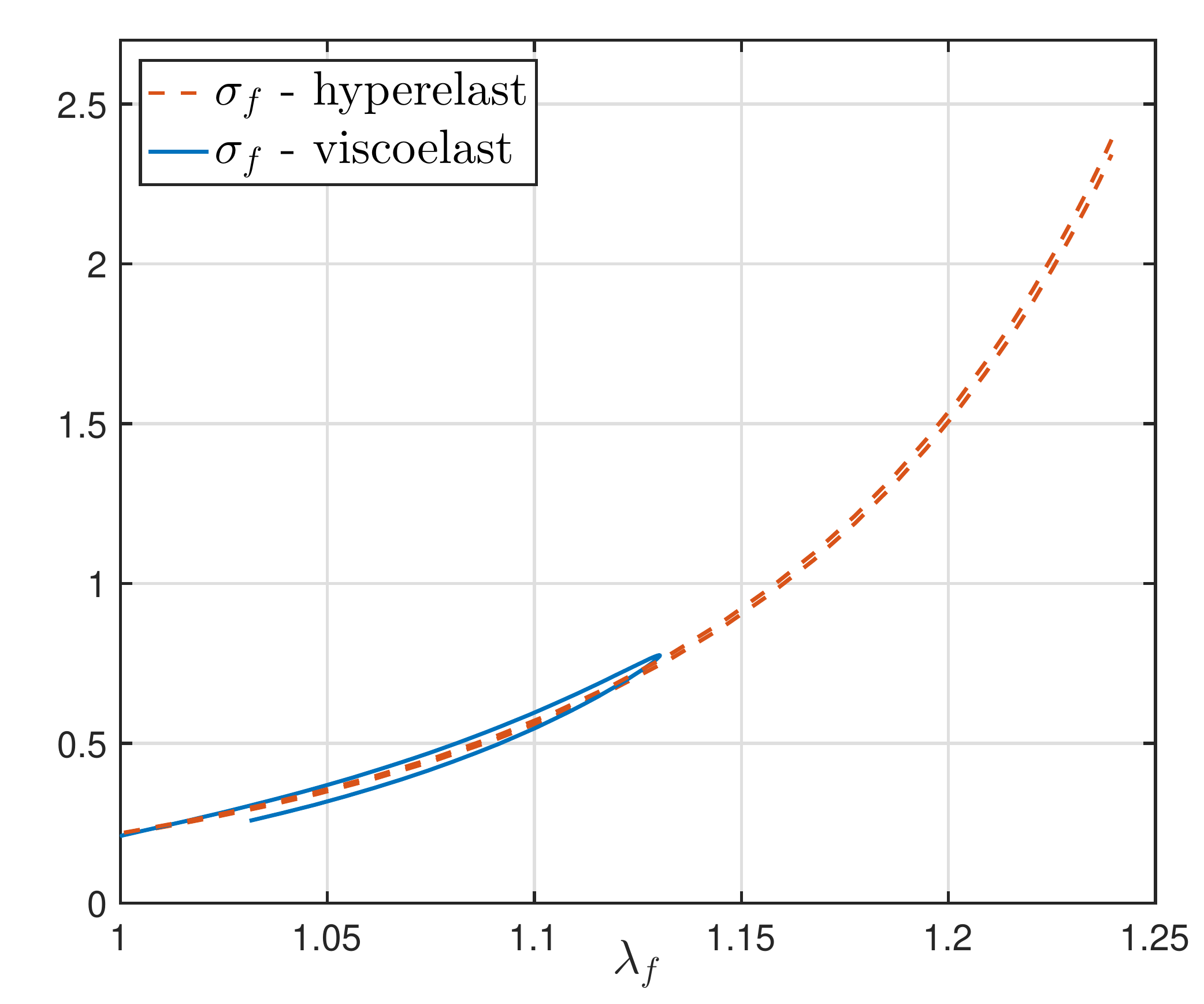}}}
\subfigure[values on $(x_1,y_1,z_1)$]{\raisebox{0.2cm}{\includegraphics[width=0.32\textwidth]{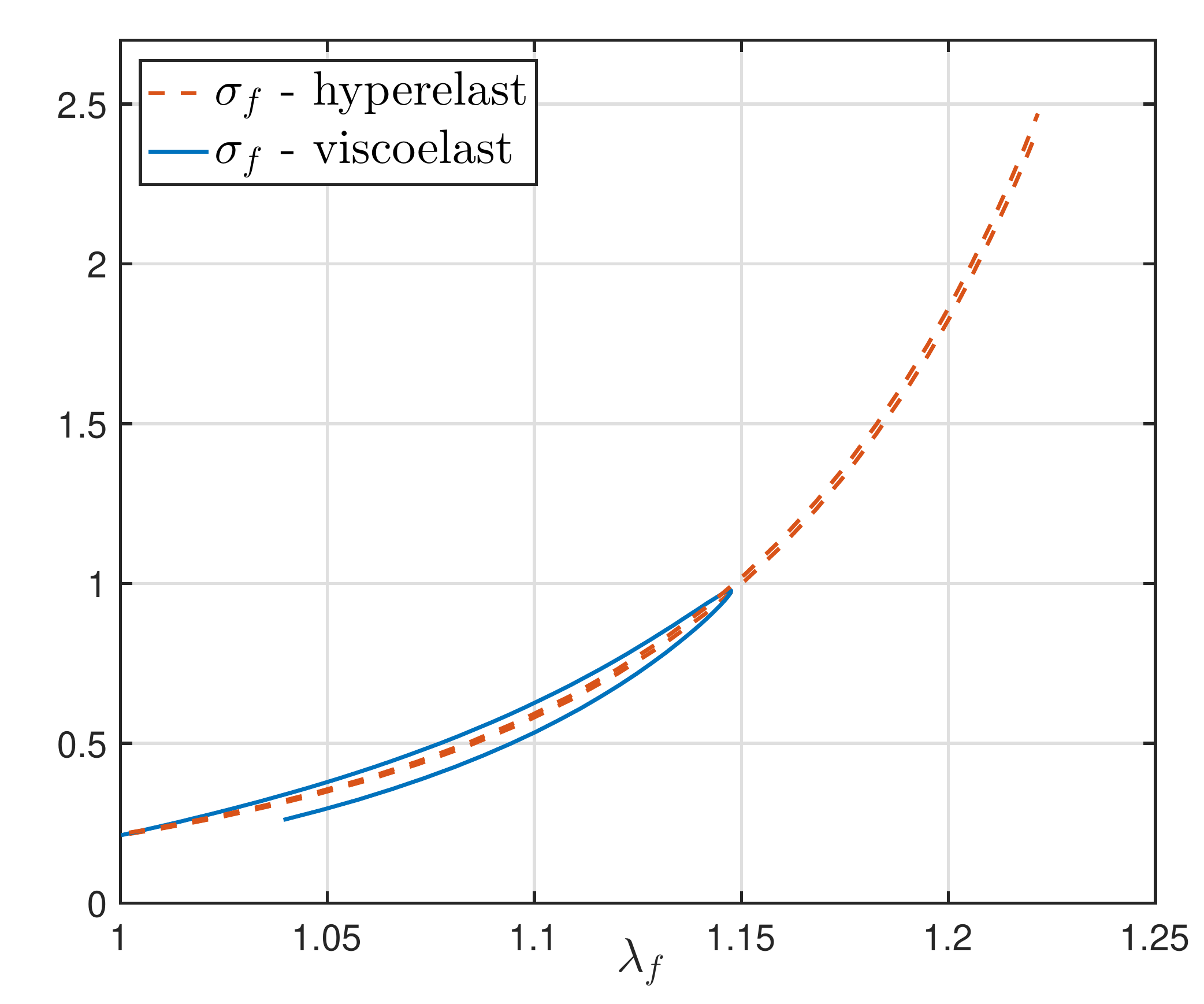}}}\\
\subfigure[values on $(x_0,y_0,z_0)$]{\includegraphics[width=0.475\textwidth]{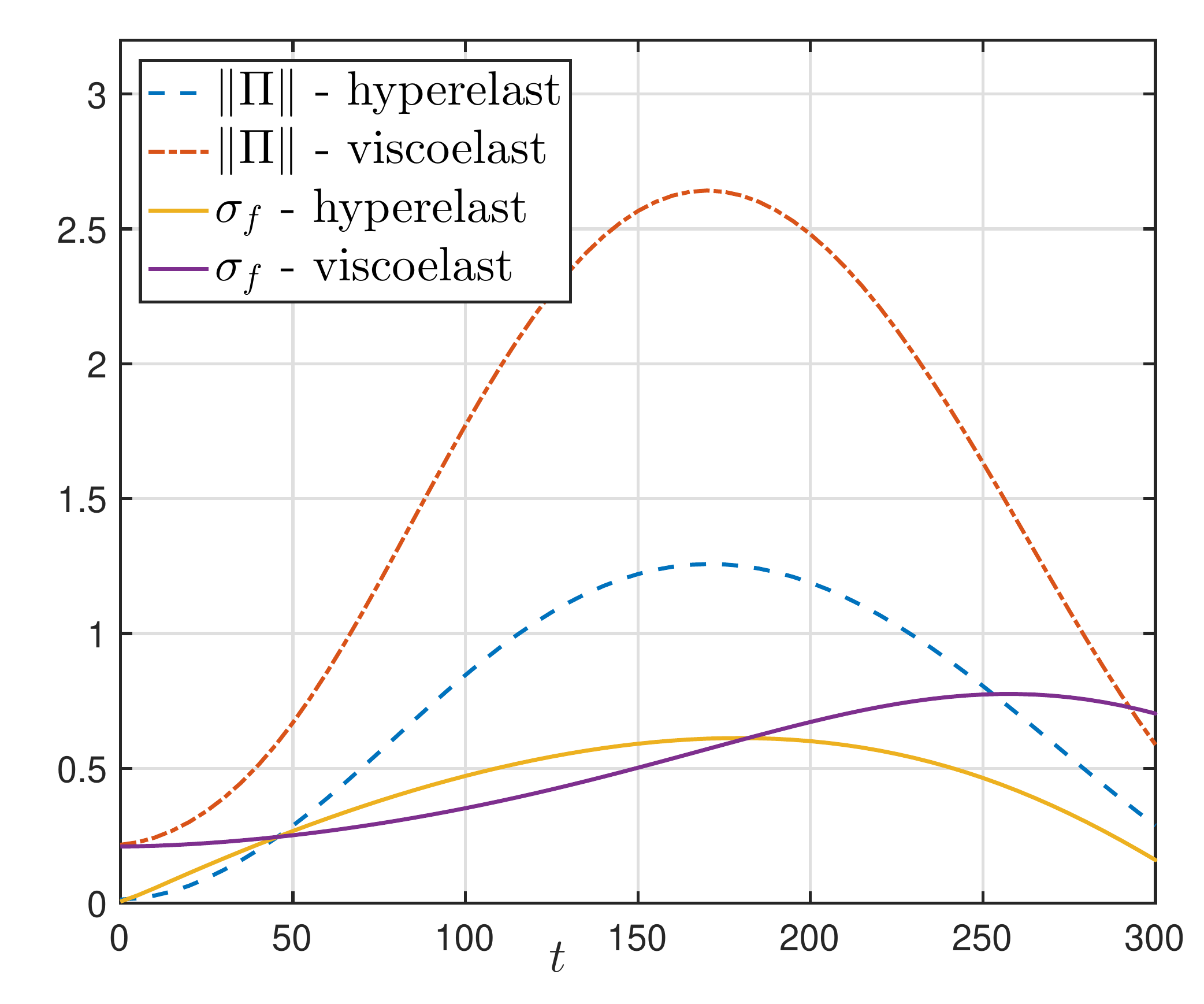}}
\subfigure[values on $(x_0,y_0,z_0)$]{\includegraphics[width=0.475\textwidth]{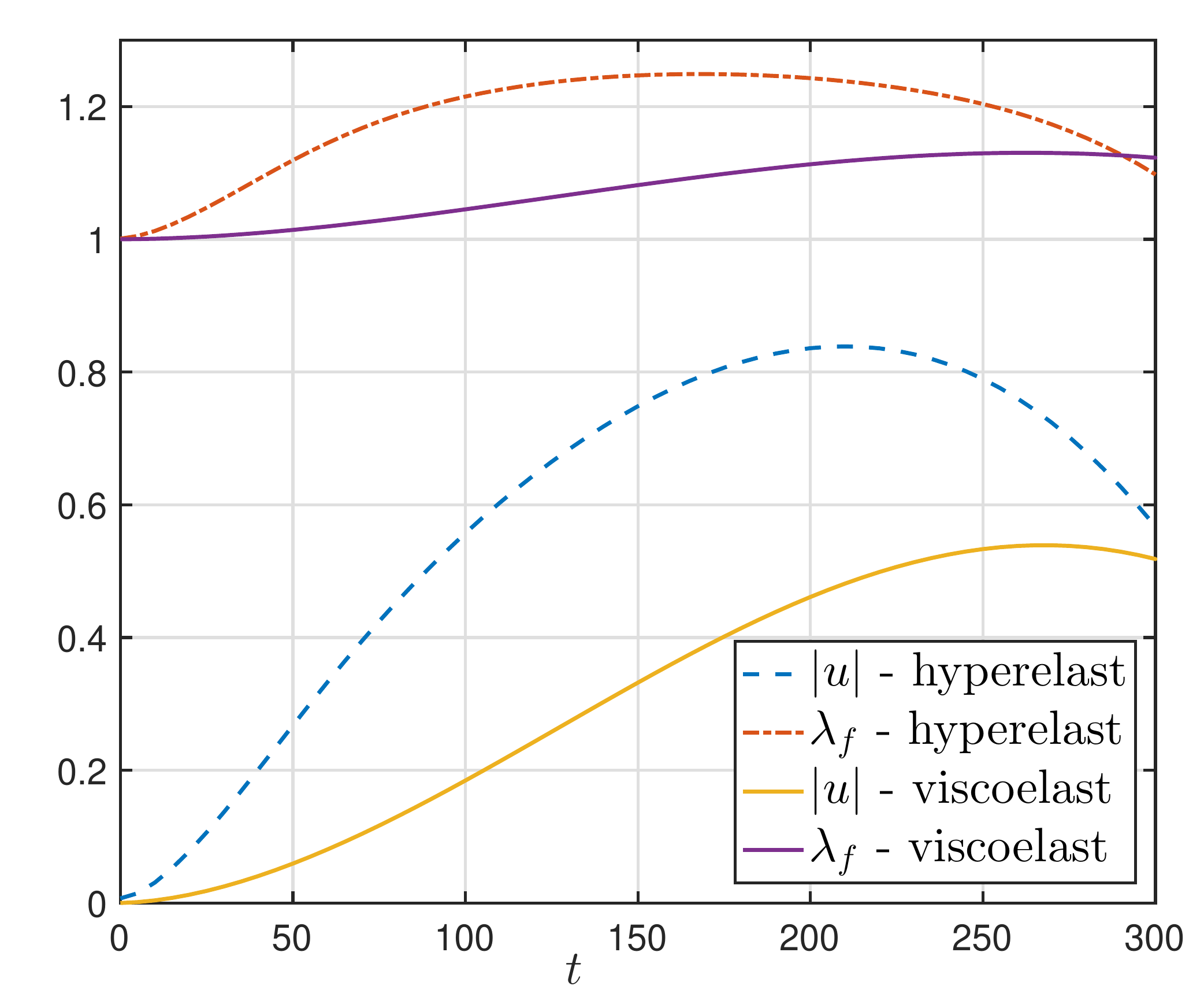}}\\
\subfigure[values on $(x_1,y_1,z_1)$]{\includegraphics[width=0.475\textwidth]{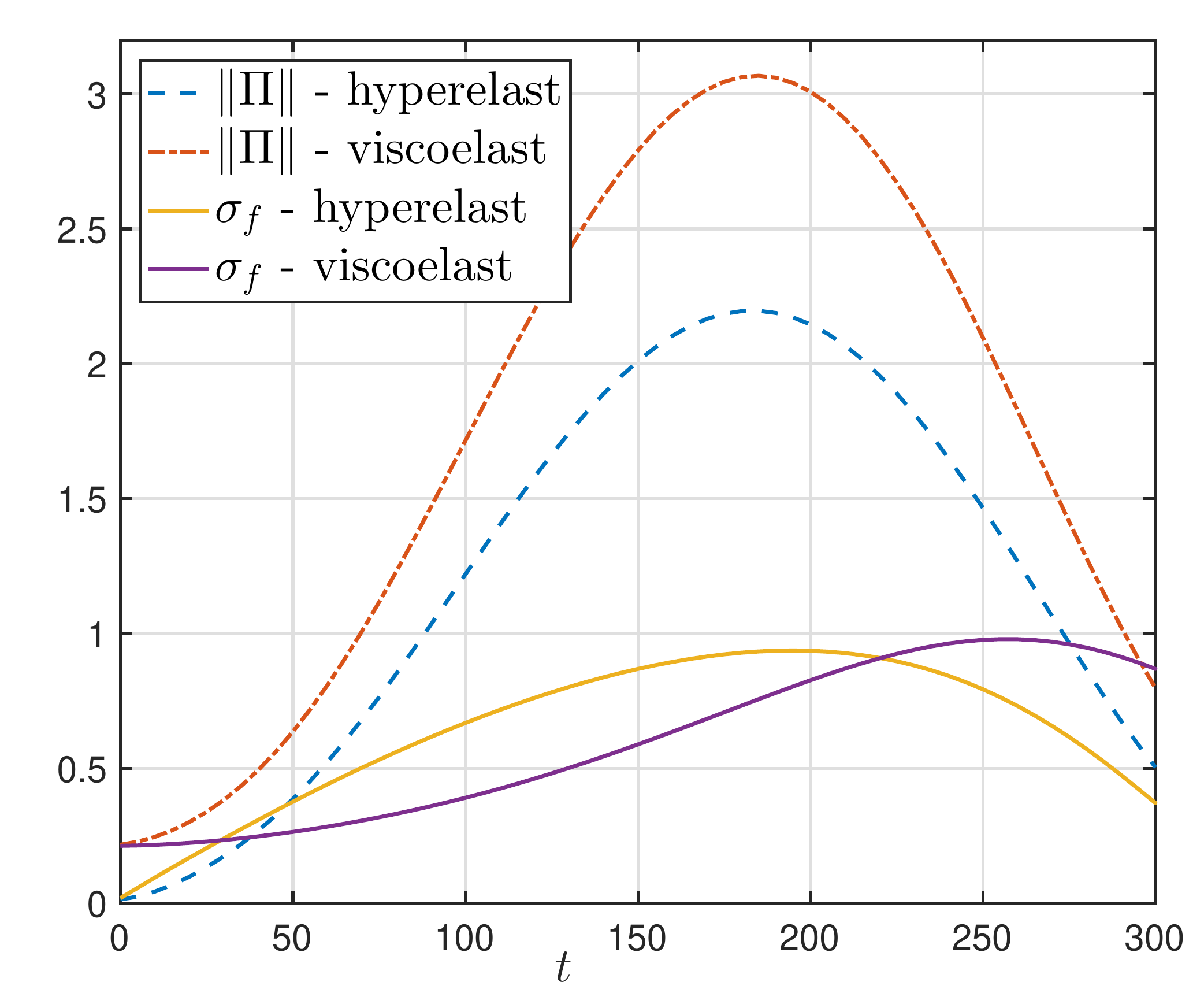}}
\subfigure[values on $(x_1,y_1,z_1)$]{\includegraphics[width=0.475\textwidth]{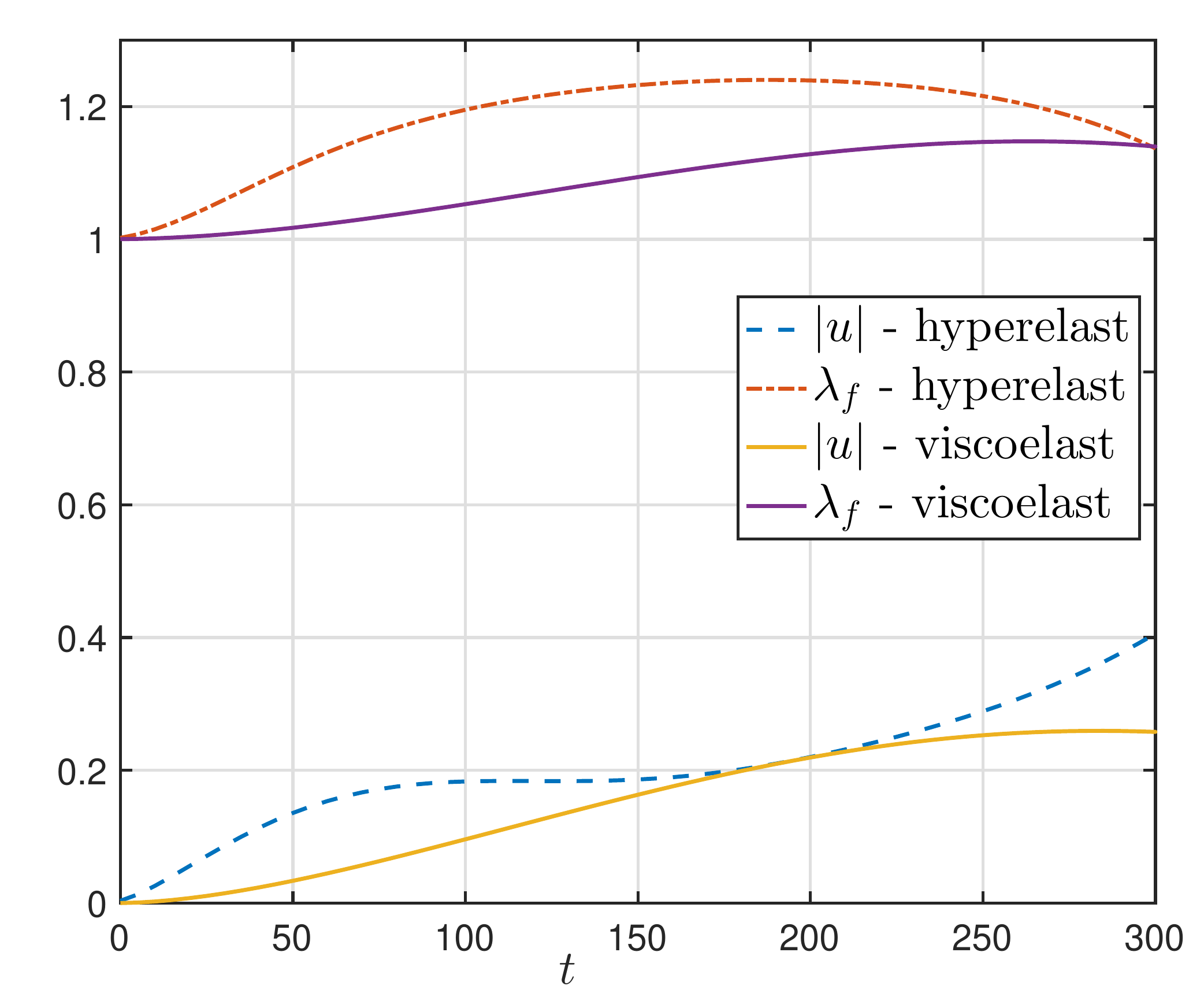}}
\end{center}
\vspace{-4mm}
\caption{Comparison between hyperelastic and viscoelastic deformation under passive inflation. True stress in the fibre direction 
$\sigma_f = \bF\fo \cdot (\bsigma \bF\fo)$, measured according to local stretch on 
two points on the epicardium (b) and endocardium (c) (points indicated in panel (a)). The plots in panels (d,e) show 
transients of mechanical outputs (Frobenius norm of the Kirchhoff stress, true stress on fibre direction, local stretch, and 
displacement magnitude) at the point $(x_0,y_0,z_0)$; and plots (f,g) display their 
counterparts in point $(x_1,y_1,z_1)$. {For these tests we have used only inertial effects and passive hyperelastic or viscoelastic contributions}.} 
\label{fig:visc-ventr}
\end{figure*}

We next proceed to investigate the effects of changing the viscosity parameters. The parameter $\beta$ from \eqref{eq:beta} exerted minimal influence over the observed dynamics. Even for the five orders of magnitude tested, from $\beta = 0.1$\,ms to $\beta = 10000$\,ms, the differences in displacement, voltage, and all other variables were of less than 0.1\%. This could be because of the low rates of change of deformation that we see in our simulations. We also tested values of $\delta$ across three orders of magnitude, from $\delta = 2.26$ to $\delta = 2260$ (in units N$/$cm$^2\cdot$ms). As expected, increasing this quantity, thereby increasing the viscoelastic contribution to the Cauchy stress, magnified the differences between the hyperelastic and viscoelastic cases (essentially magnifying the effects seen in Figure~\ref{fig:visc-slab}). Additional simulations (not reported here) also showed that higher values of $\delta$ not only reduced the magnitude of $\bPi$, $\bu$, $p$, but also smoothed their profiles, reducing the distances between peaks and troughs. 
Even if 
no substantial differences were encountered in terms of conduction velocity, the calcium transients displayed generally higher values in the viscoelastic case.

\subsection{{Viscoelastic vs. hyperelastic effects under passive inflation and active contraction}}
Much more evident differences can be observed in terms of the true stress $\sigma_f = \bF\fo \cdot (\bsigma \bF\fo)$ when plotted against the local stretch 
in the fibre direction, $\lambda_f=\sqrt{I_{4,f}}$. 
Such a comparison has been conducted in  \cite{gultekin16} for idealised geometries, and it was specifically designed to study hysteresis effects due to viscous contributions to orthotropic passive stress. 
For the inflation tests we will  restrict to $\beta = 1$\,ms and $\delta = 22.6$\,N$/$cm$^2\cdot$ms. 
These values, considered in \cite{karlsen17} (and using units of [s] and [Pa\,s], respectively), ensure that the viscoelastic component is large enough to have a visible effect, but does not completely overwhelm the dynamics of the tissue. 
Here we consider the left ventricular domain used in Section~\ref{sec:scroll} and 
proceed to analyse a stress-stretch response on two points near the basal surface on the endocardium and 
epicardium, and portrayed in Figure \ref{fig:visc-ventr}(a). The mechanical parameters were taken differently from those in Table~\ref{table:params}; here we 
focus on the patient-specific constants estimated from healthy myocardial tissue at 8 mmHg end-diastolic pressure using chamber pressure-volume and strain data taken 
in vivo \cite{gao}. The modified values for this particular test are 
$a=0.02096$\,N$/$cm$^2$, $b=3.243$, $a_f=0.30634$\,N$/$cm$^2$, $b_f=3.4595$, $a_s=0.07334$\,N$/$cm$^2$, $b_s=1.5473$, $a_{fs}=0.03646$\,N$/$cm$^2$, $b_{fs}=3.39$ and we set $\zeta_{\text{stab}} = 10$. 
In the simulation we impose a sinusoidal endocardial pressure of maximal amplitude 0.1\,N$/$cm$^2$ (approximately 8 mmHg) and run a set of transient simulations 
over the interval from 0 to 300\,ms. 
This configuration constitutes an inflation and deflation process where the majority of the fibres are acting in traction, whereas 
sheetlets work under a compression regime. 
Plots (b,c) in Figure \ref{fig:visc-ventr} illustrate the stress-stretch response (in terms of the true stress). The behaviour on the 
epicardial point shows an exponential stiffening and is quite similar to what was observed in  \cite{gultekin16}, as for both stress measures in the viscoelastic case there is evidence of hysteresis effects (that are, by definition, not present in the hyperelastic case). 
Slight deviations from the reference results in \cite{gultekin16} maybe related to the fact that we are using a full electromechanical model, a different viscoelastic contribution, and different material parameters. On the endocardial point we observe even more marked differences between the two cases, 
probably since we do not expect symmetry in the 
motion patterns for a non-ellipsoidal geometry.  Other qualitative differences in the motion patterns include a more marked 
wall thickening, and an overall lower pressure (also more evenly spread throughout the endocardium, showing a smoother profile than the one produced in the hyperelastic case). Pressure on the epicardium was higher in the viscoelastic case.

\begin{figure*}
\begin{center}
\subfigure[]{\raisebox{2mm}{\includegraphics[width=0.175\textwidth]{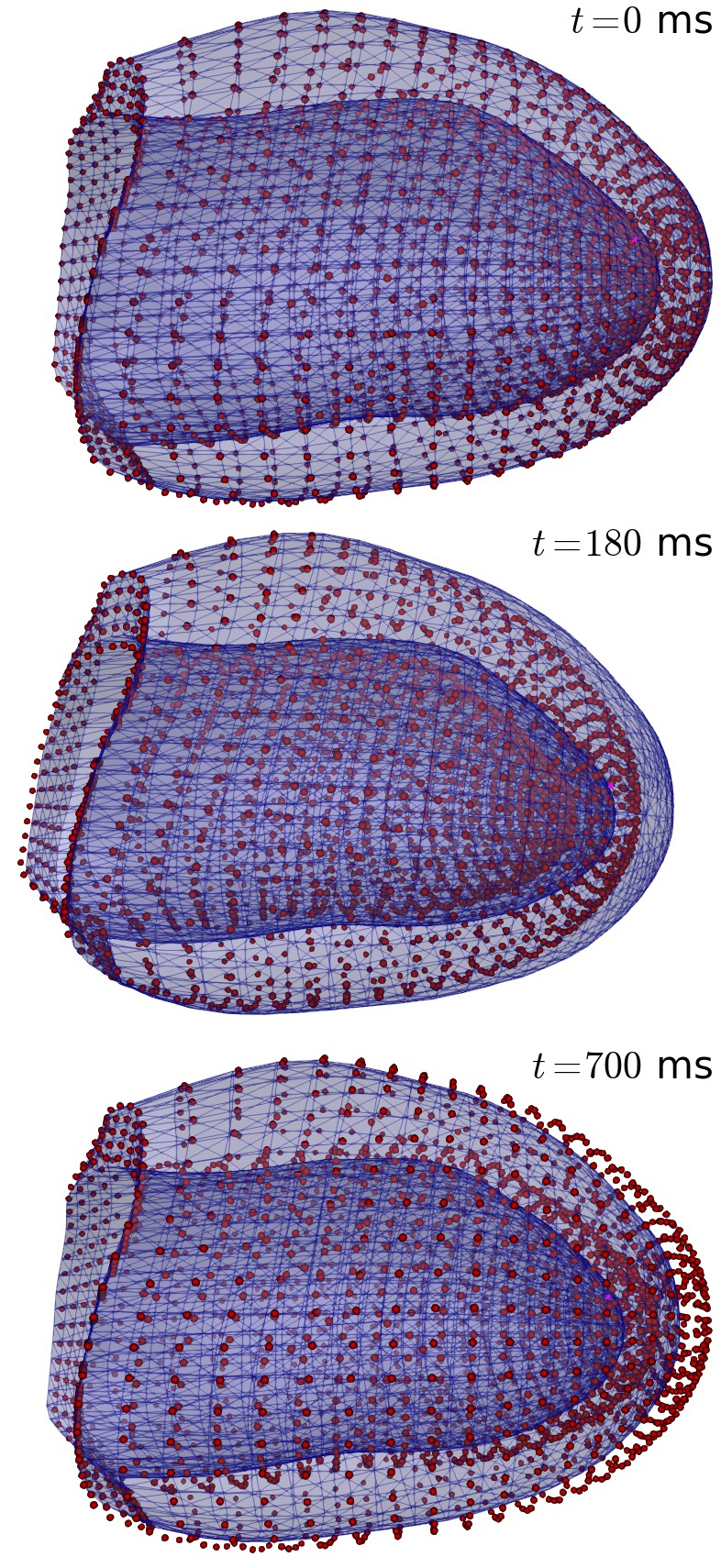}}}
\subfigure[]{\includegraphics[width=0.4\textwidth]{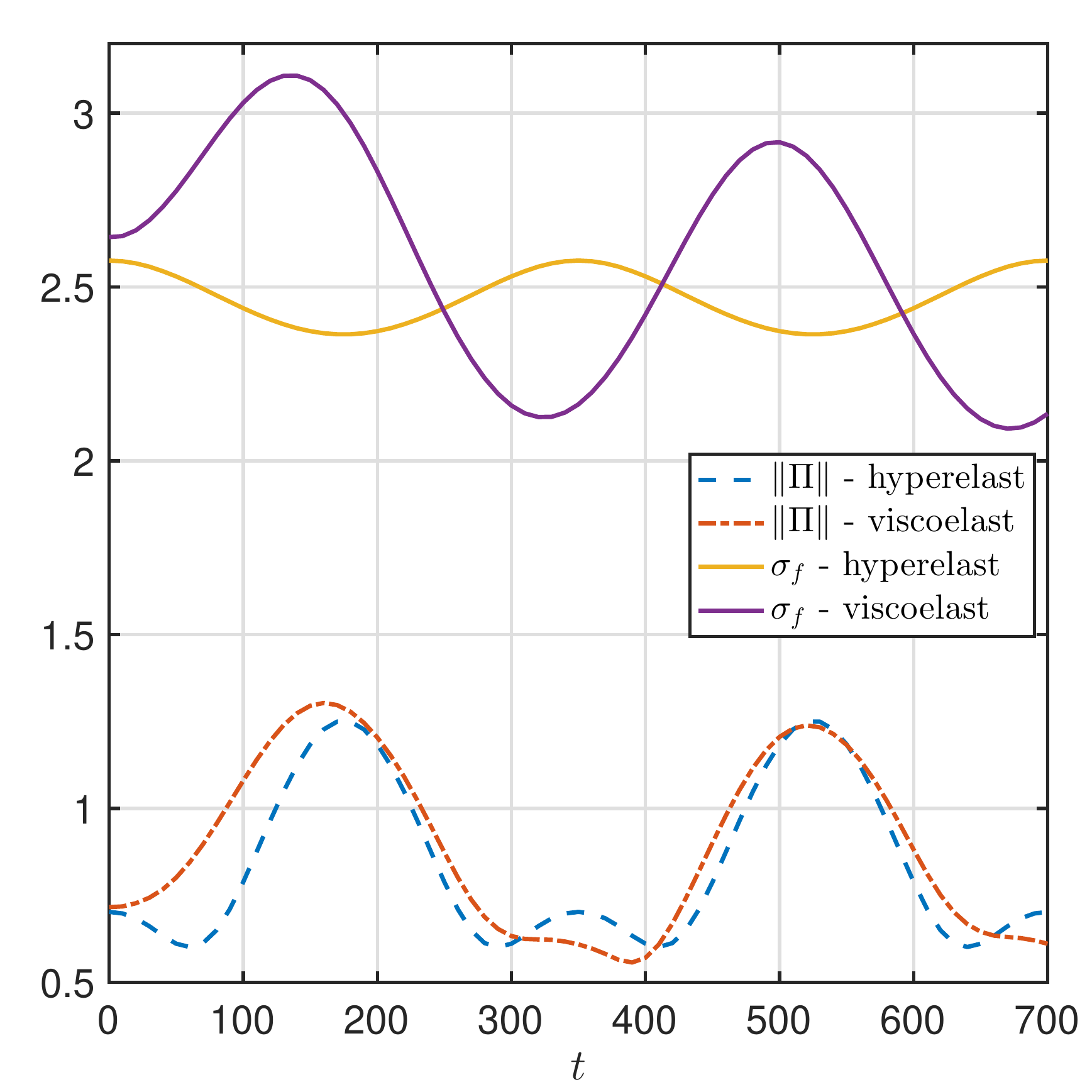}}
\subfigure[]{\includegraphics[width=0.4\textwidth]{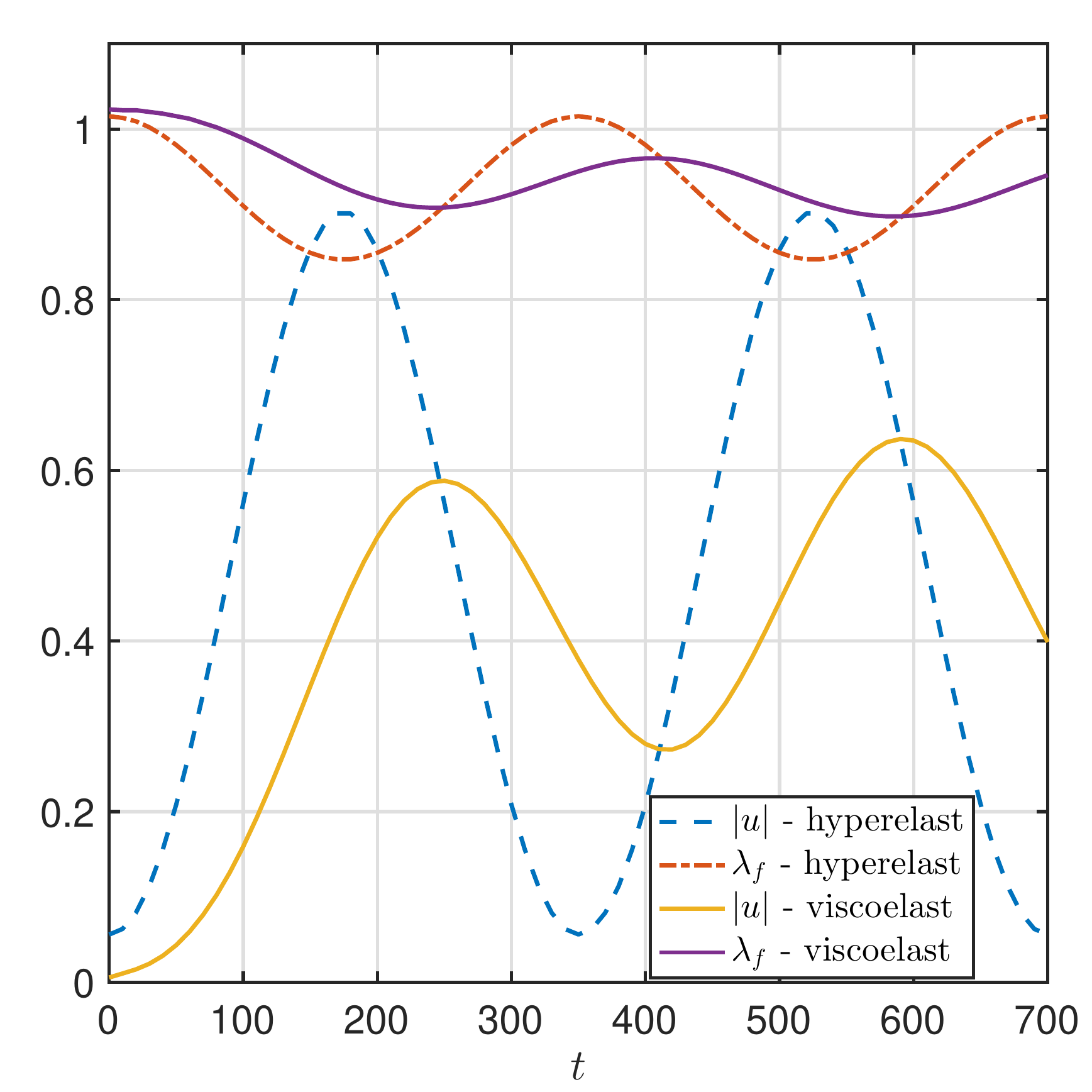}}
\end{center}
\vspace{-2mm}
\caption{{Comparison between hyperelastic and viscoelastic effects under contraction, using the active strain approach and 
$\beta = 5$\,ms and $\delta = 80.25$\,N$/$cm$^2\cdot$ms. Panel (a) has side views of the deformed domain for viscoelastic (hollow blue) and hyperelastic (dark red dots) at 
three different times, and panels (b-c) show transient of mechanical outputs extracted from a point on the lower epicardial surface.}} 
\label{fig:revised-ex}
\end{figure*}

{As a final test we analyse the differences between the viscoelastic and hyperelastic case 
under active contraction. We employ the ventricular geometry again, imposing Robin boundary 
conditions on the epicardial surface and zero normal displacements on the basal cut, and use the active 
strain approach. We use $\zeta_{\text{stab}} = 10$, and more pronounced viscoelastic effects encoded in 
$\beta = 5$\,ms and $\delta = 80.25$\,N$/$cm$^2\cdot$ms. The active contraction of the ventricle is initiated by 
an ectopic beat and impose a sinusoidal endocardial pressure but with larger amplitude 0.25\,N$/$cm$^2$ and simulate 
the process for approximately two cycles (from 0 to 700\,ms). Figure~\ref{fig:revised-ex} reveals a higher 
asynchrony in the tissue deformation and stresses between the hyperelastic and viscoelastic case. This is observed visibly 
from the motion of the ventricle in panel (a), but also from the transients extracted on an epicardial 
point near the apex and measuring 
true stress in the fibre direction, the local fibre stretch, the 
displacement magnitude. This effect was milder in the Frobenius norm of the Kirchhoff stress tensor.
}

\section{Concluding remarks}\label{sec:concl}
We have introduced a model  for 
 the active contraction of cardiac tissue. We focused on incorporating the 
 mechanoelectric feedback through stress-assisted diffusion,  
accounting for a porous-media-type nonlinear diffusivity, and including 
 inertial terms in the equations of motion. 
The three-field equations of motion of a viscoelastic orthotropic material are 
coupled with a four-variable minimal model for human ventricular action potential 
using both active strain and active stress approaches. We have also proposed a new 
stabilised mixed-primal numerical scheme written, in particular, in terms of the Kirchhoff stress. 
The nontrivial effects of both viscoelasticity and stress-assisted diffusion in our model suggest that they may play an important role in governing cardiac function and its response to external stimuli.

{An important remark is that the active strain approach seem to be much more robust than the active stress 
formulation, at least in the present context. In order to obtain comparable results with the active 
stress approach, we had to spend quite a lot of effort finding the correct scaling in \eqref{eq:active-stress}. Once this is achieved, 
we observed that there is no 
substantial difference in the output quantities. 
This is why our numerical tests have focused a little bit more on the active stress, since it is somewhat more challenging than the other case}. 
 
Further additions will be mostly focused on multiscale microstructural coupling, which will provide a more physiological justification of the model in terms of complex phenomena involved in mechano-electrical interactions. One example would be to include poroelastic effects representing perfusion of the myocardial tissue. Developing a thermodynamically consistent description of stress-assisted diffusion is also a pending task, in which electromechanical coupling with the surrounding torso and organs would represent another level of interaction. Such formulation under electromechanical coupling (and including nonlinear and stress-assisted diffusion), will require state-of-the-art tools of multiscale homogenisation \cite{cyron16} as well as dedicated multiscale numerical methods  \cite{gandhi16}.

More confident now in obtaining accurate and reliable numerical solutions, our forthcoming contributions will target an exhaustive computational analysis of restitution curves and realistic activation patterns, e.g. accounting for Purkinje fibres and cellular heterogeneity, with the purpose of characterising spatiotemporal alternans patterns in the presence of multiple mechano-electric feedback effects. Practical applications of the present study rely on the antitachycardia pacing protocols, as well as on the (still today not completely understood) effects of mechanical loads, including cardiac massage, tissue damage, and remodelling at different scales during atrial flutter \cite{mase08}. Estimates of energy dissipation and heat production would be further investigated, widening the validity of these models to non-equilibrium thermodynamical systems. The simulation of mechanically-induced ectopic activity, as well as prediction of the risk of sudden cardiac death are also part of our long-term goals.
For sake of model validation we are also interested in marrying our results to experimental observations obtained through elastography (using for instance the novel approach advanced in \cite{capilnasiu19}).  

\begin{acknowledgements}
We express our sincere thanks to Pablo Lamata (KCL) for providing the base-line surface LV meshes, as well as Figure~\ref{fig:segmentation}. 
In addition, fruitful discussions with Xiyao Li, Renee Miller, and David Nordsletten are gratefully acknowledged. 
{We also thank the efforts of two anonymous referees whose many constructive comments lead to a number of improvements to the manuscript.} 
This work has been partially supported by the Oxford Master in Mathematical Modelling and Scientific Computing, by the S19 Visiting Professorship Programme at University Campus Bio-Medico, and by the EPSRC through the research grant EP/R00207X/1.  
\end{acknowledgements}

\bigskip 
\noindent\textbf{Conflict of interest} The authors declare that they have no conflict of interest.

\appendix 
{\section{Details on the minimal ionic model}\label{sec:app-minimal}
The ionic currents 
consist of three general  terms, phenomenologically constructed (without particularisation to the ionic
species that generate them) 
\begin{align*}
	 g(v,\vec{r}) &= g^{\text{fi}}(v,\vec{r}) + g^{\text{si}}(v,\vec{r}) + g^{\text{so}}(v,\vec{r}),
\end{align*}
where the adimensional \emph{fast inward}, \emph{slow inward} and 
\emph{slow outward} currents are respectively given by 
\begin{align*}
\chi \, g^{\text{fi}}(v,\vec{r}) &= -r_1 \cH(v-\theta_1)(v-\theta_1)(v_v-v)/\tau_{fi},  \\
\chi \, g^{\text{si}}(v,\vec{r}) &= -\cH(v-\theta_2) r_2 r_3/\tau_{si}, \\
\chi \, g^{\text{so}}(v,\vec{r}) &= \frac{(v-v_0)(1-\cH(v-\theta_2))}{\tau_o}+\cH(v-\theta_2)/\tau_{so},
\end{align*}
and the kinetics of the gating variables $\vec{r}$ are given by
\begin{align*}
&\vec{m}(v,\vec{r})=\\
&\begin{pmatrix}
(1-\cH(v-\theta_1))(r_{1,\inf} - r_1)/\tau_1^- - \cH(v-\theta_1)r_1/\tau_1^+\\
(1-\cH(v-\theta_2))(r_{2,\inf} - r_2)/\tau_2^--\cH(v-\theta_2)r_2/\tau_2^+\\
((1+\tanh(k_3(v-v_3)))/2-r_3)/\tau_3 \end{pmatrix}.
\end{align*} 
Here $\cH$ is the Heaviside step function, and the time constants
and steady-state values are defined as:  
\begin{align*}
\tau_1^-&=(1-\cH(v-\theta_1^-))\tau_{1,1}^-+\cH(v-\theta_1^-)\tau_{1,2}^-,\\ 
\tau_2^-&=\tau_{2,1}^-+(\tau_{2,2}^--\tau_{2,1}^-)(1+\tanh(k_2^-(v-v_2^-)))/2,\\
\tau_{so}&=\tau_{so,1}+(\tau_{so,2}-\tau_{so,1})(1+\tanh(k_{so}(v-v_{so})))/2,\\
\tau_3&=((1-\cH(v-\theta_2))\tau_{3,1}+\cH(v-\theta_2)\tau_{3,2}, \\
\tau_o & = ((1-\cH(v-\theta_0))\tau_{o,1}+\cH(v-\theta_0)\tau_{o,2},\\
r_{1,\inf}&= \begin{cases}
1, & v<\theta_1^-\\ 
0,& u\geq \theta_1^-\end{cases},\\ 
r_{2,\inf} & = ((1-\cH(v-\theta_0))(1-v/\tau_{2,\infty})+\cH(v-\theta_0) r_{2,\infty}^*.
\end{align*}
}

{\section{Mixed-primal fully discrete finite element scheme}\label{sec:app-mixed}
We restrict the presentation to the active stress 
formulation using a smoothed model for active tension (stressing that 
the case of active strain and other boundary conditions follows a similar treatment). 
Let us denote by $\cT_h$ a regular partition of $\Omega$ into  
simplicial elements $K$ (pair-wise disjoint triangles in 2D or tetrahedra in 3D) of maximum diameter 
$h_K$, and define the mesh size as $h:=\max\{h_K: K\in \cT_h\}$. Let us also denote 
by $\mathcal{E}_h$ the set of interior facets of the mesh, and by $\jump{\cdot}_e$ the 
jump of a quantity across a given facet $e\in\mathcal{E}_h$. 
The specific finite element method we use here is based on solving  
 the discrete weak
form of the hyperelasticity equations  using piecewise constant approximations of the symmetric 
Kirchhoff stress tensor, 
piecewise linear approximation of displacements, and piecewise constant 
approximation of (solid) pressure. The transmembrane potential in the
electrophysiology equations is discretised with Lagrange finite elements
(piecewise linear and continuous functions), and the remaining ionic quantities are 
approximated by piecewise constant functions. More precisely, we use the finite dimensional 
spaces $\mathbb{H}_h\!\subset\! \mathbb{L}^2_{\text{sym}}(\Omega)$, 
$\mathbf{V}_h\subset \mathbf{H}^1(\Omega)$, $W_h\subset \mathrm{H}^1(\Omega)$, 
$Q_h\subset\mathrm{L}^2(\Omega)$, $Z_h\subset\mathrm{L}^2(\Omega)^3$ defined (for the 
case of a generic-order approximation $l\geq 0$) as follows: 
\begin{align*}
\mathbb{H}_h&:=\{ \btau_h\in  \mathbb{L}^2_{\text{sym}}(\Omega) : \btau_h|_K\in \mathbb{P}_l(K)^{d\times d},\forall K\in \cT_h\},\\
\mathbf{V}_h&:=\{\bv_h\in \mathbf{H}^1(\Omega):\bv_h|_K\in\mathbb{P}_{l+1}(K)^d,\forall  K\in \cT_h,\\ 
&\qquad\qquad  \bv_h\cdot\bn=0 \text{ on } \partial\Omega_D\},\\
Q_h&:=  \{ q_h\in \mathrm{L}^2(\Omega) : q_h|_K\in \mathbb{P}_{l}(K),\forall  K\in \cT_h\},\\
W_h&:=  \{ w_h\in \mathrm{H}^1(\Omega) : w_h|_K\in \mathbb{P}_{l+1}(K),\forall  K\in \cT_h\},\\
Z_h&:=  \{ \varphi_h\in \mathrm{L}^2(\Omega) : \varphi_h|_K\in \mathbb{P}_{l}(K),\forall  K\in \cT_h\},
\end{align*}
where $\mathbb{P}_l(K)$ denotes the space of polynomial functions of degree $s\leq l$ 
defined locally on the element $K$. 
Assuming zero body loads, and applying a backward differentiation formula (BDF) for the 
time integration, we end up with the following fully-discrete 
nonlinear electromechanical problem, starting from the discrete initial data 
$v^0_h,n^{0}_h,T_{a,h}^0$. 
For each $n=0,1,\ldots,n_{\max}$: find \\
$(\bPi_h^{n+1},\bu_h^{n+1},p_h^{n+1})$ and $(v_h^{n+1},\vec{r}_h^{n+1},T_{a,h}^{n+1})$ 
such that 
\begin{widetext}
\begin{align}
\int_\Omega [\bPi_h^{n+1} - \cG(\bu_h^{n+1}) + p_h^{n+1}J(\bu_h^{n+1})\bI] : \btau_h & = 0 & \forall \btau_h \in \mathbb{H}_h,\nonumber\\
\int_\Omega \frac{\bu_h^{n+1}-2\bu_h^n+\bu_h^{n-1}}{\Delta t^2}\cdot \bv_h +\int_\Omega \bPi_h^{n+1}\! :\! \nabla\bv_h \bF^{\tt -t}(\bu_h^{n+1}) - \int_{\partial\Omega_N} \!\!p_N \bF^{\tt -t}(\bu_h^{n+1}) \bn\cdot\bv_h 
& = 0 &  \forall \bv_h \in \mathbf{V}_h, \nonumber\\
\int_\Omega [J(\bu_h^{n+1})  - 1 ] q_h + \sum_{e\in \mathcal{E}_h} \int_e \zeta_{\text{stab}}h_e [\![p_h^{n+1}]\!]_e\jump{q_h}_e& = 0 &  \forall q_h\in Q_h, \label{eq:FE}\\
\int_\Omega \frac{v_h^{n+1}-v_h^{n}}{\Delta t} w_h + \int_\Omega  \bD(v_h^{n+1},\bPi_h^{n+1})\, \nabla v_h^{n+1} \cdot \nabla w_h 
  - \int_\Omega \biggl[g(v_h^{n},\vec{r}_h^{n}) + I_{\rm ext}\biggr] w_h &= 0 & \forall w_h \in  W_h, \nonumber\\
\int_\Omega \frac{\vec{r}_h^{n+1}-\vec{r}_h^{n}}{\Delta t} \cdot\vec{s}_h  - \int_\Omega 
\vec{m}(v_h^{n},\vec{r}_h^{n})\cdot \vec{s}_h  & = 0 & \forall \vec{s}_h \in  Z_h^3,\nonumber\\
\int_\Omega   \frac{T_{a,h}^{n+1}-T_{a,h}^{n}}{\Delta t} \varphi_h  
+ \alpha_1D_0\int_\Omega \nabla T_{a,h}^{n+1} \cdot \nabla \varphi_h 
- \int_\Omega 
\ell(T_{a,h}^{n+1},\vec{r}_h^{n}) \varphi_h  & = 0 & \forall \varphi_h\in  W_h,\nonumber
\end{align}
\end{widetext}
where $\zeta_{\text{stab}}$ is a {positive pressure} stabilisation parameter required to 
enforce solvability of the discrete problem. This is the 
tetrahedral counterpart of the finite element method for quadrilateral meshes studied in \cite{chavan07} 
and recently  used in the context of cardiac electromechanics in \cite{ruiz18}. 
Notice that the boundary condition \eqref{eq:dirichletBC} is incorporated as an essential condition 
on the displacement space, whereas the traction boundary condition \eqref{eq:pressureBC} on the remainder of the boundary 
$\partial\Omega_N$ appears naturally as the last term in the second equation of \eqref{eq:FE}. 
}

{\section{Linearisation of the mechanical problem} 
The coupling between activated mechanics and the electrophysiology
solvers  will be performed using a segregated fixed-point scheme. At each time step, the nonlinear algebraic sub-system for the 
mechanics defined by the first three equations in \eqref{eq:FE} is linearised, adopting 
the following form (where the time indices have been dropped for the sake of notation).} 

{Starting from the initial guess $(\bPi_h^{k=0},\bu_h^{k=0},p_h^{k=0}) =  (\bPi_h^n,\bu_h^n,p_h^n)$, 
for $k=0,1,\ldots$ find stress, displacement and pressure increments 
$(\delta\bPi_h^{k+1},\delta\bu_h^{k+1},\delta p_h^{k+1})$ such that 
\begin{widetext}
\begin{align}
\int_\Omega \bigl[\delta\bPi_h^{k+1} + \frac{\partial\cG_h^k}{\partial\bF_h^k}\nabla \delta\bu_h^{k+1} +(\delta p_h^{n+1}J_h^{k} + p_h^kJ_h^k\bF_h^{k,\tt -t}:\nabla\delta\bu_h^{k+1})\bI\bigr] : \btau_h & = \int_\Omega \mathcalbb{R}^k_{\bPi}:\btau_h & \forall \btau_h \in \mathbb{H}_h,\nonumber\\
\int_\Omega \frac{\delta\bu_h^{k+1}}{\Delta t^2}\cdot \bv_h 
+\int_\Omega \bigl[\delta\bPi_h^{k+1}\bF_h^{k,\tt -t}-\bPi_h^{k}\bF_h^{k,\tt -t}(\nabla \delta\bu_h^{k+1})^{\tt t} \bF_h^{k,\tt -t}\bigr] : \nabla\bv_h \qquad &
\label{eq:tangent}
 \\
+ \int_{\partial\Omega_N} \!\!p_N \bF_h^{k,\tt -t}(\nabla \delta\bu_h^{k+1})^{\tt t} \bF_h^{k,\tt -t}\bn \cdot\bv_h 
& = \int_\Omega \boldsymbol{\cR}^k_{\bu}\cdot \bv_h &  \forall \bv_h \in \mathbf{V}_h, \nonumber\\
\int_\Omega \bigl(J_h^k\bF_h^{k,\tt -t}:\nabla\delta\bu_h^{k+1}\bigr) q_h + \sum_{e\in \mathcal{E}_h} \int_e \zeta_{\text{stab}}h_e [\![ \delta p_h^{n+1}]\!]_e\jump{q_h}_e& =  \int_\Omega \cR^k_{p}q_h&  \forall q_h\in Q_h, \nonumber
\end{align}
\end{widetext}
and then update $\bPi_h^{k+1} = \bPi_h^k + \delta\bPi_h^{k+1}$,  $\bu_h^{k+1} = \bu_h^k + \delta\bu_h^{k+1}$,   $p_h^{k+1} = p_h^k + \delta p_h^{k+1}$. Here $\mathcalbb{R}^k_{\bPi}, \boldsymbol{\cR}^k_{\bu},\cR^k_{p}$ are tensor, vector, and scalar residuals associated with the Newton-Raphson linearisation at the previous step $k$, and $\bF_h^k=\bI + \nabla\bu_h^k$, $J_h^{k} = \det\bF_h^k$. 
Next we introduce the following linear maps (related to the G\^{a}teaux derivatives of the solution operator) 
\begin{align*}
\mathcalbb{L}_1^k: \mathbb{H}_h \to \mathbb{H}_h, &\btau \mapsto \mathcalbb{L}_1^k (\btau) = \frac{\partial\cG_h^k}{\partial\bF_h^k} \btau \!+\! (p_h^kJ_h^k\bF_h^{k,\tt -t}\!\!:\btau) \bI,\\
\mathcalbb{L}_2^k: \mathbb{H}_h \to \mathbb{H}_h, & \btau \mapsto \mathcalbb{L}_2^k (\btau) =- \bPi_h^{k}\bF_h^{k,\tt -t} \btau^{\tt t} \bF_h^{k,\tt -t},
\end{align*} 
as well as the bilinear forms and linear functionals 
\begin{align*} 
A_1(\bPi,\btau) &= \int_\Omega \bPi : \btau,   \\  
A_2 (\bu,\bv) &=  \frac{1}{\Delta t^2} \!\! \int_\Omega \bu\cdot\bv + \int_\Omega\mathcalbb{L}_2^k(\nabla\bu):\nabla \bv,  \\  
A_3(p,q) & =  \sum_{e\in \mathcal{E}_h}\!\! \int_e \zeta_{\text{stab}}h_e [\![ p]\!]_e\!\jump{q}_e, \\
B_1(\bu,\btau) &= \int_\Omega \mathcalbb{L}_1^k(\nabla\bu) : \btau,\\   
\widetilde{B}_1(\bPi,\bv) & = \int_\Omega \bPi\,\bF_h^{k,\tt -t} :\nabla \bv, \\   
B_2(p,\btau) & =  \int_\Omega J_h^k\, p\tr(\btau), \\   
\widetilde{B}_3(\bu,q) &=  \int_\Omega (J_h^k\bF_h^{k,\tt -t}\!:\!\nabla\bu)\,q,  \ F_1(\btau) = \int_\Omega \mathcalbb{R}^k_{\bPi}:\btau,\\ F_2(\bv) &= \int_\Omega \boldsymbol{\cR}^k_{\bu}\cdot \bv,\quad F_3(q) = \int_\Omega\cR^k_{p}q.
\end{align*}
Then, dropping all iteration indices and making abuse of notation, 
the tangent problem \eqref{eq:tangent} (now also restricted to pure displacement 
boundary conditions) can be recast as the mixed variational form
{\small\begin{alignat}{3}
A_1(\bPi_h,\btau_h)  \ +\ & B_1(\bu_h,\btau_h) & +\ B_2(p_h,\btau_h) & = F_1(\btau_h) & \forall \btau_h \in \mathbb{H}_h,\nonumber\\
\widetilde{B}_1(\bPi_h,\bv_h)  \ +\ &  A_2 (\bu_h,\bv_h) & & = F_2(\bv_h) &  \forall \bv_h \in \mathbf{V}_h, \label{eq:linearMixed}\\
 & \widetilde{B}_3(\bu_h,q_h) & +\ A_3(p_h,q_h) & = F_3(q_h) &\forall q_h \in Q_h.\nonumber
\end{alignat}}
The stability and convergence analysis of \eqref{eq:linearMixed} will be studied in a forthcoming contribution.}
\end{document}